\documentclass{amsart}
\usepackage[english]{babel}
\usepackage{amsfonts}
\usepackage{amsmath,amsthm,amssymb,latexsym}
\usepackage{mathtools}
\usepackage{authblk}
\usepackage{tikz}
\usepackage{stmaryrd}
\usepackage{mathrsfs,upgreek}
\usepackage{eucal}
\usepackage{enumerate}
\usepackage{changes}
\usepackage{amsaddr}
\usepackage{enumitem}
\usepackage{hyperref}
\usepackage{comment}
\usepackage{ulem}
\usepackage{tikz-cd}
\usepackage{mleftright}
\usepackage[capitalise]{cleveref}

\newtheorem{theorem}{Theorem}[section]
\newtheorem{lemma}[theorem]{Lemma}
\newtheorem{proposition}[theorem]{Proposition}
\newtheorem{corollary}[theorem]{Corollary}

\theoremstyle{definition}
\newtheorem{definition}[theorem]{Definition}

\newtheorem{example}[theorem]{Example}

\newtheorem{notation}[theorem]{Notation}

\newtheorem{remark}[theorem]{Remark}

\numberwithin{equation}{section}

\newcommand{\op}[1]{\textrm{\upshape #1}}
\newcommand{\cov}{\prec}

\newcommand{\join}{\vee}

\newcommand{\meet}{\wedge}

\newcommand{\la}{\langle}
\newcommand{\ra}{\rangle}
\newcommand{\alg}[1]{{\textbf{\upshape #1}}}
\newcommand{\vv}[1]{\mathcal {#1}}

\renewcommand{\d}{\delta}

\newcommand{\g}{\gamma}
\newcommand{\e}{\varepsilon}
\renewcommand{\th}{\theta}

\newcommand{\cg}{\vartheta}

\newcommand{\sse}{\subseteq}

\newcommand{\app}{\approx}
\newcommand{\HH}{{\mathbf H}}  
\newcommand{\II}{{\mathbf I}} 
\newcommand{\SU}{{\mathbf S}} 
\newcommand{\PP}{{\mathbf P}}   
\newcommand{\VV}{{\mathcal V}}   

\newcommand{\ib}{\item[$\bullet$]}

\newcommand{\Con}[1]{\operatorname{Con}(\alg #1)}

\newcommand{\con}{\operatorname{Con}}
\newcommand{\Ker}{\operatorname{Ker}}
\newcommand{\lcm}{\operatorname{lcm}}

\newcommand{\ulm}{\mathsf{u\ell M}}
\newcommand{\MVM}{\mathsf{MVM}}

\newcommand{\N}{\mathbb{N}}
\newcommand{\Z}{\mathbb{Z}}
\newcommand{\R}{\mathbb{R}}
\newcommand{\Q}{\mathbb{Q}}

\newcommand{\black}{\color{black}}

\begin{document}
	
	\title{Varieties of MV-monoids and positive MV-algebras}
	
	\author{Marco Abbadini}
	
	\address{Marco Abbadini,
		School of Computer Science,
		University of Birmingham,
		University Rd W,
		Birmingham B15 2TT,
		UK}
	\email{\tt m.abbadini@bham.ac.uk}
	
	\author{Paolo Aglian\`o}
	
	\address{Paolo Aglian\`o,
		DIISM,
		Universit\`a di Siena,
		Via Roma 56,
		53100 Siena,
		Italy}
	\email{\tt agliano@live.com}
	
	\author{Stefano Fioravanti}
	
	\address{Stefano Fioravanti,
		Department of Algebra, Faculty of Mathematics and Physics, Charles University, Sokolovská 83 18675 Praha 8,
Czech Republic}
	\email{\tt stefano.fioravanti66@gmail.com}
	
	\keywords{Lattice-ordered monoids,
Lattice-ordered groups,
Categorical equivalence,
MV-algebras,
positive MV-algebras,
Varieties of algebras,
Equational axiomatization}
    \subjclass[2020]{Primary: 06F05. Secondary: 06D35, 03C05, 08B15, 08B26, 08C15
}
	
\maketitle
\begin{abstract}
    MV-monoids are algebras $\la A,\join,\meet, \oplus,\odot, 0,1\ra$ where $\langle A, \lor, \land, 0, 1\rangle$ is a bounded distributive lattice, both $\langle A, \oplus, 0 \rangle$ and $\langle A, \odot, 1\rangle$ are commutative monoids, and some further connecting axioms are satisfied.
    Every MV-algebra in the signature $\{\oplus,\neg,0\}$ is term equivalent to an algebra that has an MV-monoid as a reduct, by defining, as standard, $1\coloneqq \neg 0$, $x \odot y \coloneqq \neg(\neg x \oplus\neg y)$, $x \join y \coloneqq (x \odot \neg y) \oplus y$ and $x \meet y \coloneqq \neg(\neg x \join \neg y)$. Particular examples of  MV-monoids are positive MV-algebras, i.e., the $\{\join, \meet, \oplus, \odot, 0, 1\}$-subreducts of MV-algebras.
    Positive MV-algebras form a peculiar quasivariety in the sense that, albeit having a logical motivation (being the quasivariety of subreducts of MV-algebras), it is not the equivalent quasivariety semantics of any logic.

    In this paper, we study the lattices of subvarieties of MV-monoids and of positive MV-algebras.
    In particular, we characterize and axiomatize all almost minimal varieties of MV-monoids, we characterize the finite subdirectly irreducible positive MV-algebras, and we characterize and axiomatize all varieties of positive MV-algebras.
\end{abstract}

\section{Introduction}
MV-algebras were introduced in \cite{chang1958algebraic} to serve as algebraic semantics for the many-valued {\L}ukasiewicz propositional logic.  In modern accounts, they are defined as algebras of the form $\langle A, \oplus, \neg, 0\rangle$ satisfying a certain list of equational axioms.
The prime example of an MV-algebra is the \emph{standard MV-algebra}, i.e., the real interval $[0,1]$ endowed with the operations:
\begin{align}\label{eq:standard}
    x\oplus y\coloneqq \min\{x+y,1\}, \quad \neg x\coloneqq 1-x, \quad \text{and the constant $0$}.
\end{align}
The standard MV-algebra $[0,1]$ generates the variety of MV-algebras.
The reader is referred to \cite{CignoliDOttavianoEtAl2000} for the basic theory of MV-algebras and to \cite{Mundici2011} for more advanced topics.

The study of MV-algebras is made easier by their tight relationship with Abelian $\ell$-groups.
An \emph{Abelian lattice-ordered group} (\emph{Abelian $\ell$-group}, for short) is an Abelian group equipped with a translation-invariant lattice-order; see \cite{BigardKeimelEtAl1977} for the theory of lattice-ordered groups.
A \emph{unital Abelian $\ell$-group} is an Abelian $\ell$-group $\alg{G}$ equipped with a \emph{strong order unit}, i.e., an element $1 \in G$ such that for every $y \in G$ there is $n \in \N$ such that $y \leq 1 + \dots + 1$ ($n$ times).
The prime example of a unital Abelian $\ell$-group is $\R$, equipped with its additive structure, its usual order, and the element $1$ as a strong order unit.
The theory of unital Abelian $\ell$-groups is far simpler than the one of MV-algebras.
Thus, the categorical equivalence between the category of MV-algebras and the category of unital Abelian $\ell$-groups, proved by D.~Mundici \cite[Theorem~3.9]{Mun1986} (see also \cite[Section~2]{CignoliDOttavianoEtAl2000}), has become a fundamental tool in the theory of MV-algebras.
In particular, every MV-algebra $\alg{A}$ is the unit interval of some unital Abelian $\ell$-group $\alg{G}$.
This means that $\alg{A}$ is isomorphic to the algebra with base set $\{x \in G \mid 0 \leq x \leq 1\}$ and operations
\begin{equation} \label{eq:op-in-Gamma}
    x \oplus y \coloneqq (x + y) \land 1, \quad \lnot x \coloneqq 1 - x, \quad \text{ and the constant $0$}.
\end{equation}
One may notice the analogy between \eqref{eq:standard} and \eqref{eq:op-in-Gamma}.
In fact, under Mundici's categorical equivalence, the MV-algebra $[0,1]$ corresponds to the unital Abelian $\ell$-group $\R$.

Further important examples of MV-algebras are Boolean algebras: every such algebra can be equipped with the structure of an MV-algebra by defining $\oplus$ as the binary join $\lor$. In fact, Boolean algebras are a subvariety of MV-algebras, defined by the equation $x \oplus x = x$.
MV-algebras are then a ``many-valued'' generalization of Boolean algebras, where the role of the 2-element Boolean algebra $\{0,1\}$ is taken by the standard MV-algebra $[0,1]$: for example, the variety of Boolean algebras is generated by $\{0,1\}$, while the variety of MV-algebras is generated by $[0,1]$.

Much of the theory of Boolean algebras can be smoothly generalized to bounded distributive lattices.
The Boolean terms that are definable in the language of bounded distributive lattices are precisely those that are \emph{order-preserving} in each argument.
For instance, negation is not order-preserving.
So, bounded distributive lattices can be seen as the \emph{order-preserving} (or \emph{positive}) fragment of Boolean algebras.

In analogy to the relationship between Boolean algebras and bounded distributive lattices, recent years have marked the beginning of the investigation of the negation-free version of MV-algebras.
To understand what an appropriate definition of a ``negation-free version of MV-algebras'' should be, we refer the reader to a fundamental result \cite[Theorem~3.5]{CintulaKroupa13}, which states that the order-preserving terms of MV-algebras are precisely those in the language $\{\lor, \land, \oplus, \odot, 0, 1\}$, where $\lor$, $\land$, $\odot$ and $1$ can be term-defined in the language of MV-algebras as $1\coloneqq \neg 0$, $x \odot y \coloneqq \neg(\neg x \oplus\neg y)$, $x \join y \coloneqq (x \odot \neg y) \oplus y$ and $x \meet y \coloneqq \neg(\neg x \join \neg y)$.
This suggests that the ``negation-free versions of MV-algebras'' should be algebras in the signature $\{\lor, \land, \oplus, \odot, 0, 1\}$.

There are two reasonable non-equivalent definitions of the negation-free versions of MV-algebras: \emph{positive MV-algebras} and \emph{MV-monoidal algebras} (called \emph{MV-monoids} in this paper).
This is similar to what happens for Abelian groups, where there are two possible candidates for ``the inverse-free version of Abelian groups'': \emph{cancellative commutative monoids} and \emph{commutative monoids}.
While the subreducts of Abelian groups are precisely the cancellative commutative monoids, a certain portion of the theory of Abelian groups holds for commutative monoids. The latter ones, in contrast to the cancellative ones, have the advantage of being defined by equations, making them a variety of algebras.

\emph{Positive MV-algebras} (introduced in \cite{Abbadinietal2022}) are defined as the $\{\lor, \land, \oplus, \odot, 0, 1\}$-subreducts of MV-algebras.
They are the quasivariety generated by $[0,1]$, considered as a positive MV-algebra, and they are precisely the unit intervals of unital cancellative commutative $\ell$-monoids (\cite[Lemma 3.8(2)]{Abbadinietal2022}).
On the other hand, \emph{MV-monoids}, introduced in \cite{Abbadini2021b} via a certain list of equations (which make them a variety of algebras), are precisely the unit intervals of unital commutative $\ell$-monoids (this is the main result of \cite{Abbadini2021b}).
In a nutshell, positive MV-algebras are precisely the \emph{cancellative} MV-monoids, where the cancellation property in this context is expressed by the following quasi-equation:
\[
(x \oplus z \app y \oplus z\ \text{and}\ x \odot z \app y \odot z) \ \Longrightarrow \ x \app y.
\]
Note the similarity with the cancellation property for commutative monoids:
\[
x + z \app y + z\  \Longrightarrow \ x \app y.
\]
MV-monoids and positive MV-algebras are peculiar classes: even if they are not quasivarieties of logic, at least in the sense of \cite{BarbourRaftery2003}, their link with logic is obvious, especially for positive MV-algebras, which constitute the quasivariety of subreducts of MV-algebras, and many standard techniques of universal algebra, commonly applied in algebraic logic, are applicable in this setting as well.

In this paper, we investigate varieties of MV-monoids and of positive MV-algebras.
Regarding the former, our main result is a characterization of the almost minimal varieties of MV-monoids (\cref{t:char-almost-minimal}).
For the latter, our main result is a characterization of the varieties of positive MV-algebras (\cref{t:VarofpMVs}): these are precisely the varieties generated by finitely many reducts of finite nontrivial MV-chains. We also prove that such reducts are precisely the subdirectly irreducible finite positive MV-algebras (\cref{thm:fin-is-sub}).
Furthermore, we provide an equational axiomatization for each such variety of positive MV-algebras and each almost minimal variety of MV-monoids (\cref{thm:axioms,t:axiomVarC2}).

The structure of the paper is as follows: in \cref{sec: prel} we introduce the main definitions and state the main results from the literature that will be used in the sequel.
In \cref{sec: si MVM} we study subdirectly irreducible MV-monoids. In \cref{sec: bottom} we study the bottom part of the lattice of subvarieties of MV-monoids, characterizing the almost minimal varieties, and in \cref{s:almostmin} we investigate some of the varieties above the almost minimal varieties. We then turn our attention to positive MV-algebras investigating those that are subdirectly irreducible in \cref{sec: si pMV}. In \cref{sec: pMV varieties} we describe the varieties of positive MV-algebras and we show that all of them are finitely generated. In fact, the varieties of positive MV-algebras are precisely the varieties generated by finitely many reducts of finite nontrivial MV-chains.
Finally, in \cref{sec:axiomatization}, we present an axiomatization for each variety of positive MV-algebras and each almost minimal variety of MV-monoids.

\section{Preliminaries: lattice-ordered monoids, MV-monoids and positive MV-algebras}\label{sec: prel}

In this section, we recall the definitions of the algebraic structures of interest in this paper.
For elementary concepts in general algebra (such as lattices, algebras, varieties, etc.), our textbook reference is \cite{BurrisSanka}, while our textbook reference for MV-algebras is \cite{Mundici2011}.

We start by defining commutative $\ell$-monoids and unital commutative $\ell$-monoids.

\begin{definition} A \emph{commutative $\ell$-monoid} is an algebra
$\alg{M} = \la M, \join,\meet,+,0\ra$ with the following properties:
\begin{enumerate}
    \ib $\la M, \join,\meet\ra$ is a distributive lattice;
    \ib $\la M,+,0\ra$ is a commutative monoid;
    \ib $+$ distributes over $\join$ and $\meet$.
\end{enumerate}

\end{definition}
We warn the reader that some authors would rather use the terminology ``commutative \emph{distributive} $\ell$-monoid''; in this paper we prefer the shorter name.

\begin{example}\label{ex:clm}\hfill
    \begin{enumerate}
        \item
        The algebra $\langle \R, \max, \min, +, 0\rangle$ is a commutative $\ell$-monoid, as well as any of its subalgebras, such as $\Q$, $\Z$, $2\Z$, $\N$, $\{0,-1, -2, -3, \dots\}$.

        \item
        For every preordered set $X$, the set of order-preserving functions from $X$ to $\R$ with pointwise defined operations is a commutative $\ell$-monoid.

        \item ($\alg{C}_2^{\Delta*}$)\label{ex:Delta}
        An example of a commutative $\ell$-monoid is the two-element chain $\{0 < \e\}$, with $\lor = \max$, $\land = \min$, $+ = \lor$, and $0 = 0$.
        We write $\alg{C}_2^{\Delta*}$ for this commutative $\ell$-monoid.
        It is not cancellative.

        \item ($\alg{C}_2^{\nabla*}$) \label{ex:nabla}
        Dually to \eqref{ex:Delta}, an example of a commutative $\ell$-monoid is the two-element chain $\{\d < 0\}$, with $\lor = \max$, $\land = \min$, $+ = \land$, and $0 = 0$.
        We write $\alg{C}_2^{\nabla*}$ for this commutative $\ell$-monoid.
        It is not cancellative.

        \item ($\alg{C}_n^{\Delta*}$ and $\alg{C}_n^{\nabla*}$) \label{ex:Cn*}
        \Cref{ex:Delta,ex:nabla} are instances of a slightly more general construction.
        For every $n \in \N$, the following defines a congruence $\sim_n$ on the additive $\ell$-monoid $\N$:
        \[
        x \sim_n y \iff x = y \text{ or } (x \geq n -1 \text{ and } y \geq n - 1).
        \]
        We let $\alg{C}_n^{\Delta*}$ denote the quotient $\N / {\sim_n}$, and we let $\alg{C}_n^{\nabla*}$ denote the order-dual of $\alg{C}_n^{\Delta*}$, i.e., the algebra obtained from $\alg{C}_n^{\Delta*}$ by swapping $\lor$ and $\land$.
        This gives a source of examples of $\ell$-monoids whose monoidal reduct is not cancellative.

        \item ($\alg{LM}_n^{\Delta*}$ and $\alg{LM}_n^{\nabla*}$)
        \label{ex:LMn*}
        \Cref{ex:Delta,ex:nabla} are instances of another slightly more general construction.
        For every $n \in \N$, we let $\alg{LM}_n^{\Delta*}$ denote the commutative $\ell$-monoid defined on the $n$-element chain, where $+$ is defined as $\lor$, and $0$ as the bottom element.
        Analogously, we let $\alg{LM}_n^{\nabla*}$ denote the commutative $\ell$-monoid defined on the $n$-element chain, where $+$ is defined as $\land$, and $0$ as the top element.
        This gives another source of examples of $\ell$-monoids whose monoidal reduct is not cancellative.
    \end{enumerate}
\end{example}

\begin{remark}
    Since $\ell$-monoids are defined by equations, they form a variety and hence they are closed under products, subalgebras, and homomorphic images.
    This allows one to obtain several examples.
\end{remark}

\begin{definition} \label{d:unital-commutative-ell-monoid}
    A \emph{unital commutative $\ell$-monoid} is an algebra $\la M, \join,\meet,+,1,0,-1\ra$ with the following properties:
    \begin{enumerate}
        \ib $\la M, \join,\meet,+,0\ra$ is a commutative $\ell$-monoid;
        \ib $-1+1 =0$;
        \ib $-1 \le 0 \le 1$;
        \ib for all $x \in M$ there is $n \in \mathbb N$ such that
        $$
        \underbrace{(-1)+ \dots + (-1)}_{n \text{ times}} \le x \le \underbrace{1+ \dots +1}_{n\text{ times}}.
        $$
    \end{enumerate}
\end{definition}
In a unital commutative $\ell$-monoid we will use the standard abbreviations
\begin{align*}
    n & \coloneqq 1 + \dots + 1\ \ \text{$n$ times},\\
    -n& \coloneqq (-1) + \dots + (-1)\ \ \text{$n$ times},\\
    x-n& \coloneqq x +(-n),
\end{align*}
omitting the usual superscript for the interpretation in the given algebra. As usual, for all $n \in \N$, we denote by $nx$ the sum of $n$ copies of $x$ and by $-nx$ the sum of $n$ copies of $-x$ if $x$ is invertible.

We let $\mathsf{u\ell M}$ denote the class of unital commutative $\ell$-monoids.

\begin{example}\hfill
    \begin{enumerate}
        \item \label{ex:ulms}
        The algebra $\langle \R, \max, \min, +, 1, 0, -1 \rangle$ is a unital commutative $\ell$-monoid, and so are all its subalgebras, such as $\Q$ and $\Z$.
    	An example of a subalgebra of $\R$ that is not an additive subgroup of $\R$ is, for any irrational element $s$ in $\R$, the algebra $\{a + bs \mid a \in \Z, b \in \N \}$.

        \item
        For every preordered set $X$, the set of bounded order-preserving functions from $X$ to $\R$ with pointwise defined operations is a unital commutative $\ell$-monoid.

        \item (Lexicographic product $\overrightarrow{\times}$)
        \label{i:lexico}
        For every totally ordered unital commutative $\ell$-monoid $\alg M$ and every commutative $\ell$-monoid $\alg L$,
        we have a unital commutative $\ell$-monoid $\alg{M} \overrightarrow{\times} \alg L$ defined as follows: the underlying set is $M \times L$, the order is lexicographic, i.e.,
        \[
            \big((n,x) \leq (m,y)\big) \Longleftrightarrow \big((n < m) \text{ or } (n = m \text{ and } x \leq y)\big),
        \]
        the operation $+$ is defined component-wise, i.e., $(n,x) + (m,y) = (n + m, x + y)$, the positive unit is $(1,0)$, the negative unit is $(-1,0)$, and $(0,0)$ is the zero.

        For example, when $\alg{M} = \Z$, the algebra $\Z \overrightarrow{\times} \alg L$ stacks $\Z$-many copies of $\alg{L}$ one on top of the other.

        When combined with the noncancellative examples of commutative $\ell$-monoids in \cref{ex:clm}, this lexicographic product produces a source of examples of unital commutative $\ell$-monoids whose monoidal reduct is not cancellative: for instance, $\Z \overrightarrow{\times}\alg{C}_n^{\Delta*}$ and $\Z \overrightarrow{\times}\alg{C}_n^{\nabla*}$, $\Z \overrightarrow{\times} \alg{LM}_n^{\nabla*}$, and $\Z \overrightarrow{\times} \alg{LM}_n^{\Delta*}$.
    \end{enumerate}
\end{example}

\begin{remark}
	All the axioms of unital commutative $\ell$-monoids are equations, except for the last one.
	However, we notice that the last one is preserved by subalgebras, homomorphic images, and finite products.
	Thus, unital commutative $\ell$-monoids are closed under these class operators.
	This is not the case for arbitrary products: for example, $\R^\N$ does not satisfy the last axiom.
\end{remark}

\begin{remark}
    In \cite{Santschi2024} the author studied varieties of commutative and distributive $\ell$-monoids, with particular emphasis on the idempotent case, defined by the equality $x + x = x$ (beware: in \cite{Santschi2024}, $\ell$-monoids were presented with the \emph{multiplicative} notation, as opposed to the \emph{additive} notation used in the present paper). However, the intersection of this setting with ours is quite small: the only idempotent unital $\ell$-monoid $\alg M$ is the trivial one: indeed, if $1 = 2$ then $-1 = 0 = 1$ and everything collapses.
\end{remark}

Given a unital commutative $\ell$-monoid $\alg{M}$, its ``unit interval''
\begin{equation} \label{eq:Gamma}
    \Gamma(\alg{M}) \coloneqq \{x \in M \mid 0 \le x \le 1\}
\end{equation}
can be turned into an algebra in the signature $\{\join,\meet,\oplus,\odot,0,1\}$: set $x \oplus y\coloneqq (x +y) \meet 1$, $x \odot y \coloneqq (x +y-1) \join 0$, and define $\lor$, $\land$, $0$ and $1$ by restriction.

To capture axiomatically the algebras of the form $\Gamma(\alg{M})$ for a unital commutative $\ell$-monoid $\alg{M}$, the first author introduced in \cite{Abbadini2021b} the notion of an \emph{MV-monoidal algebras}, or \emph{MV-monoid}, for short.
\begin{definition}
    An \emph{MV-monoid} is an algebra $\la A,\join,\meet,\oplus,\odot,0,1\ra$ where
    \begin{enumerate}
        \ib $\la A,\join,\meet,0,1\ra$ is a bounded distributive lattice;
        \ib $\la A,\oplus,0\ra$ and $\la A,\odot,1\ra$ are commutative monoids;
        \ib $\oplus$ and $\odot$ distribute over $\join$ and $\meet$;
        \ib for every $x,y,z \in A$,
        \begin{align*}
            &(x \oplus y) \odot ((x \odot y)\oplus z) = (x \odot (y \oplus z)) \oplus (y \odot z);\\
            &(x \odot y) \oplus ((x \oplus y)\odot z) = (x \oplus (y \odot z)) \odot (y \oplus z);\\
            &(x \odot y) \oplus z = ((x \oplus y) \odot ((x \odot y) \oplus z)) \join z;\\
            &(x \oplus y) \odot z = ((x \odot y) \oplus ((x \oplus y) \odot z)) \meet z.
        \end{align*}
    \end{enumerate}
\end{definition}
When there is no danger of confusion, we use the additive notation $nx \coloneqq x \oplus \dots \oplus x$ ($n$ times) and the multiplicative notation $x^n \coloneqq x \odot \dots \odot x$ ($n$ times).
For an explanation of the choice of the axioms of MV-monoids, we direct the reader to \cite[pp.\ 45 ff]{Abbadini2021b}.
Clearly, MV-monoids form a variety of algebras, which we denote by $\mathsf{MVM}$.

For every unital commutative $\ell$-monoid $\alg{M}$, its unit interval
\[
\Gamma(\alg M) = \la \Gamma(\alg{M}),\join,\meet,\oplus,\odot,0,1\ra
\]
defined in \eqref{eq:Gamma} is an MV-monoid.
The assignment $\alg{M} \mapsto \Gamma(\alg{M})$ can be extended to morphisms to define a functor
\[
\Gamma \colon \mathsf{u\ell M} \rightarrow \mathsf{MVM}
\]
from the category $\mathsf{u\ell M}$ of unital commutative $\ell$-monoids to the category $\mathsf{MVM}$ of MV-monoids: a homomorphism $f\colon \alg M \rightarrow \alg N$ between unital commutative $\ell$-monoids is mapped to its restriction $\Gamma(f) \colon \Gamma(\alg{M}) \to \Gamma(\alg{N})$.

The following is the crucial result connecting unital commutative $\ell$-monoids and MV-monoids.

\begin{theorem}[\cite{Abbadini2021b}]
    The functor $\Gamma \colon \mathsf{u\ell M} \to \mathsf{MVM}$ is an equivalence of categories.
\end{theorem}

In particular, for every MV-monoid $\alg{A}$ there is a (unique up to isomorphism) unital commutative $\ell$-monoid $\alg{M}$ such that $\Gamma(\alg{M}) \cong \alg{A}$.

In this paper, an important role is taken by those unital commutative $\ell$-monoids
that are \emph{cancellative}.
We recall that a commutative monoid $\alg M$ is \emph{cancellative} if, for all $x,y,z \in M$, $x+z=y+z$ implies $x=y$.
The MV-monoids of the form $\Gamma(\alg{M})$ for some \emph{cancellative} unital commutative $\ell$-monoid $\alg{M}$ are precisely the \emph{positive MV-algebras}.

\begin{definition}[\cite{Abbadinietal2022}]
    A \emph{positive MV-algebra} is a $\{\join,\meet,\oplus,\odot,0,1\}$-subreduct of an MV-algebra.
\end{definition}

As proved in \cite{Abbadinietal2022}, the class $\mathsf{MV}^+$ of positive MV-algebras is a proper subquasivariety of $\mathsf{MVM}$, axiomatized relatively to $\mathsf{MVM}$ by
$$
(x \oplus z \app y \oplus z\ \text{and}\ x \odot z \app y \odot z) \ \Longrightarrow \ x \app y.
$$
Note the similarity with the cancellation property in the language of commutative $\ell$-monoids: $x + z \approx y + z \Rightarrow x \approx y$.
For an MV-algebra $\alg A$, we denote by $\alg A^+$ its MV-monoid reduct.

\begin{theorem}[\cite{Abbadinietal2022}] \label{t:restricts}
    The functor $\Gamma \colon \mathsf{u\ell M} \to \mathsf{MVM}$ restricts to an equivalence between cancellative unital commutative $\ell$-monoids and positive MV-algebras.
\end{theorem}

In particular, for every positive MV-algebra $\alg{A}$ there is a (unique up to isomorphism) cancellative unital commutative $\ell$-monoid $\alg{M}$ such that $\Gamma(\alg{M}) \cong \alg{A}$.

Moreover, $\Gamma \colon \mathsf{u\ell M} \to \mathsf{MVM}$ further restricts to Mundici's equivalence between unital Abelian $\ell$-groups and MV-algebras (see \cite[Theorem~4.79]{Abbadini2021} or \cite[Theorem~A.4]{Abbadini2021b}).

Next, we give examples of MV-monoids and positive MV-algebras. For each $n \in \mathbb N \setminus \{0\}$, we set $\frac{1}{n}\mathbb Z \coloneqq \{\frac{k}{n} \mid k \in \Z\}$.
This set can be equipped in a natural way with the structure of a unital commutative $\ell$-monoid, in which the sum is the ordinary addition in $\mathbb Q$ and the order is the natural one.
We recall that $\alg {\L}_n$ denotes the MV-algebra whose universe is $\{0, \frac{1}{n}, \dots, \frac{n-1}{n}, 1\}$ and whose operations are defined in the usual way.  It is easy to check that $\Gamma(\frac{1}{n}\mathbb Z)= \alg  {\L}^+_n$.

As another example of a positive MV-algebra we recall the \emph{Chang algebra} $\alg C$, introduced by C.~C.~Chang \cite[p.~474]{chang1958algebraic}.
The \emph{Chang algebra} is an MV-algebra that we represent in the following way: the universe of $\alg C$ is
$$
\{n\e \mid n \in \N\} \cup \{\d^n \mid n \in \N\}
$$
where $\d^n$ stands for the usual $1-n\varepsilon$. Note that, here, $\e$ is simply the name for the smallest nonzero element. The choice for this particular name is motivated by the fact that its behavior with respect to $\oplus$ and $\odot$ can be loosely interpreted as the one of an infinitesimal. 
We use the shorthands $0$ for $0\e$ and $1$ for $\d^0$. 
The two commutative operations $\oplus, \odot$ and the unary operation $\lnot$ are defined as follows:
\begin{align*}
    n\e \oplus m \e  &\coloneqq (n+m)\e;\\
    \d^n \oplus \d^m  &\coloneqq 1;\\
    n\e \oplus \d^m = \d^m \oplus n\e  &\coloneqq
    \begin{cases}
        \d^{m-n} & \hbox{if $n<m$,} \\
        1 & \hbox{otherwise;}
    \end{cases}\\
    n\e \odot m \e  &\coloneqq 0;\\
    \d^n \odot  \d^m &\coloneqq \d^{n+m};\\
    n\e \odot \d^m = \d^m \odot n\e &\coloneqq
    \begin{cases}
        (n-m)\e & \hbox{if $n >m$,} \\
        0 & \hbox{otherwise;}
    \end{cases}\\
    \lnot n\e &\coloneqq \d^n;\\
    \lnot \d^n &\coloneqq n \e.
\end{align*}
It is easy to see that the induced lattice is given by
\begin{align*}
    & n \e \le  m \e\ \text{if and only if}\  n \le m\\
    & \d^n \le  \d^m\ \text{if and only if}\  m \le n\\
    & n\e < \d^m\ \text{for all } n,m \in \N.
\end{align*}
This order has $0$ as the bottom element and $1$ as the top element.
Using the lexicographic notation of \cref{i:lexico}, (the $\{\lor,\land,\oplus, \odot,  0, 1\}$-reduct of) $\alg C$ is isomorphic to $\Gamma(\Z \overrightarrow{\times} \Z)$.

\begin{definition}[$\alg{C}^\Delta$ and $\alg{C}^\nabla$]
    We denote by $\alg{C}^\Delta$ and $\alg{C}^\nabla$ the subalgebras of $\alg C^+$ whose universes are $C^\Delta \coloneqq \{0,\e,2\e, \dots,1\}$ and $C^\nabla \coloneqq \{0,\dots,\d^3,\d^2,\d,1\}$, respectively.
\end{definition}
We note that $\alg{C}^\Delta$ and $\alg{C}^\nabla$ are isomorphic to $\Gamma(\Z \overrightarrow{\times} \N)$ and $\Gamma(\Z \overrightarrow{\times} {-\N})$, respectively, where ${-\N}$ is the subalgebra of $\Z$ consisting of the nonpositive elements.

The algebras $\alg{C}^\Delta$ and $\alg{C}^\nabla$ are two proper positive MV-subalgebras of $\alg C$ which are not reducts of MV-algebras.

The construction of a quasi-inverse $\Xi \colon \MVM \to \ulm$ of $\Gamma$\label{def:quasiinverse} relies on the notion of a \emph{good $\Z$-sequence} in an MV-monoid, which is motivated by the following fact of independent interest.

\begin{lemma}[{\cite[Lemma 4.74]{Abbadini2021}}] \label{abbadinith}
    Let $\alg M$ be a unital commutative $\ell$-monoid, and let $x, y \in M$. If for all $n\in \mathbb Z$ we have
    $((x-n) \join 0) \meet 1 = ((y-n) \join 0)\meet 1$,
    then $x = y$.
\end{lemma}
This shows that the elements of $\alg{M}$ can be identified with certain functions from $\Z$ to $\Gamma(\alg{M})$.
To illustrate exactly which functions, we recall the notion of a \emph{good pair}.

\begin{definition}[{\cite[Definition~4.30]{Abbadini2021} or \cite[Definition~5.1]{Abbadini2021b}}]\label{d:goodpair}
    A \emph{good pair} in an MV-monoid $\alg{A}$ is a pair $(x_0, x_1) \in A^2$ such that $x_0 \oplus x_1 = x_0$ and $x_0 \odot x_1 = x_1$.
\end{definition}

\begin{definition}[{\cite[Definition~ 4.31]{Abbadini2021}}]
    A \emph{good $\Z$-sequence} in an MV-monoid $\alg{A}$ is a function $\mathbf{x} \colon \Z \to A$ with the following properties.
    \begin{enumerate}
        \item For all $n \in \Z$, $(\mathbf{x}(n), \mathbf{x}(n + 1))$ is a good pair in $\alg{A}$;
        \item there is $n \in \N$ such that, for every $m \in \N$ with $m \geq n$, $\mathbf{x}(-m) = 1$ and $\mathbf{x}(m) = 0$.
    \end{enumerate}
\end{definition}
For each MV-monoid $\alg{A}$, $\Xi(\alg{A})$ is defined as a certain unital commutative $\ell$-monoid whose underlying set is the set of good $\Z$-sequences in $\alg{A}$.
For the definition of the operations of unital commutative $\ell$-monoid on $\Xi(\alg{A})$ and the fact that $\Xi$ is a quasi-inverse of $\Gamma$, we direct the reader to \cite[Chapter~4]{Abbadini2021}.

For later usage, we mention that \cref{abbadinith} admits the following analogue, in which equality is replaced by $\leq$.

\begin{lemma} \label{p:good}
    Let $\mathbf{M}$ be a unital commutative $\ell$-monoid, and let $x, y \in M$.
    If for all $n \in \Z$ we have $((x-n) \lor 0) \land 1 \leq ((y-n) \lor 0) \land 1$, then $x \leq y$.
\end{lemma}

\begin{proof}
    Suppose that for all $n \in \Z$ we have $((x-n) \lor 0) \land 1 \leq ((y-n) \lor 0) \land 1$.
    Then for all $n \in \Z$ we have
    \begin{align*}
        (((x \land y)-n) \lor 0) \land 1 & = (((x -n) \land (y-n)) \lor 0) \land 1 \\
        & = (((x-n) \lor 0) \land 1) \land (((y-n) \lor 0) \land 1) \\
        & = ((x-n) \lor 0) \land 1.
    \end{align*}
    By {\cref{abbadinith}}, it follows that $x \land y = x$, i.e., $x \leq y$.
\end{proof}
We collect here some properties of $\Gamma$ and $\Xi$.

\begin{proposition} \label{preserveandreflect}
    The functors $\Gamma$ and $\Xi$ preserve and reflect injectivity, surjectivity, and bijectivity of morphisms.
\end{proposition}
\begin{proof}
    First, since $\Xi$ is the quasi-inverse of $\Gamma$, it is enough to prove the statement for $\Gamma$. Next, we observe that bijective morphisms are isomorphisms and they are always preserved by equivalences. The fact that $\Gamma$ preserves and reflects injectivity is proved in \cite[Proposition 3.7]{Abbadini2021b}.

    We prove that $\Gamma$ preserves surjectivity.	Let $f \colon \alg{M} \to \alg{N}$ be a surjective homomorphism of unital commutative $\ell$-monoids, and let $y\in \Gamma(\alg{N})$.
    By surjectivity of $f$, there is $x \in M$ such that $f(x) = y$.
    Let $x' \coloneqq (x \lor 0) \land 1$. Then $x' \in \Gamma(\alg{M})$ and
    \begin{align*}
        f(x') = f((x \lor 0) \land 1) = (f(x) \lor 0) \land 1 = (y \lor 0) \land 1 = y,
    \end{align*}
    since $0 \leq y \leq 1$. Thus, $\Gamma(f)$ is surjective.

    Finally, we prove that $\Gamma$ reflects surjectivity.
    Let $f \colon \alg{M} \to \alg{N}$ be a homomorphism of unital commutative $\ell$-monoids, and suppose that $\Gamma(f)$ is surjective.
    Let $y \in N$.
    By \cite[Proposition 4.68]{Abbadini2021} there exist $k \in \Z$, $n \in \N$, and $y_0, \dots, y_n \in \Gamma(\alg{N})$ such that $y = k + y_0 + \dots + y_n$.
    Since $\Gamma(f)$ is surjective, for every $i \in \{0, \dots, n\}$ there is $x_i \in \Gamma(\alg{M})$ such that $y_i = \Gamma(f)(x_i) = f(x_i)$.
    Then
    \begin{align*}
        f(k + x_0 + \dots + x_n) = f(k) + f(x_0) + \dots + f(x_n) = k + y_0 + \dots + y_n = y.
    \end{align*}
    Therefore, $f$ is surjective.
\end{proof}

 We let $\con(\alg{A})$ denote the congruence lattice of an algebra $\alg{A}$. Using the previous Proposition we can directly obtain an isomorphism between congruence lattices of unital commutative $\ell$-monoid and of MV-monoids.

\begin{corollary} \label{c:isomorphic-lattices-of-congruences}
    Let $\alg M$ be a unital commutative $\ell$-monoid. The congruence lattices of $\alg M$ and $\Gamma(\alg M)$ are isomorphic.
\end{corollary}

\begin{proof}
    For a congruence $\theta$ on an algebra $\alg{A}$, we denote by $\pi_\theta \colon \alg{A} \to \alg{A}/ \theta$ the quotient map induced by $\theta$.
    For a homomorphism $f \colon \alg{A} \to \alg{B}$ between similar algebras, we denote by $\Ker(f)$ the kernel congruence induced by $f$.
    We claim that the function
    \begin{align*}
        \phi \colon \con(\alg M) &\longrightarrow \con(\Gamma(\alg M))\\
        \theta & \longmapsto \text{Ker}(\Gamma(\pi_{\theta}))
    \end{align*}
    is a lattice isomorphism. Using the Homomorphism Theorem \cite[Theorem 6.12]{BurrisSanka} and the fact that $\Gamma$ preserves isomorphisms (cf.\ \cref{preserveandreflect}), it is easy to see that, for every surjective homomorphism $f \colon \alg{M} \twoheadrightarrow \alg{N}$, $\phi(\Ker(f)) = \Ker(\Gamma(f))$.
    \[
    \begin{tikzcd}[column sep = 5em]
        \alg M \arrow[two heads]{r}{\pi_{\Ker(f)}} \arrow[two heads,swap]{rd}{f}& \alg M/\Ker(f) \arrow[hook,two heads]{d}{}  & \Gamma(\alg M) \arrow[two heads]{r}{\Gamma(\pi_{\Ker(f)})} \arrow[swap, two heads]{rd}{\Gamma(f)}& \Gamma(\alg M/\Ker(f)) \arrow[hook,two heads]{d}{}\\
        & \alg N & & \Gamma(\alg N)
    \end{tikzcd}
    \]

    We prove that $\phi$ preserves the lattice order. Let $\alpha \leq \beta$ be two congruences of $\alg{M}$. By a standard consequence of the Homomorphism theorem (see e.g.\ \cite[Chapter~II.7]{BurrisSanka}), there exists a surjective homomorphism $f \colon \alg M/\alpha \rightarrow \alg M/\beta$ such that the diagram on the left-hand side below commutes \eqref{dia:1}. Applying the functor $\Gamma$ to the commutative diagram on the left-hand side of \eqref{dia:1}, we obtain the commutative diagram on the right-hand side.
    \begin{equation}\label{dia:1}
        \begin{tikzcd}[column sep = 5em]
        \alg M \arrow[two heads]{r}{\pi_{\alpha}} \arrow[two heads,swap]{rd}{\pi_{\beta}}& \alg M/\alpha \arrow[two heads]{d}{f}  & \Gamma(\alg M) \arrow[two heads]{r}{\Gamma(\pi_{\alpha})} \arrow[swap, two heads]{rd}{\Gamma(\pi_{\beta})}& \Gamma(\alg M/\alpha) \arrow[two heads]{d}{\Gamma(f)}\\
        & \alg M/\beta & & \Gamma(\alg M/\beta)
    \end{tikzcd}
    \end{equation}
    By \cref{preserveandreflect}, $\Gamma$ preserves surjectivity and thus $\Gamma(f)$, $\Gamma(\pi_{\alpha})$ and $\Gamma(\pi_{\beta})$ are surjective. It follows that $\phi(\alpha) = \text{Ker}(\Gamma(\pi_{\alpha})) \subseteq \text{Ker}(\Gamma(\pi_{\beta})) = \phi(\beta)$. This proves that $\phi$ is order-preserving.

    We prove that $\phi$ is surjective. Let $\gamma \in \con(\Gamma(\alg{M}))$.
    Since $\Gamma$ is essentially surjective, there are a unital commutative $\ell$-monoid $\alg N$ and an isomorphism $h$ from $\Gamma(\alg M)/\gamma$ to $\Gamma(\alg N)$ (for example, take $\alg{N} = \Xi(\Gamma(\alg M)/\gamma)$.
    Since $\Gamma$ is full, there is a homomorphism $f \colon \alg M \to \alg N$ such that $\Gamma(f) = h \circ \pi_\gamma$.
    \[
        \begin{tikzcd}
            \alg M \arrow{rr}{f} && \alg N \\
            \Gamma(\alg M) \arrow{r}{} \arrow[bend left = 2 em, two heads]{rr}{\Gamma(f)} \arrow[swap, two heads]{r}{\pi_{\gamma}} & \Gamma(\alg M)/\gamma \arrow[hook, two heads, swap]{r}{h}& \Gamma(\alg N)
        \end{tikzcd}
    \]
    By \cref{preserveandreflect}, $f$ is surjective.
    Since $h$ is an isomorphism, we have $\Ker(\Gamma(f)) = \Ker(h \circ \pi_\gamma)= \Ker(\pi_\gamma) =  \gamma$.
    Therefore, by the property observed right after the definition of $\phi$, we have $\phi(\Ker(f)) = \Ker(\Gamma(f)) = \gamma$, and so $\gamma$ belongs to the image of $\phi$.

    We prove that $\phi$ reflects the order.
    Let $\alpha, \beta \in \con(\alg M)$ be such that $ \phi(\alpha) \leq \phi(\beta)$, i.e., $\Ker(\Gamma(\pi_\alpha)) \leq \Ker(\Gamma(\pi_\beta))$, and let us prove $\alpha \leq \beta$.
    By a standard consequence of the Homomorphism theorem \cite[Chapter~II.7, Exercises~6(6)]{BurrisSanka}, from $\phi(\alpha) \leq \phi(\beta)$ we deduce the existence of a surjective homomorphism $h \colon \Gamma(\alg M/\alpha) \to \Gamma(\alg M/\beta)$ that makes the following diagram commute.
    \begin{equation}\label{d:diag1}
        \begin{tikzcd}[column sep = 3.5em]
        \Gamma(\alg M) \arrow[two heads]{r}{\Gamma(\pi_{\alpha})} \arrow[swap, two heads]{rd}{\Gamma(\pi_{\beta})}& \Gamma(\alg M/\alpha) \arrow[two heads]{d}{h}\\
        & \Gamma(\alg M/\beta)
    \end{tikzcd}
    \end{equation}
    Since $\Gamma$ is full, there is a homomorphism $f \colon \alg M/\alpha \to \alg M/\beta$ such that $\Gamma(f) = h$ and $f$ is surjective by \cref{preserveandreflect}.
    Since $\Gamma$ is faithful, and by commutativity of \eqref{d:diag1}, also the following diagram commutes.
    \[
    \begin{tikzcd}[column sep = 3.5em]
        \alg M \arrow[two heads]{r}{\pi_{\alpha}} \arrow[two heads,swap]{rd}{\pi_{\beta}}& \alg M/\alpha \arrow[two heads]{d}{f}\\
        & \alg M/\beta
    \end{tikzcd}
    \]
    Thus, $\alpha \leq \beta$. This concludes the proof that $\phi$ is order-reflecting.

    Since $\phi$ is a surjective order-reflecting order-preserving map between lattices, it is a lattice isomorphism.
\end{proof}

\section{Subdirectly irreducible MV-monoids}\label{sec: si MVM}

In this section, we provide a necessary condition for an MV-monoid to be subdirectly irreducible using the equivalence with unital commutative $\ell$-monoids.
We recall that an algebra is subdirectly irreducible if and only if the identity congruence is completely meet irreducible in the congruence lattice.
Our investigation is grounded on the following fact:
\begin{theorem} \label{t:subirr}
     Every subdirectly irreducible commutative $\ell$-monoid is totally ordered.
\end{theorem}
This result is a corollary of \cite[Lemma 1.4]{Repnitzkii1984}, but already in \cite[Corollary 2]{Merlier1971} the author proved that any commutative lattice-ordered monoid is a subdirect product of totally ordered ones and asserted, in Remark 3 of the same paper, that this was an unpublished result by Fuchs.

Constants do not play any role in congruences, and so in subdirect irreducibility, either. Therefore, \cref{t:subirr} has the following corollary.

\begin{corollary} \label{c:sub-totord}
    Every subdirectly irreducible unital commutative $\ell$-monoid is totally ordered.
\end{corollary}
We will use \cref{t:subirr} to study subdirect irreducibility of MV-monoids, in light of the following lemma.

\begin{lemma}\label{l:transfer of subd irr}
    A unital commutative $\ell$-monoid  $\alg M$ is subdirectly irreducible (as a $\{+,\lor,\land,0,1,-1\}$-algebra) if and only if $\Gamma(\alg M)$ is subdirectly irreducible (as a $\{\lor,\land,\oplus,\odot,0,1\}$-algebra).
\end{lemma}
\begin{proof}
    This follows from the fact that the congruence lattices of $\alg{M}$ and $\Gamma(\alg M)$ are isomorphic (\cref{c:isomorphic-lattices-of-congruences}).
\end{proof}
Let $\alg{M}$ be a unital commutative $\ell$-monoid.
Now we investigate how the property that $\alg{M}$ is totally ordered relates to properties of $\Gamma(\alg{M})$.
Of course, if $\alg{M}$ is totally ordered, then $\Gamma(\alg{M})$ is totally ordered.
However, the converse is false.
For example, consider the subalgebra $\alg{N} \coloneqq \{(a,b) \in \Z^2 \mid a \leq b\}$ of $\Z^2 = \langle \Z, \max, \min, +, -1, 0, 1\rangle ^2$.
The three-element chain $\Gamma(\alg{N}) = \{(0,0), (0, 1), (1, 1)\}$ is totally ordered, but $\alg{N}$ is not.
For $\alg{M}$ to be totally ordered, the property that $\Gamma(\alg{M})$ should satisfy besides being totally ordered is captured in item \eqref{i:comparison-Gamma} of the following Lemma.

\begin{lemma} \label{l:equivalent}
    Let $\alg{M}$ be a unital commutative $\ell$-monoid.
    The following conditions are equivalent:
    \begin{enumerate}			
        \item \label{i:comparison-exists} for every $x \in M$, there is $n \in \Z$ such that $n \leq x \leq n+1$;
        \item \label{i:comparison-all-k} for every $k \in \Z$ and $x \in M$, either $k \leq x$ or $x \leq k$;
        \item \label{i:comparison-exists-k}
        there is $k \in \Z$ such that, for every $x \in M$, either $k \leq x$ or $x \leq k$;
        \item\label{i:comparison-Gamma} for every $x,y \in \Gamma(\alg{M})$, either $x \oplus y = 1$ or $x \odot y = 0$.
    \end{enumerate}
\end{lemma}

\begin{proof}
    \eqref{i:comparison-exists} $\Rightarrow$ \eqref{i:comparison-all-k}. This is straightforward.

    \eqref{i:comparison-all-k} $\Rightarrow$ \eqref{i:comparison-exists}. Let $x \in M$. Set $I \coloneqq \{n \in \Z \mid n \leq x\}$ and $S \coloneqq \{n \in \Z \mid x \leq n\}$.
    It is immediate that $I$ is downward closed and that $S$ is upward closed in the natural order of $\Z$.
    Moreover, $I$ and $S$ are nonempty by the last axiom in the definition of a unital commutative $\ell$-monoid.
    Furthermore, by \eqref{i:comparison-all-k}, $S \cup I = \Z$.
    It follows that there is $n \in S$ such that $n + 1 \in I$.
    Then $n \leq x \leq n +1$.

    \eqref{i:comparison-all-k} $\Rightarrow$ \eqref{i:comparison-exists-k}. This follows from the fact that $\Z$ is nonempty.

    \eqref{i:comparison-exists-k} $\Rightarrow$ \eqref{i:comparison-all-k}. Suppose there is $k \in \Z$ such that, for every $x \in M$, either $k \leq x$ or $x \leq k$.
    For every $n \in \Z$ and $x \in M$, since we have either $k \leq x - n + k$ or $x - n + k \leq k$, we have either $n \leq x$ or $x \leq n$.

    \eqref{i:comparison-all-k} $\Rightarrow$ \eqref{i:comparison-Gamma}.
    Let $x, y \in \Gamma(\alg{M})$.
    By \eqref{i:comparison-all-k} (applied to $k = 1$), either $1 \leq x + y$ or $x + y \leq 1$.
    In the first case, $x \oplus y = 1$. In the second case, $x \odot y = 0$.

    \eqref{i:comparison-Gamma} $\Rightarrow$ \eqref{i:comparison-exists-k}.
    We prove \eqref{i:comparison-exists-k} with $k = 1$.
    Let $x \in M$.
    For $n \in \Z$, let $x_n \coloneqq ((x-n) \lor 0) \land 1$.
    By \cite[Proposition 4.64]{Abbadini2021}, for every $n \in \Z$ we have $x_n \oplus x_{n + 1} = x_n$ and $x_n \odot x_{n + 1} = x_{n + 1}$.
    In particular, $x_0 \oplus x_1 = x_0$ and $x_0 \odot x_1 = x_1$.
    By \eqref{i:comparison-Gamma}, either $x_0 \oplus x_1 = 1$ or $x_0 \odot x_1 = 0$.
    In the first case, we have $x_0 = 1$, which implies that for every $n \in \Z$ with $n \leq 0$ we have $x_n = 1$; by \cref{p:good}, it follows that $x \geq 1$.
    In the second case, we have $x_1 = 0$, which implies that for every $n \in \Z$ with $n \in \mathbb N \setminus\{0\}$ we have $x_n = 0$; by \cref{p:good}, it follows that $x \leq 1$.
\end{proof}
Note that the equivalent conditions of \cref{l:equivalent} fail in the example $\alg{N}$ seen above. We are now ready to provide conditions on $\Gamma(\alg{M})$ that are equivalent to $\alg{M}$ being totally ordered.

\begin{proposition}\label{prop: M tot ord}
    Let $\alg{M}$ be a unital commutative $\ell$-monoid.
    The following conditions are equivalent:
    \begin{enumerate}
        \item \label{i:tot-ord} $\alg{M}$ is totally ordered;
        \item \label{i:two-conditions} $\Gamma(\alg{M})$ is totally ordered, and, for every $x,y \in \Gamma(\alg{M})$, $x \oplus y = 1$ or $x \odot y = 0$.
    \end{enumerate}
\end{proposition}

\begin{proof}
    \eqref{i:tot-ord} $\Rightarrow$ \eqref{i:two-conditions}.
    If $\alg{M}$ is totally ordered, then $\Gamma(\alg{M})$ is totally ordered.
    Moreover, by the implication \eqref{i:comparison-exists} $\Rightarrow$ \eqref{i:comparison-Gamma} in \cref{l:equivalent}, for every $x,y \in \Gamma(\alg{M})$ we have $x \oplus y = 1$ or $x \odot y = 0$.

    \eqref{i:two-conditions} $\Rightarrow$ \eqref{i:tot-ord}.
    For every $n \in \Z$, the map $x \mapsto x + n$ is an order-automorphism of $\alg{M}$.
    Therefore, from the fact that $\Gamma(\alg{M})$ is totally ordered, it follows that for every $n \in \Z$ the interval $\{x \in M \mid n \leq x \leq n +1\}$ is totally ordered.
    By the implication \eqref{i:comparison-Gamma} $\Rightarrow$ \eqref{i:comparison-exists} in \cref{l:equivalent}, for every $x \in M$ there is $n \in \Z$ such that $n \leq x \leq n+1$.
    It follows that $\alg{M}$ is totally ordered.
\end{proof}
At the beginning of the section we observed that there are totally ordered MV-monoids $\alg{A}$ such that $\Xi(\alg{A})$ is not totally ordered. Using \cref{prop: M tot ord} we can produce other examples of this fact by considering MV-monoids with elements $x,y$ such that $x \oplus y\not=1$ and $x \odot y \not= 0$.
For example, this happens in any bounded chain $\alg L$ with $\lvert L \rvert >2$ in which $\oplus$ and $\odot$ coincide with the standard join and meet, respectively; in this case, it suffices to take $x,y \in L \setminus \{0,1\}$. By \cref{prop: M tot ord}, $\Xi(\alg{L})$ is not totally ordered.

\begin{theorem}\label{subdMVM}
    If an MV-monoid $\alg A$ is subdirectly irreducible, then it is nontrivial, totally ordered, and such that, for all $x,y \in A$, $x \oplus y =1$ or $x \odot y =0$.
\end{theorem}
\begin{proof}
        Nontriviality is obvious. Let $\alg A$ be a subdirectly irreducible MV-monoid. Then there is a unital commutative $\ell$-monoid  $\alg M$ such that $\alg{A} \cong \Gamma(\alg M)$. By \cref{l:transfer of subd irr}, $\alg M$ is subdirectly irreducible. By \cref{t:subirr}, $\alg M$ is totally ordered and thus $\Gamma(\alg M)$ and $\alg A$ are totally ordered. By
        \cref{prop: M tot ord}, for all $x,y \in A$ we have $x \oplus y =1$ or $x \odot y =0$.
\end{proof}
The previous Theorem was proved also in \cite[Theorems~B.3 and C.4]{Abbadini2021b} (see also \cite[Theorem~4.42 and Corollary~4.48]{Abbadini2021}) with a different strategy.
Here we included a new proof, which uses the categorical equivalence between unital commutative $\ell$-monoids and MV-monoids \cite{Abbadini2021b}, allowing us to piggyback on the known fact that subdirectly irreducible commutative $\ell$-monoids are totally ordered \cite{Repnitzkii1984}.

We will see later (\cref{l:LMn-notsubirr}) that the converse of \cref{subdMVM} does not hold, in the sense that there is a nontrivial totally ordered MV-monoid $\alg{A}$ satisfying $x \oplus y =1$ or $x \odot y =0$ for all $x,y \in A$ which is not subdirectly irreducible. However, we will prove that if we restrict to finite positive MV-algebras, the conditions in \cref{subdMVM} are also sufficient for subdirect irreducibility.

\section{The almost minimal varieties of MV-monoids}\label{sec: bottom}

Given any variety $\vv V$, we denote by $\Lambda(\vv V)$ the set of all subvarieties of $\vv V$. It is well-known that ordering $\Lambda(\vv V)$ by inclusion makes $\Lambda(\vv V)$ a complete lattice, and so we refer to $\Lambda(\vv V)$ as the {\em lattice of subvarieties} of $\vv V$.  The description of this lattice is always a very important step in understanding the variety.
The main result of this section is a description of the bottom of the lattice $\Lambda(\MVM)$ of all varieties of MV-monoids (Figure~\ref{f:almost-minimal}), which turns out to be no simple task.
First, we can observe that $\alg {\L}_1^+$ is a subalgebra of any nontrivial MV-monoid, and so the variety $\VV(\alg {\L}_1^+)$ generated by $\alg{\L}_1^+$ is the only atom in $\Lambda(\mathsf{MVM})$. The covers of an atom are usually called the \emph{almost minimal} subvarieties. These often play an important role in the description of the lattice of all subvarieties.
In \cref{t:char-almost-minimal} we characterize the almost minimal varieties of MV-monoids: these are precisely the varieties $\VV(\alg{\L}_p^+)$ (for $p$ prime), $\VV(\alg{C}_2^\Delta)$ and $\VV(\alg{C}_2^\nabla)$; definitions to follow.

\begin{figure}[h]
    \begin{center}
        \begin{tikzpicture}[scale=1.3]
            \draw (0,0) -- (0,1);
            \draw (0,1) -- (0.75,2) -- (0,1) -- (1.75,2) -- (0,1) -- (2.75,2);
            \draw (0,1) -- (-.75,2) -- (0,1) -- (-1.75,2);
            \draw[fill] (0,0) circle [radius=0.05];
            \draw[fill] (0,1) circle [radius=0.05];
            \draw[fill] (0.75,2) circle [radius=0.05];
            \draw[fill] (1.75,2) circle [radius=0.05];
            \draw[fill] (2.75,2) circle [radius=0.05];
            \draw[fill] (-.75,2) circle [radius=0.05];
            \draw[fill] (-1.75,2) circle [radius=0.05];
            \node at (3.85,2){\footnotesize $\cdots$};
            \node[right] at (0,1){\ \ \tiny $\VV(\alg{\L}_1^+)$};
            \node[right] at (-.75,2){\tiny $\VV(\alg{C}^\nabla_2)$};
            \node[right] at (-1.75,2){\tiny $\VV(\alg{C}^\Delta_2)$};
            \node[right] at (.75,2){\tiny $\VV(\alg {\L}_2^+)$};
            \node[right] at (1.75,2){\tiny $\VV(\alg {\L}_3^+)$};
            \node[right] at (2.75,2){\tiny $\VV(\alg {\L}_5^+)$};
        \end{tikzpicture}
        \caption{The bottom part of $\Lambda(\mathsf{MVM})$}\label{f:almost-minimal}
    \end{center}
\end{figure}

We recall that, as a consequence of Birkhoff's Subdirect Representation Theorem (\cite[Theorem 8.6]{BurrisSanka}), every variety of algebras is generated by its subdirectly irreducible members.
We show that every $\alg{\L}_n^+$ is subdirectly irreducible. For this, we recall that an algebra $\alg{A}$ is \emph{simple} if its congruence lattice is the 2-element chain, and that any simple algebra is subdirectly irreducible.

\begin{proposition}\label{l:hsimple}
    Every subalgebra of the unital commutative $\ell$-monoid $\R$ is simple.
\end{proposition}

\begin{proof}
    Let $\alg{M}$ be a subalgebra of $\R$.
    The algebra $\alg{M}$ is nontrivial because it contains $\Z$, and thus the congruence lattice of $\alg{M}$ has at least two elements.
    Let $\theta \in \con(\alg{M})$ with $0_{\con(\alg{M})} < \theta$.
    Then there is a pair $(a,b) \in \theta$ with $a < b$.
    By the classical Archimedean property of $\R$, there exists $n \in \N$ such that $nb - na \geq 2$, and hence there is $k \in \Z$ such that $n a \leq k$ and $k + 1 \leq nb$.
    Since $\theta$ is a lattice congruence, from $(a,b) \in \theta$ we deduce $(na, nb) \in \theta$.
    From $na \leq k \leq k+1 \leq nb$ and $(na, nb) \in \theta$ we deduce $(k, k + 1) \in \theta$.
    It follows that $(0, 1) = (k , k + 1) +(-k, -k) \in \theta$ and then it is easy to prove that for all $c,d \in \Z$ we have $(c, d) \in \theta$.
    Since $\theta$ is a lattice congruence, we conclude that $1_{\con(\alg{M})} = \theta$.
\end{proof}

Using \cref{preserveandreflect,c:isomorphic-lattices-of-congruences} we get the next corollary as an immediate consequence of \cref{l:hsimple}.

\begin{corollary}\label{lem:LnSimple}
    Every subalgebra of the MV-monoid $[0,1]^+$ is simple.
\end{corollary}

In particular, for every $n \in \N$, the algebra $\alg{\L}_n^+$ is simple and hence subdirectly irreducible.

For our description of the almost minimal varieties of MV-monoids, we need to introduce the three-element algebras $\alg{C}_2^\Delta$ and $\alg C_2^\nabla$.

The algebra $\alg{C}_2^\Delta$ is the unique MV-monoid on the 3-element chain $0 < \e < 1$ satisfying $\e \oplus \e = \e$ and $\e \odot \e = 0$. (Roughly speaking, $\e$ is an infinitesimal element).
Equivalently (up to isomorphism), $\alg{C}_2^\Delta$ can be defined as $\Gamma(\Z \overrightarrow{\times} \alg{C}_2^{\Delta*})$ or as the quotient of $\alg{C}^\Delta = \la \{0, \e, 2\e, \dots, 1\}, \join,\meet,\oplus,\odot,0,1\ra$ by the equivalence relation that identifies all elements that are neither $0$ nor $1$.

The algebra $\alg C_2^\nabla$ is defined ``dually'' as the unique MV-monoid on the 3-element chain $0 < \d < 1$ satisfying $\d \oplus \d = 1$ and $\d \odot \d = \d$.
Equivalently, $\alg C_2^\nabla$ can be defined as $\Gamma(\Z \overrightarrow{\times} \alg{C}_2^{\nabla*})$ or as the quotient of $\alg{C}^\nabla = \la \{0, \dots,\delta^3, \delta^2, \delta,  1\}, \join,\meet,\oplus,\odot,0,1 \ra$ by the equivalence relation that identifies all elements that are neither $0$ nor $1$.

We will prove that $\alg C^\Delta_2$ and $\alg C^\nabla_2$ are subdirectly irreducible, and that $\alg C^\Delta$ and $\alg C^\nabla$ are relatively subdirectly irreducible positive MV-algebras that are not subdirectly irreducible in the absolute sense (\cref{sec: si pMV}).
To do so, we characterize the congruence lattices of $\alg C^\Delta$ and $\alg C^\nabla$, a fact of independent interest.
The lattice reducts of $\alg C^\Delta$ and $\alg C^\nabla$ are chains.
For every chain $\alg L$ regarded as a distributive lattice, and for all $a, b \in {L}$ with $a < b$, it is easily seen that the principal congruence $\cg_\alg L(a,b)$ is the congruence in which the only nontrivial block is the interval $[a,b]$.
It is a nice exercise in lattice theory to prove that:

\begin{lemma}
    If $\alg L$ is the $n+1$-element chain, then $\Con L$ is isomorphic to the $2^n$-element complemented distributive lattice.
\end{lemma}

In particular, if $\alg A$ is a finite totally ordered MV-monoid, then any of its congruences is a lattice congruence, and so it is the join of finitely many lattice congruences of the form described above. In fact, the congruence classes of a congruence on a finite lattice are always intervals.
    \begin{proposition}\label{techlemma}  Let $\th$ be a congruence of $\alg{C}^\Delta$.
    \begin{enumerate}
        \item \label{i:notintheta}
        If $\th <1_{\con(\alg{C}^\Delta)}$, then $(n\e,1) \notin \th$ for all $n \in \mathbb N$;

        \item \label{i:intheta}
        for all $n,m \in \mathbb N$ with $n < m$, if $(n\e,m\e) \in \th$, then $(n\e,l\e) \in \th$ for all $l \ge n$.
    \end{enumerate}
\end{proposition}
\begin{proof}
    \eqref{i:notintheta}.
    We prove the contrapositive.
    Suppose $(n\e,1) \in \th$ for some $n \in \mathbb N$.
    Then
    $$
    (0,1) = (n\e \odot n\e,1 \odot1) \in \th.
    $$
    Since $\th$ is a lattice congruence, from $(0,1) \in \th$ we deduce $\th = 1_{\con(\alg{C}^\Delta)}$.

    \eqref{i:intheta}.
    If $(n\e,m\e) \in \th$ with $n < m$, then $n < n+1 \le m$ and thus $(n\e,(n+1)\e) \in \th$, since $\th$ is a lattice congruence. Next we observe that
    $$
    ((n+1)\e,(n+2)\e) = (n\e \oplus \e,(n+1)\e \oplus \e) \in\th.
    $$
    By transitivity of $\theta$, from $(n\e, (n + 1)\e) \in \theta$ and $((n+1)\e,(n+2)\e) \in \th$  we deduce $(n\e,(n+2)\e) \in \th$. The conclusion follows by induction.
\end{proof}

\begin{corollary} \label{cor:Cdeltacon}
    The congruence lattice of $\alg{C}^\Delta$ is an infinite chain, ordered as the order dual of $\omega + 1$, and the proper nontrivial congruences of $\alg{C}^\Delta$ are exactly those whose associated partition is
    $$
    \g^\Delta_n \coloneqq \{\{0\},\{\e\},\{2\e\},\dots,\{(n-1)\e\},\{n\e,\dots\},\{1\}\}
    $$
    for some $n\ge 0$.
\end{corollary}

\begin{proof}
    If $\th$ is proper, then $1/\th=\{1\}$ by \cref{techlemma}.
    Since $\th$ is nontrivial we have $(n \e,m\e) \in \th$ for some $n <m$.
    Since the set $\{k \in \N \mid (k\e,m\e) \in \th\}$ is nonempty, it has minimum $\bar n$. By \cref{techlemma}\eqref{i:intheta}, $\th = \g_{\bar n}$.
\end{proof}

With a totally similar argument we can prove that the nontrivial proper congruences of $\alg{C}^\nabla$ are exactly those whose associated partition is
\begin{align*}
    &\g^\nabla_n = \{\{0\},\{\dots,\d^n\},\{\d^{n-1}\},\dots ,\{\d\},\{1\}\},
\end{align*}
and thus also the congruence lattice of $\alg{C}^\nabla$ is an infinite chain.

For each $n \in \N \setminus \{0\}$, we let $\alg L_n $ be the $(n+1)$-element chain seen as a distributive lattice.
The lattice $\alg{L}_n$ is the lattice reduct of various MV-monoids.
We have already introduced one of these: $\alg{\L}_n^+ = \Gamma(\frac{1}{n}\Z)$, which is a simple and hence subdirectly irreducible positive MV-algebra.
The following are further examples.
\begin{enumerate}
    \ib $\alg{C}^\Delta_n \coloneqq \alg{C}^\Delta/\g^\Delta_{n-1} \cong \Gamma(\Z \overrightarrow{\times}\alg{C}^{\Delta*}_n)$;
    \ib $\alg{C}^\nabla_n \coloneqq \alg{C}^\nabla/\g^\nabla_{n-1} \cong \Gamma(\Z \overrightarrow{\times}\alg{C}^{\nabla*}_n)$.
\end{enumerate}

We note that, for $n = 1$, $\alg{\L}_n^+$, $\alg{C}^\Delta_n$ and $\alg{C}^\nabla_n$ are all isomorphic to the 2-element bounded distributive lattice.
For all $n \ge 2$, however, $\alg {\L}^+_n$, $\alg{C}^\Delta_n$ and $\alg{C}^\nabla_n$ are pairwise nonisomorphic, and $\alg{C}^\Delta_n$ and $\alg{C}^\nabla_n$ are not positive MV-algebras because the monoids $\Z \overrightarrow{\times}\alg{C}^{\Delta*}_n$ and $\Z \overrightarrow{\times}\alg{C}^{\nabla*}_n$ are not cancellative. Finally, we show that $\alg{C}^\Delta_n$ and $\alg{C}^\nabla_n$ share with $\alg{\L}_n^+$ the property of being subdirectly irreducible.
\begin{lemma}\label{l:Bn-subirr}
    For every $n \in \N \setminus \{0\}$, both $\alg{C}^\Delta_n$ and $\alg{C}^\nabla_n$ are subdirectly irreducible.
\end{lemma}
\begin{proof}
    By \cref{cor:Cdeltacon} the congruence lattice of $\alg{C}^\Delta$ is an infinite chain. As $\alg{C}^\Delta_n \cong \alg{C}^\Delta/\g^\Delta_{n-1}$, by the Correspondence Theorem \cite[Theorem $6.20$]{BurrisSanka} we get that $\op{Con}(\alg{C}^\Delta_n)$ is the $n+1$-element chain isomorphic to the interval $[\g^\Delta_{n-1}, 1]$ of $\con(\alg{C}^\Delta)$, and the same holds {\em mutatis mutandis} for $\alg{C}^\nabla_n$. Hence they are both subdirectly irreducible.
\end{proof}
The following important result will be used in the proof of the main result of the section (\cref{t:char-almost-minimal}) and other statements in the remaining part of the paper.
Moreover, as we will discuss after its corollary for MV-monoids, it implies a version for commutative $\ell$-monoids of the celebrated H\"older's theorem for Abelian $\ell$-groups.

\begin{theorem}\label{Holder}
    Let $\alg{M}$ be a nontrivial totally ordered unital commutative $\ell$-monoid. There is a unique homomorphism from $\alg{M}$ to $\R$.
\end{theorem}

\begin{proof}
    For all $x \in M$, we define
    \begin{align*}
        I_x & \coloneqq \left\{\frac{k}{n} \in \Q \mid k \in \Z,\, n \in \N \setminus \{0\},\, nx \geq k\right\};\\
        S_x & \coloneqq \left\{ \frac{k}{n} \in \Q \mid k \in \Z,\, n \in \N \setminus \{0\},\, nx \leq k \right\}.
    \end{align*}
    We prove $\sup I_x = \inf S_x \in \R$.
    To do so, it is enough to prove $I_x \neq \emptyset \neq S_x$, $I_x \cup S_x = \Q$, that for all $x_1 \in I_x$ and $x_2 \in S_x$ $x_1 \leq x_2$, and that $I_x \cap S_x$ has at most one element.
    Since every element is bounded from above and below by some integer multiples of $1$, $I_x \neq \emptyset \neq S_x$.
    Furthermore, let $\frac{k_1}{n_1} \in I_x$ and $\frac{k_2}{n_2} \in S_x$. By monotonicity of $+$ we have $\frac{n_2k_1}{n_1n_2} \in I_x$ and $\frac{n_1k_2}{n_1n_2} \in S_x$. Then $n_2k_1 \leq n_1n_2x \leq n_1k_2$ and thus $\frac{k_1}{n_1} \leq \frac{k_2}{n_2}$.
    Moreover, since $\alg{M}$ is totally ordered, $I_x \cup S_x = \Q$. We prove that $I_x \cap S_x$ has at most an element.
    Suppose there are $k,k' \in \Z$ and $n, n' \in \N \setminus \{0\}$ such that $k = nx $ and $k' = n' x$.
    We shall prove that the rational numbers $\frac{k}{n}$ and $\frac{k'}{n'}$ coincide, i.e., that $kn' = k'n$.
    Without loss of generality, we may assume $kn' \leq k'n$.
    Since $k = nx $ and $k' = n' x$, we have $kk' = k'nx = kn'x$, and thus $(k'n - kn')x = 0$,
    since the element $kk' = k'nx$ is invertible.
    Then, either $x = 0$ or $k'n - kn' = 0$. If $x = 0$, then the interpretations of $k$ and $k'$ in $\alg{M}$ are $0$, which implies that the integers $k$ and $k'$ are $0$, by nontriviality of $\alg{M}$, and hence $k'n = 0 = kn'$. Thus in both cases we have $k'n = kn'$ and this proves that $I_x \cap S_x$ has at most an element.
    Therefore, $\sup I_x = \inf S_x \in \R$.
    We define the map
    \begin{align*}
        \varphi \colon \alg{M} & \longrightarrow \R\\
        x & \longmapsto \sup I_x = \inf S_x.
    \end{align*}
    We prove that $\varphi$ is a homomorphism. It is not difficult to see that it preserves all constant symbols.

    Let $x, y \in M$.
    First, we prove $\varphi(x) + \varphi(y) = \varphi(x + y)$.
    To do so, we prove that the element-wise sum $I_x + I_y$ is contained in $I_{x + y}$.
    Let $\frac{k}{n} \in I_x$ (i.e., $nx \geq k$) and $\frac{k'}{n'} \in I_y$ (i.e., $n'y \geq k'$).
    We shall prove $\frac{k}{n} + \frac{k'}{n'} \in I_{x + y}$, i.e., $\frac{n'k + nk'}{nn'} \in I_{x + y}$, i.e., $nn' (x + y) \geq n'k + nk'$. Since $nx \geq k$ and $n'y \geq k'$, we have respectively $n'nx \geq n'k$ and $nn'y \geq nk'$, by monotonicity of $+$. Thus $nn'(x + y) \geq n'k + nk'$, still by monotonicity.
    This proves $I_x + I_y \subseteq I_{x + y}$.
    Therefore, $\sup (I_x + I_y) \leq I_{x + y}$.
    Since $+ \colon \R^2 \to \R$ is continuous, we have $\sup (I_x + I_y) = \sup I_x + \sup I_y$.
    Therefore,
    \[
    \varphi(x) + \varphi(y) = \sup I_x + \sup I_y = \sup(I_x + I_y) \leq \sup I_{x+ y} = \varphi(x + y).
    \]
    Replacing in the proof above $I_x$ with $S_x$, $\sup$ with $\inf$, and reversing the inequalities, we obtain $\varphi(x) + \varphi(y) \geq \varphi(x + y)$ and hence $\varphi(x) + \varphi(y) = \varphi(x + y)$.

    For the compatibility with the lattice operations let $x', y' \in M$. If $x' \leq y'$ then $I_{x'} \subseteq I_{y'}$ and hence $\varphi(x') \leq \varphi(y')$.
    Thus $\varphi$ respects the order and hence $\varphi(x) \lor \varphi(y) = \varphi(x \lor y)$ and $\varphi(x) \land \varphi(y) = \varphi(x \land y)$, since $\alg M$ and $\R$ are totally ordered. Therefore, $\varphi$ is a homomorphism.

    The uniqueness follows from the fact that for every homomorphism $\psi \colon \alg{M} \to \R$ and every $x \in M$ the condition $nx  \geq k$ implies $n\psi(x) \geq k$ and the condition $nx \leq k$ implies $n\psi(x) \leq k$.
    Therefore, for all $\frac{k}{n} \in I_x$, we have $nx  \geq k$ and thus $\psi(x) \geq \frac{k}{n}$. Hence, $\psi(x) \geq \sup I_x = \varphi(x)$. With a symmetric proof using $S_x$ we see that $\psi(x) \leq \varphi(x)$ and thus $\psi = \varphi$.
\end{proof}

\begin{corollary}\label{c:Holder-MV}
    Let $\alg{A}$ be a nontrivial totally ordered MV-monoid such that for all $x, y \in A$ we have $x \oplus y = 1$ or $x \odot y = 0$. There is a unique homomorphism from $\alg{A}$ to $[0,1]^+$.
\end{corollary}

\begin{proof}
    This follows from \cref{prop: M tot ord}, \cref{Holder}, the categorical equivalence between MV-monoids and unital commutative $\ell$-monoids, and the fact that $\Gamma(\R) = [0,1]^+$.
\end{proof}

As a corollary of \cref{Holder}, since every monoid homomorphism between groups is a group homomorphism, we get the corresponding result for unital Abelian $\ell$-groups:

\begin{corollary} \label{c:cor-abel}
    Let $\alg{G}$ be a nontrivial totally ordered unital Abelian $\ell$-group. Then, there is a unique homomorphism of unital Abelian $\ell$-groups from $\alg{G}$ to $\R$.
\end{corollary}

Likewise, every MV-monoid homomorphism between MV-algebras is an MV-algebra homomorphism because, for every $a$ in an MV-algebra, $\lnot a$ is the unique $c$ such that $a \oplus c = 1$ and $a \odot c = 0$. Therefore, \cref{c:Holder-MV} implies the corresponding result for MV-algebras:

\begin{corollary}
    Let $\alg{A}$ be a nontrivial totally ordered MV-algebra. There is a unique MV-algebra homomorphism from $\alg{A}$ to $[0,1]$.
\end{corollary}

\begin{proof}
    In light of \cref{c:Holder-MV} and the fact that every MV-monoid homomorphism between MV-algebras is an MV-algebra homomorphism, it is left to observe that, for every totally ordered MV-algebra $\alg{A}$, for all $x, y \in A$ we have $x \oplus y = 1$ or $x \odot y = 0$.
    This is true because $x \oplus y = 1$ is equivalent to $y \geq \lnot x$, and $x \odot y = 0$ is equivalent to $y \leq \lnot x$.
\end{proof}

\Cref{Holder} has also the following consequences.

\begin{corollary}\label{c:if-mon-simple-then-unique}
    For every simple unital commutative $\ell$-monoid $\alg{M}$ there is a unique homomorphism from $\alg{M}$ to $\R$, and it is injective.
\end{corollary}

\begin{proof}
    Let $\alg{M}$ be a simple unital commutative $\ell$-monoid.
    Since it is simple, it is subdirectly irreducible.
    Thus, by \cref{c:sub-totord}, it is totally ordered (and clearly also nontrivial).
    Hence, by \cref{Holder}, there is a unique homomorphism $\varphi$ from $\alg{M}$ to $\R$.
    The image of $\varphi$ contains $\Z$ and hence is not trivial.
    Therefore, by simplicity of $\alg{M}$, the kernel congruence of $\varphi$ is the identity on $\alg{M}$, i.e., $\varphi$ is injective.
\end{proof}

\begin{corollary} \label{c:simple-char-mon}
    The simple unital commutative $\ell$-monoids are precisely the algebras isomorphic to a subalgebra of the unital commutative $\ell$-monoid $\R$.
\end{corollary}

\begin{proof}
    The right-to-left inclusion is \cref{l:hsimple} and the left-to-right inclusion is an immediate consequence of \cref{c:if-mon-simple-then-unique}.
\end{proof}

From \cref{c:if-mon-simple-then-unique,c:simple-char-mon}, using the categorical equivalence between MV-monoids and unital commutative $\ell$-monoids, we get the following two corollaries.

\begin{corollary} \label{c:if-MV-simple-then-unique}
    For every simple MV-monoid $\alg{A}$ there is a unique homomorphism from $\alg{A}$ to $[0,1]^+$, and it is injective.
\end{corollary}

\begin{corollary} \label{c:simple}
    The simple MV-monoids are precisely the algebras isomorphic to a subalgebra of $[0,1]^+$.
\end{corollary}

\Cref{c:simple} has the following immediate consequence. We recall that an algebra is \emph{semisimple} if it is isomorphic to a subdirect product of simple algebras.

\begin{corollary} \label{c:semisimple}
    The semisimple MV-monoids are precisely the algebras isomorphic to a subalgebra of a power of $[0,1]^+$.
\end{corollary}

\begin{remark}
    \Cref{c:if-mon-simple-then-unique,c:simple-char-mon,c:if-MV-simple-then-unique,c:simple,c:semisimple} immediately entail their versions for unital Abelian $\ell$-groups and MV-algebras (see, e.g., \cite[Lemma~3.8]{MarSpad2012}, \cite[Theorem~3.5.1]{CignoliDOttavianoEtAl2000}, and \cite[Proposition~3.6.1]{CignoliDOttavianoEtAl2000} for the versions for MV-algebras).
    
    Indeed, we observe that every monoid congruence on a group is a group congruence, since $x \sim y$ implies $xy^{-1} \sim yy^{-1} = 1$ which implies $y^{-1} = x^{-1}xy^{-1} \sim x^{-1}1$.
    Thus, every unital $\ell$-monoid congruence on a unital Abelian $\ell$-group is a unital $\ell$-group congruence.
    From this, we can deduce that every MV-monoid congruence on an MV-algebra is an MV-algebra congruence; indeed, one can use the equivalence between unital commutative $\ell$-monoids and MV-monoids, the fact that it restricts to Mundici's equivalence between unital Abelian $\ell$-groups and MV-algebras, and the fact that it preserves and reflects surjectivity of maps.
    Alternatively, one can give a simple direct proof analogous to the one for monoids and groups: $x \sim y$ implies $x \odot \lnot y \sim y \odot \lnot y = 0$, which implies $\lnot y = \lnot x \oplus (x \odot \lnot y) \sim \lnot x \oplus 0 = \lnot x$.
    Therefore, every (semi)simple unital Abelian $\ell$-group has a (semi)simple unital commutative $\ell$-monoid reduct, and, likewise, every (semi)simple MV-algebra has a (semi)simple MV-monoid reduct.
\end{remark}

If we replace MV-monoids with MV-algebras, \cref{c:if-MV-simple-then-unique} is what some authors refer to as H\"{o}lder's theorem for MV-algebras \cite[Lemma~3.8]{MarSpad2012}.
H\"older's Theorem for Archimedean ordered groups (\cite[Theorem~2.6.3]{BigardKeimelEtAl1977}), which goes back to \cite{Hoelder1901}, asserts that every totally ordered Archimedean group is isomorphic to a subalgebra of $\R$.
To make the connection between \cref{Holder} and H\"older's theorem clearer, we show that \cref{Holder} implies a version of H\"older's Theorem for unital commutative $\ell$-monoids.

We recall that a unital Abelian $\ell$-group $\alg{G}$ is \emph{Archimedean} provided that, for every $y \in G$, if for every $n \in \N$ we have $ny \leq 1$, then $y \leq 0$.
In the context of unital commutative $\ell$-monoids, we propose the following definition.

\begin{definition}\label{d:arch}
    A unital commutative $\ell$-monoid $\alg{M}$ is \emph{Archimedean} provided that, for all $x,y \in M$, if for all $n \in \N$ we have $nx \leq ny + 1$, then $x \leq y$.
\end{definition}

This definition is motivated by the following observation.

\begin{proposition} \label{p:Arch}
    Archimedean cancellative unital commutative $\ell$-monoids are precisely the subreducts of Archimedean unital Abelian $\ell$-groups.
\end{proposition}

\begin{proof}
    It is easy to see that every $\{\vee, \wedge, +, 1, 0, -1\}$-subreduct of an Archimedean unital Abelian $\ell$-group is an Archimedean unital commutative $\ell$-monoid.

    Let $\alg{M}$ be an Archimedean cancellative unital commutative $\ell$-monoid.
    Since $\alg{M}$ is a cancellative unital commutative $\ell$-monoid, it is the subreduct of a unital Abelian $\ell$-group $\alg{G}$ \cite[Proposition~4.4]{Abbadinietal2022}.
    Without loss of generality, we can assume that $\alg{M}$ generates $\alg{G}$.
    Therefore, every element of $G$ is of the form $x -y$ for $x,y \in M$ (cf.\ \cite[Lemma~4.2]{Abbadinietal2022}).
    Let us prove that $\alg{G}$ is Archimedean.
    Let $z \in G$ and suppose that for every $n \in \N$ we have $nz \leq 1$.
    We should prove $z \leq 0$.
    There are $x,y \in M$ such that $z = x - y$.
    For every $n \in \N$ we have
    \[
    nx - ny = n(x -y) = nz \leq 1,
    \]
    which implies $nx \leq ny + 1$.
    Since $\alg{M}$ is Archimedean as a unital commutative $\ell$-monoid, we deduce $x \leq y$.
    Thus, $z = x - y \leq 0$.
\end{proof}

We are now ready to prove a version of H\"older's theorem for unital commutative $\ell$-monoids.

\begin{theorem}[H\"older's theorem for unital commutative $\ell$-monoids]\label{c:Holder}
    Let $\alg{M}$ be an Archimedean nontrivial totally ordered unital commutative $\ell$-monoid. The unique homomorphism from $\alg{M}$ to $\R$ is injective, and so $\alg{M}$ is isomorphic to a subalgebra of $\R$.
\end{theorem}

\begin{proof}
    Let $\varphi$ be the unique homomorphism from $\alg{M}$ to $\R$.
    Let $x, y \in M$ be such that $\varphi(x) = \varphi(y)$, and let us prove that $x = y$.
    Since $\alg{M}$ is totally ordered, either $x \leq y$ or $y \leq x$.
    Without loss of generality, we can suppose $x \leq y$.
    We shall prove $y \leq x$.
    By the definition of Archimedean unital commutative $\ell$-monoid (\cref{d:arch}), it suffices to prove that, for every $n \in \N$, $ny \leq nx + 1$.
    Let $n \in \N$.
    We have
    \begin{equation}\label{eq:strict}
        \varphi(ny)  = n \varphi(y) < n\varphi(y) + 1 = n\varphi(x) + 1 = \varphi(nx +1).
    \end{equation}
    Since $\alg{M}$ is totally ordered, either $ny \leq nx + 1$ or $nx + 1 \leq ny$.
    The second case is not possible because otherwise we would have $\varphi(nx + 1) \leq \varphi(ny)$ by monotonicity of $\varphi$, contradicting \eqref{eq:strict}.
    Thus, the first case holds, i.e., $ny \leq nx + 1$.
\end{proof}

\begin{remark}
    \cref{c:Holder} implies H\"older's Theorem for Archimedean unital ordered groups, since monoid homomorphisms between groups are group homomorphisms.
\end{remark}
To proceed further, we need some definitions.
\begin{notation}
    Let $\alg{M}$ be a nontrivial totally ordered unital commutative $\ell$-monoid, and let $\varphi$ be the unique homomorphism from $\alg{M}$ to $\R$ (which exists by \cref{Holder}).
    Recalling from \cref{d:unital-commutative-ell-monoid} that every integer $n$ has an interpretation in $\alg{M}$ (the $n$-fold sum of $1$), for every $z \in \varphi^{-1}[\Z]$ we denote by $K_z$ the interpretation in $\alg{M}$ of the integer $\varphi(z)$.
\end{notation}

\begin{lemma}\label{l: T subalgebra}
    Let $\alg{M}$ be a nontrivial totally ordered unital commutative $\ell$-monoid, and let $\varphi$ be the unique homomorphism from $\alg{M}$ to $\R$ (which exists by \cref{Holder}).
    The following are universes of subalgebras of $\alg{M}$:
    \begin{enumerate}
        \item \label{i:S} $S \coloneqq \varphi^{-1}[\Z]$;
        \item \label{i:T}$T \coloneqq \{z \in S \mid K_z \leq z\}$;
        \item \label{i:Z}$Z \coloneqq \{z \in T \mid K_z = z\}$.
    \end{enumerate}
\end{lemma}

\begin{proof}
    \eqref{i:S}.
    $\alg{S}$ is a subalgebra of $\alg{M}$ because it is the preimage of $\Z$ (which is a subalgebra of $\R$) under the homomorphism $\varphi$.

    \eqref{i:T}. Since $\alg{M}$ is totally ordered, $T$ is closed under $\lor$ and $\land$.
    Moreover, $T$ contains $0$, $1$ and $-1$. To prove closure under $+$, let $z_1, z_2 \in T$.
    We already know that $z_1 + z_2 \in S$.
    Moreover,
    \[
    K_{z_1 + z_2} = K_{z_1} + K_{z_2} \leq z_1 + z_2,
    \]
    and thus $z_1 + z_2 \in T$, proving the closure of $T$ under $+$. This proves \eqref{i:T}.

    \eqref{i:Z}. This can be proven using an argument similar to the one for \eqref{i:T}, or using the fact that $Z$ is the image of the unique homomorphism from $\Z$ to $\alg{M}$.
\end{proof}	

\begin{theorem} \label{t:for-almost-minimal}
    Let $\alg{M}$ be a nontrivial totally ordered unital commutative $\ell$-monoid which is not isomorphic to $\mathbb Z$ (regarded as a unital commutative $\ell$-monoid).
    Then, at least one of the following conditions holds:
    \begin{enumerate}
        \item \label{th:l1}
        The image of the unique homomorphism from $\alg{M}$ to $\R$ is not $\Z$;
        \item \label{th:l2}$\Z \overrightarrow{\times}\alg{C}^{\Delta*}_2 \in \HH\SU(\alg{M})$;
        \item \label{th:l3} $\Z \overrightarrow{\times}\alg{C}^{\nabla*}_2 \in \HH\SU(\alg{M})$.
    \end{enumerate}
\end{theorem}

\begin{proof}
    Since $\alg{M}$ is a nontrivial totally ordered unital commutative $\ell$-monoid, by \cref{Holder} there is a unique homomorphism $\varphi \colon \alg{M} \to \R$.
    Clearly, either the image $\alg{A}$ of $\alg{M}$ is $\mathbb Z$, or it strictly contains $\mathbb Z$.

    In the second case, \eqref{th:l1} holds.
    So, we are left to analyze the case in which the image $\alg{A}$ of $\alg{M}$ under $\varphi$ is $\mathbb Z$.
    In this case, $\varphi$ is not injective, because $\alg{M}$ is not isomorphic to $\mathbb Z$.
    Therefore, there are $a \in M$ and $n \in \mathbb Z$ such that $\varphi(a) = n$ and $a \neq n$.
    Up to replacing $a$ with $a - n$, we can suppose $n = 0$ and,
    since $\alg{M}$ is totally ordered, either $a \geq 0$ or $a \leq 0$.
    Let us suppose first that $a \geq 0$.
    Recall that, for every $x \in \varphi^{-1}[\Z]$, $K_x$ denotes the interpretation in $\alg{M}$ of the integer $\varphi(x)$.
    Recall also that $T \coloneqq \{x \in \varphi^{-1}[\Z] \mid x \geq K_x\}$.
    Then $\alg{T}$ is a subalgebra of $\alg{M}$, by \cref{l: T subalgebra}.
    We define the function
    \begin{align*}
        \psi \colon T & \longrightarrow \Z \overrightarrow{\times}\alg{C}^{\Delta*}_2\\
        x &\longmapsto \begin{cases}
        (K_x,0) & \text{if }x = K_x,\\
        (K_x, \e) & \text{otherwise.}
    \end{cases}
    \end{align*}
    We show that $\psi$ is a homomorphism.
    We have $K_0 = 0$, $K_1 = 1$, and $K_{-1} = -1$, and so $\psi$ preserves $0$, $1$ and $-1$.
    We prove that $\psi$ is order-preserving.
    Let $x, y \in T$ be such that $x \leq y$. Then $\varphi(x) \leq \varphi(y)$, and so $K_x \leq K_y$. If $K_x < K_y$, then $\psi(x) < \psi(y)$ by inspection on the first coordinates of $\psi(x)$ and $\psi(y)$.
    If $K_x = K_y$, either $y = K_x$ and then $x = K_x$, since $x \leq y$ and $\varphi(x) = \varphi(y) = K_x$, or $y > K_x$ and thus $\psi(x) = (K_x, z) \leq (K_x, \e) = \psi(y)$, for some $z \in \{0,\e\}$.
    This proves that $\psi$ is order-preserving. Since $T$ is totally ordered, it follows that $\psi$ is a lattice homomorphism.

    We prove that $\psi$ preserves $+$. Let $x,y \in T$ and let $(K_x,z_x) = \psi(x)$ and $(K_y,z_y) = \psi(y)$ for some $z_x,z_y \in \{0,\e\}$. We prove that
    \[
    \psi(x + y) = (K_x + K_y, \max(\{z_x,z_y\})) = \psi(x) + \psi(y).
    \]
    If $x = K_x$ and $y = K_y$, then
    \[
    \psi(x + y) = \psi(K_x + K_y) =  (K_x + K_y, 0) = (K_x + K_y, \max(\{z_x,z_y\})) = \psi(x) + \psi(y).
    \]
    If $x > K_x$ or $y > K_y$, then
    \[
    \varphi(x+y) = \varphi(x) + \varphi(y) = K_x + K_y
    \]
    and $x+y > K_x+K_y$. Thus
    \[
    \psi(x + y) =  (K_x + K_y, \e) = (K_x + K_y, \max(\{z_x,z_y\})) = \psi(x) + \psi(y).
    \]

    We then show that $\psi$ is surjective. Let $(K,y) \in \Z \overrightarrow{\times} \alg C^{\Delta*}_2$. We have two cases. If $y = 0$, then $(K,y) = (K,0) = \psi(K)$. If $y = \e$, then $\varphi(a + K) = \varphi(a) + \varphi(K) = K$ and thus $\psi(a + K) = (K,\e) = (K,y)$ because $a + K \neq K$ (since $K$ is invertible and $a \neq 0$).

    Therefore, $\Z \overrightarrow{\times}\alg{C}^{\Delta*}_2 \in \HH(\alg{T}) \subseteq \HH\SU(\alg{M})$.
    If instead $a \leq 0$, we can prove in a similar way that $\Z \overrightarrow{\times}\alg{C}^{\nabla*}_2 \in \HH\SU(\alg{M})$.
\end{proof}

\begin{corollary} \label{c:almost-minimal-MV}
    Let $\alg{A}$ be a totally ordered MV-monoid with at least 3 elements such that for all $x,y \in A$ we have $x \oplus y = 1$ or $x \odot y = 0$.
    At least one of the following conditions holds:
    \begin{enumerate}
        \item The image of the unique homomorphism from $\alg{A}$ to $[0,1]$ is not $\{0,1\}$;
        \item $ \alg{C}^\Delta_2 \in \HH\SU(\alg{A})$;
        \item $\alg{C}^\nabla_2 \in \HH\SU(\alg{A})$.
    \end{enumerate}
\end{corollary}

\begin{proof}
    The statement follows from \cref{t:for-almost-minimal} using the equivalence between $\ulm$ and $\MVM$, \cref{prop: M tot ord} and the fact that $\Gamma$ preserves and reflects injectivity and surjectivity (\cref{preserveandreflect}).
\end{proof}

It is a known fact that every additive subgroup of $\R$ is either dense or cyclic; see e.g.\ \cite[Lemma~4.21]{Goodearl1986}.
Something similar holds for unital commutative $\ell$-monoids.

\begin{lemma} \label{l:dense-or-discrete}
    Let $\alg{M}$ be a subalgebra of the unital commutative $\ell$-monoid $\R$.
    Either $M$ is dense in $\R$ or there is $n \in \N \setminus \{0\}$ such that $\alg{M} = \frac{1}{n}\Z$.
\end{lemma}

\begin{proof}
    Suppose that $M$ is not dense in $\R$.
    The set $\{x \in M \mid x > 0\}$ is nonempty (because $1 \in M$) and is bounded from below by $0$; therefore, it admits an infimum
    \[
    \e \coloneqq \inf\{x \in M \mid x > 0\}
    \]
    in $\R$.
    There is a sequence $(a_n)_{n \in \N}$ of elements in $M$ greater than $\e$ that converges to $\e$.

    We claim that $\e$ is strictly greater than $0$.
    Suppose, by way of contradiction, that this is not the case.
    Then we prove that $M$ is dense.
    Let $x \in \R$ and $a > 0$, and let us show the existence of an element $y \in M$ such that $\lvert x - y \rvert < a$.
    There is $n \in \N$ such that $a_n < 2a$.
    Then, the set $\{k a_n + l \mid k \in \N, l \in \Z \}$ intersects the interval $(x-a, x + a)$ nontrivially. Indeed, the nonemptiness of this intersection amounts to the existence of an integer $l$ for which there is $k \in \N$ in the interval $[\frac{x -a - l}{a_n},\frac{x + a - l}{a_n}]$. Note that this interval has length $\frac{2a}{a_n}$, which is greater than $1$. Moreover, the endpoints of this interval are positive if $l \leq x - a$.
    Therefore, it is enough to take any integer $l$ below $x - a$ and a natural number $k$ in $[\frac{x -a - l}{a_n},\frac{x +a - l}{a_n}]$, which exists because the latter is an interval of length greater than $1$ with positive endpoints.
    Since $\{k a_n + l \mid k \in \N, l \in \Z \}$ is a subset of $M$, there is $y \in M$ such that $\lvert x - y \rvert < a$.
    This proves our claim that $\e$ is strictly greater than $0$.
    Analogously, the set $\{x \in M \mid x < 0\}$ has a strictly negative supremum
    \[
    \d \coloneqq \sup\{x \in M \mid x < 0\}
    \]
    in $\R$.
    There is a sequence $(b_n)_{n \in \N}$ of elements in $M$ smaller than $\d$ that converges to $\d$.

    We claim that $\d = - \e$ and that $\e$ and $\d$ belong to $M$.
    Indeed, the sequence $(a_n + b_n)_{n \in \N}$ of elements of $M$ converges to $\e + \d$, and $\e + \d$ is strictly between $\d$ and $\e$.
    Since $0$ is the only element in $M$ in the interval $(\d, \e)$, the sequence $(a_n + b_n)$ is eventually equal to $0$, and thus $\e + \d = \lim_{n \to \infty} a_n + \lim_{n \to \infty} b_n = \lim_{n \to \infty} (a_n + b_n) = 0$.
    This proves $\d = -\e$.
    Analogously, the sequence $(a_n + b_{n+1})_{n \in \N}$ is eventually equal to $0$.
    Therefore, eventually, $b_{n} = - a_n = b_{n + 1} = -a_{n + 1}$.
    Therefore, the sequences $(a_n)_{n \in \N}$ and $(b_n)_{n \in \N}$ are eventually constant.
    Therefore, $\e$ and $\d$ belong to $M$, proving our claim.

    It is now clear that every element of the form $k \e$ for some $k \in \Z$ belongs to $M$.
    We claim that also the converse holds, i.e., that every element of $M$ is of the form $k \e$ for some $k \in \Z$.
    Indeed, by way of contradiction, suppose this were not the case; there would exist $x \in M$ which is not of the form $k \e$ ($k \in \Z$).
    Since $\e >0$, there is $k \in \Z$ such that $k\e < x < (k+1) \e$.
    But then the element $0 < x - k\e < \e$ and $x - k\e = x + k\d \in M$, contradicting the definition of $\e$.
    This proves our claim that every element of $M$ is of the form $k \e$ for some $k \in \Z$.

    In particular, also $1$ is of this form, and so there is $n \in \N$ such that $1 = n\e$, i.e., $\e = \frac{1}{n}$.
    Thus, $M = \frac{1}{n}\Z$.
\end{proof}

Note that a subalgebra $\alg{M}$ of the unital commutative $\ell$-group $\R$ is dense in $\R$ if and only if $\Gamma(\alg{M})$ is dense in $[0,1]$. Indeed, the left-to-right implication is immediate. The right-to-left implication, instead, follows from the fact that if $\Gamma(\alg{M})$ is dense in $[0,1]$, then for every $k \in \Z$ the set $M \cap [k, k+1]$ is dense in $[k, k+1]$ because it is the image of $\Gamma(\alg{M})$ under the homomorphism $[0,1] \to [k, k+1]$ mapping $x$ to $x + k$.
Therefore, from \cref{l:dense-or-discrete} we deduce the following fact, which is known for MV-algebras \cite[Proposition~3.5.3]{CignoliDOttavianoEtAl2000}.

\begin{lemma}\label{l:dense-ord-dicrete-unital}
    Let $\alg{A}$ be a subalgebra of the MV-monoid $[0,1]^+$.
    Either $A$ is dense in $[0,1]$ or there is $n \in \N \setminus \{0\}$ such that $\alg A = \alg \L_n^+$.
\end{lemma}

We have the following result, whose version for MV-algebras is also true \cite[Proposition 8.1.1]{CignoliDOttavianoEtAl2000}.
\begin{lemma} \label{l:generates-as-much-as-01}
    Let $\alg{A}$ be an infinite subalgebra of the MV-monoid $[0,1]^+$.
    Then $\HH\SU\PP(\alg{A}) = \HH\SU \PP([0,1]^+)$.
\end{lemma}

\begin{proof}
    The left-to-right inclusion is trivial since $\alg{A}$ is a subalgebra of $[0,1]^+$.
    For the converse inclusion, we prove that $[0,1]^+ \in \HH\SU\PP(\alg{A})$.
    Let $\alg{S}$ be the subalgebra of the $\omega$-power $\alg{N}^\omega$ consisting of the converging sequences.
    Let $\varphi \colon \alg{S} \to [0,1]^+$ be the function that maps an element of $S$ to its limit as a sequence.
    This is a well-defined homomorphism.
    Since there is no $n \in \mathbb N \setminus \{0\}$ such that $\alg{A} = \alg{\L}_n^+$, $A$ is dense in $[0,1]$ by \cref{l:dense-ord-dicrete-unital} and thus $\varphi$ is surjective.
    Therefore, $[0,1]^+$ is a homomorphic image of $\alg{S}$.
\end{proof}

\begin{theorem} \label{t:char-almost-minimal}
    The almost minimal varieties in $\mathsf{MVM}$ are precisely $\VV(\alg{C}^\Delta_2)$, $\VV(\alg{C}^\nabla_2)$ and $\VV(\alg{\L}_p^+)$ (for $p$ prime), and they are all pairwise distinct.
\end{theorem}
\begin{proof}
    First of all we prove that $\VV(\alg{C}^\Delta_2)$, $\VV(\alg{C}^\nabla_2)$ and $\VV(\alg{\L}_p^+)$ (for $p$ prime)
    are all pairwise distinct and distinct from the atom $\VV(\alg{\L}_1^+)$ (i.e., the variety of bounded distributive lattices).
    Consider the equations
    \begin{equation} \label{eq:idem-oplus}
        x \oplus x \approx x
    \end{equation}
    and
    \begin{equation} \label{eq:idem-odot}
        x \odot x \approx x.
    \end{equation}
    We can see that:
    \begin{enumerate}
        \item  $\alg{\L}_1^+$ satisfies both \eqref{eq:idem-oplus} and \eqref{eq:idem-odot};
        \item $\alg{\L}_p$ ($p$ prime) satisfies neither \eqref{eq:idem-oplus} nor \eqref{eq:idem-odot};
        \item $\alg{C}^\Delta_2$ satisfies \eqref{eq:idem-oplus} but not \eqref{eq:idem-odot};
        \item $\alg{C}^\nabla_2$ satisfies \eqref{eq:idem-odot} but not \eqref{eq:idem-oplus}.
    \end{enumerate}
    Therefore, $\VV(\alg{C}^\Delta_2)$, $\VV(\alg{C}^\nabla_2)$ and $\VV(\alg{\L}_p^+)$ (for $p$ prime) are not minimal varieties.
    Moreover, $\VV(\alg{C}^\Delta_2)$ and $\VV(\alg{C}^\nabla_2)$ are distinct, and both of them are distinct from any $\VV(\alg{\L}_p^+)$ (for $p$ prime).
    For distinct primes $p < q$, $\VV(\alg{\L}_p^+)$ and $\VV(\alg{\L}_q^+)$ are distinct because $(p+1)x \app px$  is satisfied by $\alg{\L}_p^+$ but not by $\alg{\L}_q^+$.
    This proves that $\VV(\alg{C}^\Delta_2)$, $\VV(\alg{C}^\nabla_2)$ and $\VV(\alg{\L}_p^+)$ (for $p$ prime)
    are all pairwise distinct and distinct from the atom $\VV(\alg{\L}_1^+)$.

    Let $\vv V$ be a variety of MV-monoids such that $\VV(\alg{\L}_1^+) \subsetneq \vv V$.
    We prove that $\vv V$ is above one of the varieties in the statement.
    By Birkhoff's Subdirect Representation Theorem \cite[Theorem 8.6]{BurrisSanka}, $\vv V$ is generated by its subdirectly irreducible members.
    Therefore, from $\VV(\alg{\L}_1^+) \subsetneq \vv V$ we deduce that there is a subdirectly irreducible $\alg A \in \vv V \setminus \VV(\alg{\L}_1^+)$.
    Since $\alg A \notin \VV(\alg{\L}_1^+)$, $\alg A$ has at least 3 elements.
    Moreover, by \cref{subdMVM}, $\alg{A}$ is totally ordered and such that, for all $x,y \in A$, $x \oplus y = 1$ or $x \odot y = 0$.

    Hence, by \cref{c:almost-minimal-MV}, one of the following conditions holds: (i) the image of the unique homomorphism from $\alg{A}$ to $[0,1]^+$ is not $\{0,1\}$, (ii) $\alg{C}^\Delta_2 \in \HH\SU(\alg A)$, (iii) $\alg{C}^\nabla_2 \in \HH\SU(\alg A)$.

    If (ii) holds then $\VV(\alg{C}^\Delta_2) \sse \vv V$. Similarly, if (iii) holds then $\VV(\alg{C}^\nabla_2) \sse \vv V$.

    We are left with the case (i).
    In this case, $\HH(\alg A)$ contains a subalgebra $\alg{B}$ of $[0,1]^+$ with at least 3 elements. By \cref{l:dense-ord-dicrete-unital}, either
    $\alg {B} =\alg {\L}_n^+$ for some $n \in \N \setminus\{0,1\}$ (where we have excluded $1$ because $B$ has at least three elements) or $B$ is dense in $[0,1]$.
    In the first case, picking a prime $p$ that divides $n$ (which exists because $n \neq 1$), we get $\alg{\L}_p^+ \in \SU(\alg{B}) \subseteq \HH\SU\PP(\alg A) \subseteq \VV$ and thus $\VV(\alg {\L}^+_p) \sse \vv V$.
    In the second case, by \cref{l:generates-as-much-as-01} we have $\HH\SU\PP(\alg{B}) = \HH\SU\PP([0,1]^+)$, and thus for any prime $p$ we get $\alg{\L}_p^+ \in \SU([0,1]^+) \subseteq \HH\SU\PP([0,1]) = \HH\SU\PP(\alg{B}) \subseteq \HH\SU\PP(\alg{A}) \subseteq \VV$, which implies $\VV(\alg {\L}^+_p) \sse \vv V$.
\end{proof}
This proves that Figure~\ref{f:almost-minimal} above is indeed a faithful representation of the bottom of the lattice $\Lambda(\MVM)$ of all varieties of MV-monoids.

\section{Above the almost minimal varieties of MV-monoids}\label{s:almostmin}

In this section we study the subposet of $\Lambda(\mathsf{MVM})$ given by the subvarieties of MV-monoids with ``small'' generators, improving our picture of the bottom part of $\Lambda(\mathsf{MVM})$.
Indeed, the inspection of the subvarieties generated by ``small'' generators is often a useful tool in investigating a given variety.
The main result of this section is a description of the sublattice of $\Lambda(\mathsf{MVM})$ consisting of all varieties generated by an arbitrary set of MV-monoids with at most four elements (\cref{t:with-at-most-4}).

It is clear that $\alg {\L}_1^+$ is a subalgebra of any MV-monoid, and so $\VV(\alg {\L}_1^+)$ is the only atom in $\Lambda(\mathsf{MVM})$.
Its covers are described in \cref{t:char-almost-minimal}.

Above $\VV(\alg {\L}_1^+)$ we can observe that we have two countable chains
$\{\VV(\alg{C}^\Delta_n) \mid n \in \N \setminus \{0\}\}$ and $\{\VV(\alg{C}^\nabla_n) \mid n \in \N \setminus \{0\}\}$.
The inclusion $\VV(\alg{C}_n^\Delta) \subseteq \VV(\alg{C}_{n+1}^\Delta)$ is proper, since $nx \approx (n-1)x$ holds in $\alg{C}_n^\Delta$ but not in $\alg{C}_{n+1}^\Delta$.
Similarly, the inclusion $\VV(\alg{C}_n^\nabla) \subseteq \VV(\alg{C}_{n+1}^\nabla)$ is proper, since the equation $x^n \approx x^{n-1}$ holds in $\alg{C}_n^\nabla$ but not in $\alg{C}_{n+1}^\nabla$.
We show that the joins of the chains $\{\VV(\alg{C}^\Delta_n) \mid n \in \N \setminus \{0\}\}$ and $\{\VV(\alg{C}^\nabla_n) \mid n \in \N \setminus \{0\}\}$ are $\VV(\alg{C}^\Delta)$ and $\VV(\alg{C}^\nabla)$, respectively.
It is clear that $\VV(\alg{C}^\Delta)$ is above every $\VV(\alg{C}_n^\Delta)$, since $\alg{C}_n^\Delta$ is a quotient of $\alg{C}^\Delta$.
Then, to show that $\VV(\alg{C}^\Delta)$ is the least upper bound of $\{\VV(\alg{C}^\Delta_n) \mid n \in \N \setminus \{0\}\}$ it is enough to show that $\alg{C}^\Delta \in \II\SU\PP(\{\alg{C}^\Delta_{n} \mid n \in \N \setminus \{0\}\})$; for this it is enough to observe that the homomorphism
\begin{align*}
    C^\Delta & \longrightarrow \prod_{n \in \N \setminus \{0\}}C^\Delta_n\\
    x & \longmapsto (x/\gamma_{n-1})_{n \in \N \setminus \{0\}},
\end{align*}
is injective since $\bigcap_{n \in \N\setminus \{0\}}\gamma_{n-1}$ is the trivial congruence (cf.\ \cref{cor:Cdeltacon}).
\begin{lemma}
    $\alg{C}^\Delta \in \II\SU\PP_u(\{\alg{C}^\Delta_{n} \mid n \in \N \setminus \{0\}\})$.
\end{lemma}
\begin{proof}
    Let $\mathcal{U}$ be an ultrafilter that extends the filter of cofinite subsets of $\N \setminus \{0\}$, which exists by the Ultrafilter Lemma.
    The homomorphism
    \begin{align*}
        C^\Delta & \longrightarrow \mleft(\prod_{n \in \N \setminus \{0\}}C^\Delta_n\mright)/\mathcal{U}\\
        x & \longmapsto ((x/\gamma_{n-1})_{n \in \N \setminus \{0\}})/\mathcal{U}
    \end{align*}
    is injective because, for all distinct $x,y \in C^\Delta$, the set $\{n \in \N \setminus \{0\} \mid x/\gamma_{n - 1} = y/ \gamma_{n-1}\}$ is finite and hence not in $\mathcal{U}$.
\end{proof}

\color{black}

Therefore, $\bigvee_{n \in \N \setminus \{0\}} \VV(\alg{C}^\Delta_n) = \VV(\alg{C}^\Delta)$; analogously, $\bigvee_{n \in \N \setminus \{0\}} \VV(\alg{C}_n^\nabla) = \VV(\alg{C}^\nabla)$.

\begin{figure}[h]
    \begin{center}
        \begin{tikzpicture}[scale=1.3]
            \draw (0,0) -- (0,1);
            \draw (0,1) -- (0.75,2) -- (0,1) -- (1.75,2) -- (0,1) -- (2.75,2);
            \draw (0,1) -- (-.75,2) -- (0,1) -- (-1.75,2);
            \draw[fill] (0,0) circle [radius=0.05];
            \draw[fill] (0,1) circle [radius=0.05];
            \draw[fill] (0.75,2) circle [radius=0.05];
            \draw[fill] (1.75,2) circle [radius=0.05];
            \draw[fill] (2.75,2) circle [radius=0.05];
            \draw[fill] (-.75,2) circle [radius=0.05];
            \draw[fill] (-1.75,2) circle [radius=0.05];
            \draw (-1.75,2) -- (-1.75,3.4);
            \draw[dashed] (-1.75,3.4) -- (-1.75,4.4);
            \draw (-.75,2) -- (-.75,3.4);
            \draw[dashed] (-.75,3.4) -- (-.75,4.4);
            \draw[dashed] (-1.75,4.4) -- (-1.25,5.1) -- (-.75,4.4);
            \draw[fill] (-1.75,2.7) circle [radius=0.05];
            \draw[fill] (-1.75,3.4) circle [radius=0.05];
            \draw[fill] (-.75,2.7) circle [radius=0.05];
            \draw[fill] (-.75,3.4) circle [radius=0.05];
            \draw[fill] (-1.75,4.4) circle [radius=0.05];
            \draw[fill] (-.75,4.4) circle [radius=0.05];
            \node at (3.85,2){\footnotesize $\cdots$};
            \node[right] at (0,1){\ \ \tiny $\VV(\alg{\L}_1^+)$};
            \node[right] at (-.75,2){\tiny $\VV(\alg{C}^\nabla_2)$};
            \node[left] at (-1.75,2){\tiny $\VV(\alg{C}^\Delta_2)$};
            \node[right] at (.75,2){\tiny $\VV(\alg {\L}_2^+)$};
            \node[right] at (1.75,2){\tiny $\VV(\alg {\L}_3^+)$};
            \node[right] at (2.75,2){\tiny $\VV(\alg {\L}_5^+)$};
            \node[left] at (-1.75,2.7){\tiny $\VV(\alg{C}^\Delta_3)$};
            \node[right] at (-.75,2.7){\tiny $\VV(\alg{C}^\nabla_3)$};
            \node[left] at (-1.75,4.4){\tiny $\VV(\alg{C}^\Delta)$};
            \node[right] at (-.75,4.4){\tiny $\VV(\alg{C}^\nabla)$};
            \node[above] at (-1.25,5.1){\tiny $\VV(\alg C^+)$};
        \end{tikzpicture}
        \caption{A portion of $\Lambda(\mathsf{MVM})$}\label{MVMlattice}
    \end{center}
\end{figure}

We now start looking at ``small'' MV-monoids and the varieties generated by them. It is obvious that the only 2-element MV-monoid is the 2-element bounded distributive lattice. There are exactly four 3-element MV-monoids, represented in Figure~\ref{3elements};
the right-most one is the 3-element bounded chain, which clearly does not satisfy the conclusion of \cref{subdMVM} and hence is not subdirectly irreducible.
The remaining three ones are all subdirectly irreducible.
In Figures~\ref{3elements} and~\ref{4elements} we use the convention of specifying the outputs of $\oplus$ and $\odot$ when they are not derivable from other outputs, i.e., when they cannot be inferred by the monotonicity of $\oplus$ and $\odot$ and by the properties of $1$ and $0$ with respect to $\oplus$ and $\odot$.

\begin{figure}[h]
    \begin{center}
        \begin{tikzpicture}[scale=.7]
            \draw (0,0) -- (0,1) -- (0,2);
            \draw (4,0) -- (4,1) -- (4,2);
            \draw  (8,0) --(8,1) -- (8,2);
            \draw  (12,0) --(12,1) -- (12,2);
            \draw[fill] (0,0) circle [radius=0.05];
            \draw[fill] (0,1) circle [radius=0.05];
            \draw[fill] (0,2) circle [radius=0.05];
            \draw[fill] (4,0) circle [radius=0.05];
            \draw[fill] (4,1) circle [radius=0.05];
            \draw[fill] (4,2) circle [radius=0.05];
            \draw[fill] (8,0) circle [radius=0.05];
            \draw[fill] (8,1) circle [radius=0.05];
            \draw[fill] (8,2) circle [radius=0.05];
            \draw[fill] (12,0) circle [radius=0.05];
            \draw[fill] (12,1) circle [radius=0.05];
            \draw[fill] (12,2) circle [radius=0.05];
            \node[right] at (0,0){\footnotesize $0=a \odot a$};
            \node[right] at(0,1){\footnotesize $a$};
            \node[right] at (0,2){\footnotesize $1=a \oplus a$};
            \node[right] at (4,0){\footnotesize $0 =  a \odot a$};
            \node[right] at (4,1){\footnotesize $a =a \oplus a$};
            \node[right] at (4,2){\footnotesize $1$};
            \node[right] at (8,0){\footnotesize $0$};
            \node[right] at (8,1){\footnotesize $a =a \odot a$};
            \node[right] at (8,2){\footnotesize $1 = a \oplus a$};
            \node[below] at (0,-.5) {\footnotesize $\alg {\L}_{2}^+$};
            \node[below] at (4,-.5) {\footnotesize $ \alg{C}^\Delta_2$};
            \node[below] at (8,-.5) {\footnotesize $ \alg{C}^\nabla_2$};
            \node[right] at (12,0){\footnotesize $0$};
            \node[right] at (12,1){\footnotesize $a =a \odot a = a \oplus a$};
            \node[right] at (12,2){\footnotesize $1$};
            \node[below] at (12,-.5) {\footnotesize $ \alg{L}_2$};

        \end{tikzpicture}
        \caption{the 3-element MV-monoids}\label{3elements}
    \end{center}
\end{figure}
Furthermore, for each $n \in \N$ we define the MV-monoids $\alg{LM}_n^\Delta$ and $\alg{LM}_n^\nabla$ whose lattice reduct is the $(n+1)$-element chain and whose monoidal operations are defined as follows:
\begin{itemize}
    \item $\alg{LM}_n^\Delta = \Gamma(\Z \overrightarrow{\times} \alg{LM}_n^{\Delta*})$  is the MV-monoid in which $x \oplus y = x \join y$ and
    $$
    x \odot y =
    \begin{cases}
        y & \hbox{if $x=1$,} \\
        0 & \hbox{otherwise;}
    \end{cases}
    $$
    \item $\alg{LM}_n^\nabla = \Gamma(\Z \overrightarrow{\times} \alg{LM}_n^{\nabla*})$ is the MV-monoid in which $x \odot y = x \meet y$ and
    $$
    x \oplus y =
    \begin{cases}
        y & \hbox{if $x=0$,} \\
        1 & \hbox{otherwise.}
    \end{cases}
    $$
\end{itemize}
We can observe that $\alg{LM}^\Delta_1 \cong \alg{LM}^\nabla_1 \cong \alg{\L}_1^+$, $\alg{LM}^\Delta_2 \cong  \alg{C}^\Delta_2$, and $\alg{LM}^\nabla_2 \cong \alg{C}^\nabla_2$.
Thus, the new cases are for $n \geq 3$.
In fact, for $n \ge 3$, $\alg {\L}^+_n$, $\alg{C}^\Delta_n$, $\alg{C}^\nabla_n$, $\alg{LM}^\Delta_n$ and $\alg{LM}^\nabla_n$ are pairwise nonisomorphic.
Moreover, for all $n \geq 2$, $\alg{LM}^\Delta_n$ and $\alg{LM}^\nabla_n$ are not positive MV-algebras (i.e., they are not cancellative).
Finally, we note that for all $n \in \N \setminus\{0\}$, both $\alg{LM}^\Delta_n$ and $\alg{LM}^\nabla_n$ satisfy the conclusion of \cref{subdMVM} (i.e., the necessary condition for subdirect irreducibility). However:

\begin{lemma}\label{l:LMn-notsubirr} If $n \ge 3$, neither $\alg{LM}^\Delta_n$ nor $\alg{LM}^\nabla_n$ are subdirectly irreducible.
\end{lemma}
\begin{proof}
    Suppose $n \ge 3$; then there are $a,b \in LM^\Delta_n$ such that $0 \cov a \cov b < 1$.
    We claim that $\{a,b\}$ and $\{0,a\}$ identify two congruences (meaning that they are the only nontrivial blocks) whose meet is the bottom congruence.
    To prove that $\{a,b\}$ is the only nontrivial block of a congruence we have to show that, for all $x \in LM^\Delta_n$,
    \begin{align*}
        &\text{either $a \oplus x = b \oplus x$ or $\{a \oplus x, b\oplus x\} = \{a,b\}$, and}\\
        &\text{either $a \odot x = b \odot x$ or $\{a \odot x, b\odot x\} = \{a,b\}$}.
    \end{align*}
    The proof is really straightforward given the definition of the operations. The argument for $\{0,a\}$ is similar, but easier.
    It is immediate that the meet of these two congruences is the bottom congruence.

    Finally, the argument for $\alg{LM}^\nabla_n$ is analogous.
\end{proof}
This shows that even a finite MV-monoid satisfying the conclusion of \cref{subdMVM} need not be subdirectly irreducible.

Using the softwares \emph{Prover-9} and \emph{Mace4} \cite{prover9-mace4}, we can observe that there are 19 totally ordered 4-element MV-monoids; only nine of them satisfy the conclusion of \cref{subdMVM}.

\begin{figure}[h]
    \begin{center}
        \begin{tikzpicture}[scale=.7]
        \draw (0,12) -- (0,13) -- (0,14) --(0,15);
            \draw (4,12) -- (4,13) -- (4,14) -- (4,15);
            \draw  (8,12) --(8,13) -- (8,14) -- (8,15);
            \draw (12,12) --(12,13) -- (12,14) -- (12,15);
            \draw[fill] (0,12) circle [radius=0.05];
            \draw[fill] (0,13) circle [radius=0.05];
            \draw[fill] (0,14) circle [radius=0.05];
            \draw[fill] (0,15) circle [radius=0.05];
            \draw[fill] (4,12) circle [radius=0.05];
            \draw[fill] (4,13) circle [radius=0.05];
            \draw[fill] (4,14) circle [radius=0.05];
            \draw[fill] (4,15) circle [radius=0.05];
            \draw[fill] (8,12) circle [radius=0.05];
            \draw[fill] (8,13) circle [radius=0.05];
            \draw[fill] (8,14) circle [radius=0.05];
            \draw[fill] (8,15) circle [radius=0.05];
            \draw[fill] (12,12) circle [radius=0.05];
            \draw[fill] (12,13) circle [radius=0.05];
            \draw[fill] (12,14) circle [radius=0.05];
            \draw[fill] (12,15) circle [radius=0.05];

            \node[below] at (0,11.5) {\footnotesize $\alg A^{\Delta}_3$};
            \node[below] at (4,11.5) {\footnotesize $\alg A^{\nabla}_3$};
            \node[below] at (8,11.5) {\footnotesize $\alg B^{\Delta}_3$};
            \node[below] at (12,11.5) {\footnotesize $\alg B^{\nabla}_3$};

            \node[right] at (0,12){\footnotesize $0 = b \odot b$};
            \node[right] at(0,13){\footnotesize $b =b \oplus b= a \odot b$};
            \node[right] at (0,14){\footnotesize $a = a \odot a$ };
            \node[right] at (0,15){\footnotesize $1 = a \oplus b$};

            \node[right] at (4,12){\footnotesize $0 = a \odot b$};
            \node[right] at(4,13){\footnotesize $b =b \oplus b$};
            \node[right] at (4,14){\footnotesize $a = a\oplus b= a \odot a$};
            \node[right] at (4,15){\footnotesize $1 = a \oplus a$};
ì
            \node[right] at (8,12){\footnotesize $0 = a \odot b$};
            \node[right] at(8,13){\footnotesize $b=b \oplus b= a \odot a$};
            \node[right] at (8,14){\footnotesize $a=a\oplus b$};
            \node[right] at (8,15){\footnotesize $1 = a \oplus a$};
            \node[right] at (12,12){\footnotesize $0 = b \odot b$};
            \node[right] at(12,13){\footnotesize $b =a\odot b$};
            \node[right] at (12,14){\footnotesize $a=b \oplus b = a \odot a$};
            \node[right] at (12,15){\footnotesize $1 = a \oplus b$};

            \draw (0,6) -- (0,7) -- (0,8) --(0,9);
            \draw (4,6) -- (4,7) -- (4,8) -- (4,9);
            \draw (8,6) --(8,7) -- (8,8) -- (8,9);
            \draw[fill] (0,6) circle [radius=0.05];
            \draw[fill] (0,7) circle [radius=0.05];
            \draw[fill] (0,8) circle [radius=0.05];
            \draw[fill] (0,9) circle [radius=0.05];
            \draw[fill] (4,6) circle [radius=0.05];
            \draw[fill] (4,7) circle [radius=0.05];
            \draw[fill] (4,8) circle [radius=0.05];
            \draw[fill] (4,9) circle [radius=0.05];
            \draw[fill] (8,6) circle [radius=0.05];
            \draw[fill] (8,7) circle [radius=0.05];
            \draw[fill] (8,8) circle [radius=0.05];
            \draw[fill] (8,9) circle [radius=0.05];

            \node[below] at (0,5.5) {\footnotesize $\alg {\L}_{3}^+$};
            \node[below] at (8,5.5) {\footnotesize $\alg{C}^\Delta_3$};
            \node[below] at (4,5.5) {\footnotesize $\alg{C}^\nabla_3$};

            \node[right] at (0,6){\footnotesize $0 =a \odot b$};
            \node[right] at(0,7){\footnotesize $b =a \odot a$};
            \node[right] at (0,8){\footnotesize $a = b\oplus b$};
            \node[right] at (0,9){\footnotesize $1=a\oplus b$};

            \node[right] at (8,6){\footnotesize $0 = a \odot a$};
            \node[right] at(8,7){\footnotesize $b$};
            \node[right] at (8,8){\footnotesize $a=a \oplus b=b \oplus b=a \oplus a$};
            \node[right] at (8,9){\footnotesize $1$};

            \node[right] at (4,6){\footnotesize $0$};
            \node[right] at(4,7){\footnotesize $b =a\odot a = b \odot b$};
            \node[right] at (4,8){\footnotesize $a$};
            \node[right] at (4,9){\footnotesize $1 = b \oplus b$};

            \draw (0,0) -- (0,1) -- (0,2) --(0,3);
            \draw (4,0) -- (4,3);
            \draw[fill] (0,0) circle [radius=0.05];
            \draw[fill] (0,1) circle [radius=0.05];
            \draw[fill] (0,2) circle [radius=0.05];
            \draw[fill] (0,3) circle [radius=0.05];
            \draw[fill] (4,0) circle [radius=0.05];
            \draw[fill] (4,1) circle [radius=0.05];
            \draw[fill] (4,2) circle [radius=0.05];
            \draw[fill] (4,3) circle [radius=0.05];

            \node[below] at (0,-0.5) {\footnotesize $\alg{LM}^\Delta_3$};
            \node[below] at (4,-0.5) {\footnotesize $\alg{LM}^\nabla_3$};

            \node[right] at (0,0){\footnotesize $0 = a \odot a$};
            \node[right] at (0,1){\footnotesize $b=b \oplus b$};
            \node[right] at (0,2){\footnotesize $a= a \oplus a =a \oplus b$};
            \node[right] at (0,3){\footnotesize $1$};
            \node[right] at (4,0){\footnotesize $0$};
            \node[right] at(4,1){\footnotesize $b= a \odot b = b \odot b$};
            \node[right] at (4,2){\footnotesize $a=a \odot a$};
            \node[right] at (4,3){\footnotesize $1 = b \oplus b$};
        \end{tikzpicture}
        \caption{the 4-element MV-monoids satisfying the conclusion of \cref{subdMVM}}\label{4elements}
    \end{center}
\end{figure}
We already know that $\alg {\L}_3^+$ is simple (\cref{lem:LnSimple}) and that $\alg{C}^\Delta_3$ and $\alg{C}^\nabla_3$ are subdirectly irreducible (\cref{l:Bn-subirr}). On the other hand, neither $\alg{LM}^\Delta_3$ nor $\alg{LM}^\nabla_3$ are subdirectly irreducible (\cref{l:LMn-notsubirr}).
What about the others? As they are all 4-element chains, the congruence lattice of their lattice reduct is the $2^3=8$-element distributive complemented lattice, depicted in Figure~\ref{conL}.
In the sequel we will denote a congruence by its nontrivial blocks.
\begin{figure}[h]
    \begin{center}
        \begin{tikzpicture}[scale=1.4]

            \draw (0,0) -- (-1,1) -- (0,2) --(1,1) -- (0,0) -- (0,1)--(-1,2) -- (0,3) -- (1,2) -- (0,1);
            \draw (-1,1)--(-1,2);
            \draw (1,1) --(1,2);
            \draw (0,2) -- (0,3);
            \draw[fill] (0,0) circle [radius=0.05];
            \draw[fill] (-1,1) circle [radius=0.05];
            \draw[fill] (0,2) circle [radius=0.05];
            \draw[fill] (1,1) circle [radius=0.05];
            \draw[fill] (0,1) circle [radius=0.05];
            \draw[fill] (-1,2) circle [radius=0.05];
            \draw[fill] (0,3) circle [radius=0.05];
            \draw[fill] (1,2) circle [radius=0.05];
            \node[left] at (-1,1){\footnotesize $\{0,b\}$};
            \node[right] at (0,1){\footnotesize $\{a,b\}$};
            \node[right] at (1,1){\footnotesize $\{a,1\}$};
            \node[left] at (-1,2){\footnotesize $\{0,a,b\}$};
            \node[left] at (0,2){\footnotesize $\{0,b\}$};
            \node[right] at (0,2) {\footnotesize $\{a,1\}$};
            \node[right] at (1,2) {\footnotesize $\{a,b,1\}$};
        \end{tikzpicture}
        \caption{$\con(\alg L_3)$}\label{conL}
    \end{center}
\end{figure}

\begin{enumerate}
    \ib We can see that $\{0,b\}$ is a congruence of $\alg A^{\nabla}_3$. Moreover, if $\th \in \op{Con}(\alg A^{\nabla}_3)$ and $(a,b) \in \th$, then
    \begin{align*}
        &(1,a) = (a \oplus a,b \oplus a) \in \th\\
        &(a,0) = (a \odot a, b \odot a) \in \th.
    \end{align*}
    Thus $(0,1) \in \th$, and therefore $\th$ is the total congruence. If instead $(a,1) \in \th$, then
    $$
    (0,b) = (a \odot b,1 \odot b) \in \th.
    $$
    It follows that $\op{Con}(\alg A^{\nabla}_3)$ is the 4-element chain
    $$
    0_{\alg A^{\nabla}_3} < \{0,b\} < \{\{0,b\}\{1,a\}\} < 1_{\alg A^{\nabla}_3}
    $$
    and hence $\alg A^{\nabla}_3$ is subdirectly irreducible.
    Similar computations show that $\op{Con}(\alg A^{\Delta}_3)$ is the 4-element chain with minimal congruence $\{a,1\}$.
    So $\alg{A}^{\nabla}_3$ is subdirectly irreducible.

    \ib For $\alg B^{\Delta}_3$ a similar computation proves that $\op{Con}(\alg B^{\Delta}_3)$ is the 3-element chain and the minimal congruence is again $\{0,b\}$. So $\alg B^{\Delta}_3$ is subdirectly irreducible as well.
    Similar computations show that $\op{Con}(\alg B^{\nabla}_3)$ is the 3-element chain with minimal congruence $\{a,1\}$. So $\alg{B}^{\Delta}_3$ is subdirectly irreducible.
\end{enumerate}

If we want the 4-element chains to come into play, then we have to zoom in; first, we observe that using the description of the congruences we can show that
\begin{align*}
    &\HH\SU(\alg A^{\nabla}_3) = \II(\{*\}, \alg {\L}_1^+, \alg{C}^\Delta_2,\alg{C}^\nabla_2, \alg A^{\nabla}_3),\\
    &\HH\SU(\alg B^{\Delta}_3) = \II(\{*\},\alg {\L}_1^+, \alg{C}^\Delta_2, \alg {\L}_2^+,\alg B^{\Delta}_3),\\
    &\HH\SU(\alg B^{\nabla}_3) = \II(\{*\},\alg {\L}_1^+, \alg{C}^\nabla_2, \alg {\L}_2^+, \alg B^{\nabla}_3),\\
    &\HH\SU(\alg A^{\Delta}_3) = \II(\{*\}, \alg {\L}_1^+, \alg{C}^\Delta_2,\alg{C}^\nabla_2, \alg A^{\Delta}_3),
\end{align*}
where $\{*\}$ is the trivial algebra.
Now we can draw the picture in Figure~\ref{firstzoomin}, which describes the inclusion relations of the other 4-element chains.

\begin{figure}[h]
    \begin{center}
        \begin{tikzpicture}[scale=1.5]
            \draw (0,.6) -- (0,1.3) -- (-1,2) -- (-1.5,2.7) -- (0,2) -- (0,1.3);
            \draw (0,1.3) -- (1,2) -- (1.5,2.7) -- (0,2);
            \draw (-1,2) -- (0,2.7) -- (1,2);
            \draw (-1.5,2.7) -- (-1.5,3.4);
            \draw (.5,3.4) -- (0,2.7);
            \draw (1.5,2.7) -- (1.5,3.4);
            \draw (-.5,3.4) -- (0,2.7);
            \draw[fill] (0,.6) circle [radius=0.05];
            \draw[fill] (0,1.3) circle [radius=0.05];
            \draw[fill] (-1,2) circle [radius=0.05];
            \draw[fill] (0,2) circle [radius=0.05];
            \draw[fill] (1,2) circle [radius=0.05];
            \draw[fill] (-1.5,2.7) circle [radius=0.05];
            \draw[fill] (0,2.7) circle [radius=0.05];
            \draw[fill] (1.5,2.7) circle [radius=0.05];
            \draw[fill] (-1.5,3.4) circle [radius=0.05];
            \draw[fill] (-.5,3.4) circle [radius=0.05];
            \draw[fill] (.5,3.4) circle [radius=0.05];
            \draw[fill] (1.5,3.4) circle [radius=0.05];
            \node[right] at (0,1.3){\tiny $\VV(\alg{\L}_1^+)$};
            \node[left] at (-1,2){\tiny $\VV(\alg{C}^\Delta_2$)};
            \node[right] at (0,2){\tiny $\VV(\alg {\L}_2^+)$};
            \node[right] at (1,2){\tiny $\VV(\alg{C}^\nabla_2)$};
            \node[left] at (-1.5,2.7){\tiny $\VV(\alg{C}^\Delta_2,\alg {\L}^+_2)$};
            \node[right] at (0,2.7){\tiny $\VV(\alg{C}^\Delta_2,\alg{C}^\nabla_2)$};
            \node[right] at (1.5,2.7){\tiny $\VV(\alg {\L}_2^+,\alg{C}^\nabla_2)$};
            \node[above] at (-1.5,3.4){\tiny $\VV(\alg B^{\Delta}_3)$};
            \node[above] at (-.5,3.4){\tiny $\VV(\alg A^{\nabla}_3)$};
            \node[above] at (.5,3.4){\tiny $\VV(\alg A^{\Delta}_3)$};
            \node[above] at (1.5,3.4){\tiny $\VV(\alg B^{\nabla}_3)$};
        \end{tikzpicture}
        \caption{The first zoom in}\label{firstzoomin}
    \end{center}
\end{figure}

However, this picture can be improved. Indeed, if we zoom in more we can see that the interval
$[\vv T, \VV(\alg A^{\Delta}_3, \alg A^{\nabla}_3, \alg B^{\Delta}_3,\alg B^{\nabla}_3)]$ ($\vv T$ is the trivial variety) consists of three isomorphic posets glued together. In Figure~\ref{secondzoomin} we display one of the posets (the labels for the unlabeled points are obvious).

In particular, by an easy application of J\'onsson's Lemma \cite[Corollary 3.4]{Jonsson1967}, we have that every variety generated by a set of MV-monoids with cardinality less or equal four is generated by a subset of the subdirectly irreducible MV-monoids with at most four elements $I = \{\alg \L_1^+, \alg \L_2^+, \alg \L_3^+, \alg{C}^\Delta_2, \alg{C}^\nabla_2, \alg{C}^\Delta_3, \alg{C}^\nabla_3, \alg A^{\Delta}_3, \alg A^{\nabla}_3, \alg B^{\Delta}_3, \alg B^{\nabla}_3\}$. The poset formed by $\HH\SU$ of the elements in $I$ is depicted in Figure~\ref{posetsmall}.

\begin{figure}[h]
    \begin{center}
        \begin{tikzpicture}[scale=1.5]
            \draw (0,0) -- (-1.5,1.5) -- (-2.5,3) -- (-1.5,1.5)-- (-1.5,3) -- (-1.5,1.5) -- (-0.5,3) -- (-1.5,1.5) -- (0.5,3) -- (-1.5,1.5) -- (0,0) -- (0,1.5) -- (-1.5,3) -- (0,1.5) -- (1.5,3) -- (1.5,1.5) -- (2.5,3) -- (1.5,1.5) -- (0.5,3) -- (1.5,1.5) -- (-0.5,3) -- (1.5,1.5) -- (0,0) -- (3.5,3);

            \draw[fill] (0,0) circle [radius=0.05];
             \node[left] at (0,0){\tiny $\alg \L_1^+$};

            \draw[fill] (1.5,1.5) circle [radius=0.05];
            \node[left] at (1.5,1.5){\tiny $\alg C_2^{\nabla}$};
            \draw[fill] (-1.5,1.5) circle [radius=0.05];
            \node[left] at (-1.5,1.5){\tiny $\alg C_2^{\Delta}$};
            \draw[fill] (0,1.5) circle [radius=0.05];
            \node[right] at (0,1.5){\tiny $\alg \L_2^+$};
            \draw[fill] (3.5,3) circle [radius=0.05];
            \node[right] at (3.5,3.1){\tiny $\alg \L_3^+$};

            \draw[fill] (1.5,3) circle[radius=0.05];
            \node[right] at (1.5,3.1){\tiny $\alg B^{\nabla}_3$};
            \draw[fill] (-1.5,3) circle[radius=0.05];
            \node[left] at (-1.5,3.1){\tiny $\alg B^{\Delta}_3$};
            \draw[fill] (-2.5,3) circle[radius=0.05];
            \node[left] at (-2.5,3.1){\tiny $\alg C_3^{\Delta}$};
            \draw[fill] (2.5,3) circle[radius=0.05];
            \node[right] at (2.5,3.1){\tiny $\alg C_3^{\nabla}$};
            \draw[fill] (0.5,3) circle[radius=0.05];
            \node[right] at (0.5,3.1){\tiny $\alg A^{\nabla}_3$};
            \draw[fill] (-0.5,3) circle[radius=0.05];
            \node[left] at (-0.5,3.1){\tiny $\alg A^{\Delta}_3$};
\end{tikzpicture}
        \caption{Subdirectly irreducible MV-monoids with cardinality $\leq 4$ ordered by: $\alg{A} \leq \alg{B}$ iff $\alg{A} \in \HH\SU(\alg{B})$}\label{posetsmall}
    \end{center}
\end{figure}

\begin{theorem} \label{t:with-at-most-4}
    The sublattice (in fact, the ideal) of $\Lambda(\mathsf{MVM})$ consisting of all varieties generated by an arbitrary set of MV-monoids with at most four elements is isomorphic to the lattice of downward-closed subsets of the poset in Figure~\ref{posetsmall}.
\end{theorem}

\begin{proof}
    The proof is a direct consequence of J\'onsson's Lemma \cite[Corollary 3.4]{Jonsson1967}. Indeed let $\vv V(\vv K)$ be a variety generated by a set $\vv K$ of MV-monoids with at most 4 elements. By \cite[Corollary 3.4]{Jonsson1967} the set $\vv I$ of subdirectly irreducible algebras in $\vv V(\vv K)$ is contained in $\HH\SU(\vv K)$ and thus it is composed by algebras with at most 4 elements. Moreover, it is clear that $\vv I$ contains all the subdirectly irreducible MV-monoids in $\HH\SU(\vv I)$. We already observed that the subdirectly irreducible MV-monoids of at most 4 elements are $\{\alg \L_1^+, \alg \L_2^+, \alg \L_3^+, \alg{C}^\Delta_2, \alg{C}^\nabla_2, \alg{C}^\Delta_3, \alg{C}^\nabla_3, \alg A^{\Delta}_3, \alg A^{\nabla}_3, \alg B^{\Delta}_3, \alg B^{\nabla}_3\}$.
    The desired statement follows from this and the fact that, as a consequence of Birkhoff's Subdirect Representation Theorem \cite[Theorem 8.6]{BurrisSanka}, a variety $\VV_1$ is contained in a variety $\VV_2$ if and only if each subdirectly irreducible member in $\VV_1$ is in $\VV_2$.
\end{proof}

It is clear that making an understandable drawing of the ideal of all varieties generated by MV-monoids with cardinality less than or equal to $4$ is a difficult task. If we restrict further the bound to three elements this is feasible (see Figure~\ref{fig:3elm-var}).

\begin{figure}[h]
    \begin{center}
        \begin{tikzpicture}[scale=1.5]
            \draw (0,.6) -- (0,2) ;
            \draw (0,1.3) -- (-1.5,2) -- (0,3.4) -- (-.75,4.1) -- (0,3.4) --(1.5,2) -- (0,1.3);
            \draw (.75,4.1) -- (0,3.4);
            \draw (0,2) -- (-1.5,3.4) -- (0,4.8) -- (1.5,3.4) -- (0,2);
            \draw[fill] (0,.6) circle [radius=0.05];
            \draw[fill] (0,1.3) circle [radius=0.05];
            \draw[fill] (-1.5,2) circle [radius=0.05];
            \draw[fill] (0,2) circle [radius=0.05];
            \draw[fill] (1.5,2) circle [radius=0.05];
            \draw[fill] (-.75,2.7) circle [radius=0.05];
            \draw[fill] (.75,2.7) circle [radius=0.05];
            \draw[fill] (-1.5,3.4) circle [radius=0.05];
            \draw[fill] (0,3.4) circle [radius=0.05];
            \draw[fill] (1.5,3.4) circle [radius=0.05];
            \draw[fill] (-.75,4.1) circle [radius=0.05];
            \draw[fill] (.75,4.1) circle [radius=0.05];
            \draw[fill] (0,4.8) circle [radius=0.05];
            \node[right] at (0,1.3){\tiny $\VV(\alg{\L}_1^+)$};
            \node[left] at (-1.5,2){\tiny $\VV(\alg{C}^\nabla_2$)};
            \node[right] at (0,2){\tiny $\VV(\alg{C}^\Delta_2)$};
            \node[right] at (1.5,2){\tiny $\VV(\alg {\L}_2^+)$};
            \node[left] at (-1.5,3.4){\tiny $\VV(\alg A^{\nabla}_3)$};
            \node[left] at (-.75,4.1){\tiny $\VV(\alg A^{\nabla}_3,\alg {\L}_2^+)$};
            \node[right] at (1.5,3.4){\tiny $\VV(\alg B^{\Delta}_3)$};
            \node[right] at (.75,4.1){\tiny $\VV(\alg B^{\Delta}_3,\alg{C}^\nabla_2)$};
            \node[above] at (0,4.8) {\tiny $\VV(\alg A^{\nabla}_3,\alg B^{\Delta}_3)$};
        \end{tikzpicture}
        \caption{The second zoom in}\label{secondzoomin}
    \end{center}
\end{figure}

A similar pattern repeats itself for each cardinality $n \in \N$ but it gets hopelessly complicated very quickly. There are already thirty-five
5-element MV-monoids satisfying the conclusion of \cref{subdMVM}, but we already know that two of them ($\alg{LM}^\Delta_4$ and $\alg{LM}^\nabla_4$) are not subdirectly irreducible and three of them ($\alg {\L}_4^+$, $\alg{C}^\Delta_4$ and $\alg{C}^\nabla_4$) are subdirectly irreducible. The remaining thirty could be checked by hand (or one could implement a script) and in principle we can establish all the covering relations among them. Drawing the proper interval would be a real work of art.

\begin{figure}[h]
    \begin{center}
        \begin{tikzpicture}[scale=1.1]

            \draw (0,-1) -- (0,0) -- (-1.5,1.5) -- (0,3) --(1.5,1.5) -- (0,0) -- (0,1.5)--(-1.5,3) -- (0,4.5) -- (1.5,3) -- (0,1.5);
            \draw (-1.5,1.5)--(-1.5,3);
            \draw (1.5,1.5) --(1.5,3);
            \draw (0,3) -- (0,4.5);
            \draw[fill] (0,-1) circle [radius=0.05];
            \draw[fill] (0,0) circle [radius=0.05];
            \draw[fill] (-1.5,1.5) circle [radius=0.05];
            \draw[fill] (0,3) circle [radius=0.05];
            \draw[fill] (1.5,1.5) circle [radius=0.05];
            \draw[fill] (0,1.5) circle [radius=0.05];
            \draw[fill] (-1.5,3) circle [radius=0.05];
            \draw[fill] (0,4.5) circle [radius=0.05];
            \draw[fill] (1.5,3) circle [radius=0.05];
            \node[right] at (0,-1){\tiny $\vv T$};
            \node[left] at (-1.5,1.5){\tiny $\vv V(\alg C_2^{\Delta})$};
             \node[right] at (0,0){\tiny $\vv V(\alg \L_1^+)$};
            \node[right] at (0,1.5){\tiny $\vv V(\alg \L_2^+)$};
            \node[right] at (1.5,1.5){\tiny $\vv V(\alg C_2^{\nabla})$};
            \node[left] at (-1.5,3){\tiny $\vv V(\alg C_2^{\Delta}, \alg \L_2^+)$};
            \node[right] at (0,3){\tiny $\vv V(\alg C_2^{\Delta}, \alg C_2^{\nabla})$};;
            \node[right] at (1.5,3) {\tiny $\vv V(\alg \L_2^+, \alg C_2^{\nabla})$};
            \node[right] at (0,4.5) {\tiny $\vv V(\alg \L_2^+, \alg C_2^{\Delta}, \alg C_2^{\nabla})$};
        \end{tikzpicture}
        \caption{The lattice of all the varieties of MV-monoids generated by algebras with at most three elements}\label{fig:3elm-var}
    \end{center}
\end{figure}

\section{Subdirectly irreducible positive MV-algebras}\label{sec: si pMV}

In this section we begin our investigation of positive MV-algebras. As is usual for many other classes of algebraic structures, one of the first pivotal steps is to characterize the subdirectly irreducible algebras, when feasible.
The principal aim of this section, achieved in \cref{l:fin-is-sub}, is to provide such a characterization in the finite case.
Furthermore, in the last part of the section we provide two examples of (totally ordered) relatively subdirectly irreducible positive MV-algebras that are not subdirectly irreducible in the absolute sense.

We start recording some basic properties of cancellative commutative $\ell$-monoids. We recall that a commutative monoid $\alg{M}$ is called \emph{torsion-free} if for every $n \in \N \setminus \{0\}$ and every $x,y \in M$ such that $n x = ny$ we have $x = y$.

 \begin{lemma}\label{l:torsion-free-monoid}
     Any cancellative commutative $\ell$-monoid is torsion-free.
 \end{lemma}

 \begin{proof}
    It is a standard fact that any Abelian $\ell$-group is torsion-free \cite[Corollary~1.2.13]{BigardKeimelEtAl1977}.
    Cancellative commutative $\ell$-monoids are precisely the subreducts of Abelian $\ell$-groups (see e.g.\ \cite[Proposition~4.3]{Abbadinietal2022}), and thus they are torsion-free.
 \end{proof}

\begin{lemma} \label{l:cancellativity}
    Let $\alg{M}$ be a cancellative commutative $\ell$-monoid.
    For all $x,x',y,y' \in M$, if $x + y = x' + y'$, $x \leq x'$ and $y \leq y'$, then $x = x'$ and $y = y'$.
\end{lemma}

\begin{proof}
    $\alg{M}$ is a subreduct of an Abelian $\ell$-group (see \cite[Proposition~4.3]{Abbadinietal2022}), and these quasi-equations hold in Abelian $\ell$-groups.
\end{proof}

We then turn our attention to the following consequence of \cref{Holder}.
\begin{proposition} \label{p:finite-Gamma-implies-Ln}
    Let $\alg{M}$ be a nontrivial cancellative totally ordered unital commutative $\ell$-monoid such that $\Gamma(\alg{M})$ is finite.
    There is $n \in \N \setminus \{0\}$ such that $\alg{M} \cong \frac{1}{n} \mathbb Z$.
\end{proposition}

\begin{proof}
    Since $\alg{M}$ is a nontrivial totally ordered unital commutative $\ell$-monoid, by \cref{Holder} there is a unique homomorphism $h \colon \alg{M} \to \R$.
    Since $\Gamma(\alg{M})$ is finite, also $h[\Gamma(\alg{M})]$ is finite. Furthermore, we claim that $\Gamma(h[\alg{M}]) = h[\Gamma(\alg{M})]$.
    Indeed, let $y \in \Gamma(h[\alg{M}])$.
    Then there is $x \in M$ such that $h(x) = y$ and $0 \leq y \leq 1$. Set $x' \coloneqq (x \vee 0) \wedge 1$. We observe that $x' \in \Gamma(\alg{M})$ and $h(x') = (h(x) \vee 0) \wedge 1 = (y \vee 0) \wedge 1 = y$. Thus $y \in h[\Gamma(\alg{M})]$. For the other inclusion, let $y \in h[\Gamma(\alg{M})]$. Then, there is $x \in M$ such that $0 \leq x \leq 1$ and $h(x) = y$. We prove that $0 \leq y \leq 1$ and thus that $y \in \Gamma(h[\alg{M}])$:
    \[
    y = h(x) = h((x \vee 0) \wedge 1) = (h(x) \vee 0) \wedge 1 = (y \vee 0) \wedge 1.
    \]
    Thus, $\Gamma(h[\alg{M}]) = h[\Gamma(\alg{M})]$, and hence $\Gamma(h[\alg{M}])$ is finite.
    Therefore, the subalgebra $h[\alg{M}]$ of $\R$ is not dense and hence, by \cref{l:dense-or-discrete} there is $n \in \N \setminus \{0\}$ such that $h[\alg{M}] = \frac{1}{n} \mathbb Z$. We prove that $h$ is injective and thus that $\alg{M} \cong \frac{1}{n} \mathbb Z$. To do so we claim that, for all $x \in M$ and $n,k \in \mathbb Z$, if $h(nx) = k$, then $nx = k$.
    Let $x \in M$ and $n,k \in \mathbb Z$ and suppose $h(nx) = k$;
    then $h(nx - k) = 0$ and, since $\alg{M}$ is totally ordered, either $nx - k \leq 0$ or $nx - k \geq 0$.
    We only deal with the case $nx - k \geq 0$ since the other one is analogous.
    Let us then assume $nx - k \geq 0$;
    for every $l \in \N \setminus \{0\}$ we have $h(l(nx - k)) = 0$ and $l(nx - k) \geq 0$, from which we conclude $l(nx - k) \in \Gamma(\alg{M})$.
    From the facts that $\{l(nx - k) \mid k \in \N \setminus \{0\} \} \subseteq \Gamma(\alg{M})$, that $\Gamma(\alg{M})$ is finite, and that $\alg{M}$ is cancellative, it follows that $nx - k = 0$, and thus $nx = k$.
    This proves our claim, i.e., that for all $x \in M$ and $n,k \in \mathbb Z$, if $h(nx) = k$, then $nx = k$.
    It follows that, for all $x, y \in M$ with $h(x) = h(y)$, we have $h(x) \in \frac{1}{n} \mathbb Z$ and thus $nh(x) = nh(y) = k$ for some $k \in \Z$. Using the claim we obtain that $nx = ny =k$, which implies $x = y$ since
    $\alg M$ is torsion-free (\cref{l:torsion-free-monoid}). Thus $h$ is injective and $\alg{M} \cong \frac{1}{n} \mathbb Z$.
\end{proof}

\begin{corollary}\label{iso-to-Ln}
    Let $\alg{A}$ be a finite nontrivial totally ordered positive MV-algebra such that for all $x, y \in A$ we have $x \oplus y = 1$ or $x \odot y = 0$. There is $n \in \N \setminus \{0\}$ such that $\alg{A} \cong \alg{\L}^+_n$.
\end{corollary}

\begin{proof}
    This follows from \cref{prop: M tot ord,p:finite-Gamma-implies-Ln}, the equivalence between positive MV-algebras and cancellative unital commutative $\ell$-monoids, and the fact that $\Gamma(\frac{1}{n}\mathbb Z) = \alg{\L}_n^+$.
\end{proof}

\begin{theorem}
     \label{l:fin-is-sub}
    Every finite positive MV-algebra is a finite subdirect product of positive MV-algebras of the form $\alg{\L}^+_n$ for some $n  \in \mathbb{N} \setminus\{0\}$.
\end{theorem}

\begin{proof}
    Let $\alg{A}$ be a finite positive MV-algebra; this means that $\alg{A}$ is a subreduct of an MV-algebra $\alg{B}$.
    By Birkhoff's Subdirect Representation Theorem \cite[Theorem 8.6]{BurrisSanka}, $\alg{B}$ is a subdirect product of a family $\{\alg{B}_i\}_{i \in I}$ of subdirectly irreducible MV-algebras.
    Thus, $\alg{A}$ is a subalgebra of a product of a family $\{\alg{B}_i\}_{i \in I}$ of subdirectly irreducible MV-algebras.
    Every subdirectly irreducible MV-algebra is totally ordered (in fact, a nontrivial MV-algebra is totally ordered if and only if it is finitely subdirectly irreducible \cite[Theorem~15]{FontRodriguez1984}).
    Therefore, for every $i \in I$, $\alg{B}_i$ is totally ordered.
    Moreover, for every totally ordered MV-algebra $\alg{D}$ and all $x, y \in D$, we have $x \oplus y = 1$ (if $\lnot y \leq x$) or $x \odot y = 0$ (if $x \leq \lnot y$).
    Therefore, for all $i \in I$ and all $x, y \in B_i$, we have $x \oplus y = 1$ or $x \odot y = 0$.
    For every $i \in I$, let $p_i\colon \alg{B} \to \alg{B}_i$ be the $i$-th projection, and set $\alg A_i \coloneqq p_i(\alg A)$.
    For each $i \in I$, $\alg{A}_i$ is a finite nontrivial totally ordered positive subalgebra of $\alg{B}_i^+$ such that for all $x, y \in \alg{A}_i$ we have $x \oplus y = 1$ or $x \odot y = 0$.
    Therefore, by \cref{iso-to-Ln}, for every $i \in I$ there is $n_i \in \N \setminus \{0\}$ such that $\alg{A}_i \cong \alg{\L}_{n_i}^+$.
    This shows that $\alg{A}$ is a subdirect product of positive MV-algebras of the form $\alg{\L}_n^+$ for some $n \in \N \setminus \{0\}$.
    Since $\alg{A}$ is finite, finitely many of them suffice.
\end{proof}	

\begin{theorem}\label{thm:fin-is-sub}
    For a finite positive MV-algebra $\alg A$, the following are equivalent:
    \begin{enumerate}
        \item \label{i:subirr}
        $\alg A$ is subdirectly irreducible;
        \item \label{i:iso}
        $\alg A \cong \alg{\L}_n^+$ for some $n \in \N \setminus \{0\}$;
        \item \label{i:various}
        $\alg A$ is nontrivial, totally ordered and, for all $x,y \in A$, either $x \oplus y =1$ or $x \odot y=0$.
    \end{enumerate}
\end{theorem}
\begin{proof}
    The implication \eqref{i:subirr} $\Rightarrow$ \eqref{i:iso} follows from \cref{l:fin-is-sub}, the implication \eqref{i:iso} $\Rightarrow$ \eqref{i:subirr} follows from the fact that each $\alg{\L}_n^+$ is simple (\cref{lem:LnSimple}), the implication \eqref{i:iso} $\Rightarrow$ \eqref{i:various} is straightforward, and the implication \eqref{i:various} $\Rightarrow$ \eqref{i:iso} is \cref{iso-to-Ln}.
\end{proof}

By Di Nola’s representation theorem \cite{DiNola1991}, the standard MV-algebra $[0,1]$ generates the class of MV-algebras as a quasi-variety.
As observed in \cite[Proposition~2.5]{Abbadinietal2022}, an immediate consequence of this fact is that the class of positive MV-algebras is the quasivariety generated by the reduct $[0,1]^+$ of the standard MV-algebra $[0,1]$.
By \cite[Lemma 1.5]{CzelakowskiDziobiak1990}, every relatively subdirectly irreducible in $\mathsf{MV}^+$ is in $\II\SU\PP_u([0,1]^+)$ and hence is totally ordered.

\begin{remark}\label{r:01}
    The variety of MV-monoids generated by $[0,1]^+$ contains MV-monoids that are not positive MV-algebras.
    For example, $\alg{C}_2^\Delta \in \HH\SU\PP([0,1]^+)$.
    Indeed, let $\alg{S}$ be the subalgebra of $[0,1]^\omega$ consisting of the sequence constantly equal to $1$ and the sequences that converge to $0$.
    Let $\varphi \colon \alg{S} \to \alg{C}_2^\Delta$ be the homomorphism that maps the sequence constantly equal to $1$ to $1$, the sequences that are eventually equal to $0$ to $0$ and all the remaining ones to $\e$.
    This is a surjective homomorphism.
    This shows that $\alg{C}_2^\Delta \in \HH\SU\PP([0,1]^+)$.
    Analogously, also $\alg{C}_2^\nabla$ belongs to $\HH\SU\PP([0,1]^+)$.
    Note also that $\alg{C}_2^\Delta$ is an MV-monoid that is not cancellative. This gives another proof of the fact that $\mathsf{MV}^+ = \II\SU\PP_u([0,1]^+)$ is not a variety.
\end{remark}

To conclude the section we observe that we have produced two examples of (totally ordered) relatively subdirectly irreducible positive MV-algebras that are not subdirectly irreducible in the absolute sense, namely $\alg C^\Delta$ and $\alg C^\nabla$, i.e., the positive MV-subalgebras of $\alg C^+$ introduced in \cref{sec: prel}. Indeed, we can observe that $\bigwedge_{n\in \mathbb N} \gamma^\Delta_n = 0_{\con(\alg{C}^\Delta)}$. Hence $\alg{C}^\Delta$ is a positive MV-algebra which satisfies the conclusion of \cref{subdMVM} without being subdirectly irreducible.
We next prove that $\alg{C}^\Delta$ is relatively subdirectly irreducible.
For each $n \in \N \setminus \{0\}$ there is an element $a \in \alg{C}^\Delta/\g^\Delta_n$  (namely  $n\e/\g^\Delta_n$) such that $a \oplus a = a$ but $a \odot a =0$. Now the quasi-equations
$ x \oplus x \app x \Longleftrightarrow x \odot x \app x$ hold in MV-algebras and hence in positive MV-algebras. It follows that  $\alg{C}^\Delta/\g^\Delta_n \notin \mathsf{MV}^+$ for all $n \in \N \setminus \{0\}$. However
$\alg{C}^\Delta/\g^\Delta_0$ is the 2-element bounded distributive lattice, so $\alg{C}^\Delta/\g_0 \in \mathsf{MV}^+$; this implies that  the relative congruence lattice of $\alg{C}^\Delta$ (i.e., the lattice of congruences $\th \in \con(\alg{C}^\Delta)$  such that $\alg A/\th \in \mathsf{MVM}$) is the 3-element chain. So $\alg{C}^\Delta$ is {\em relatively subdirectly irreducible} in $\mathsf{MV}^+$ and it is also finitely subdirectly irreducible in $\mathsf{MVM}$.

Clearly also $\alg{C}^\nabla \in \mathsf{MV}^+$ and with a totally similar argument we can prove that $\bigwedge_{n \in \mathbb N}\g^\nabla_n = 0_{\con(\alg{C}^\nabla)}$; so $\op{Con}(\alg{C}^\nabla)$ is an infinite chain and $\alg{C}^\nabla$ is not subdirectly irreducible.
However, as above,  $\alg{C}^\nabla/\g^\nabla_n \notin \mathsf{MV}^+$ and so $\alg{C}^\nabla$ is relatively subdirectly irreducible in $\mathsf{MV}^+$.

\section{Varieties of positive MV-algebras}\label{sec: pMV varieties}

In this section we investigate varieties of positive MV-algebras. We will provide a characterization of the varieties of positive MV-algebras (\cref{t:VarofpMVs}): these are precisely the varieties generated by finitely many reducts of finite nontrivial MV-chains.
In fact, the varieties of positive MV-algebras are in bijection with the divisor-closed finite subsets of $\N \setminus\{0\}$ (\cref{d:div-closed}, \cref{cor:vredsets}).
Furthermore, we prove that positive MV-algebras form an unbounded sublattice of $\Lambda(\mathsf{MVM})$ (\cref{c:pMVvarieties}).

First we observe that not every positive MV-algebra generates a variety of positive MV-algebras: we have seen that $\alg C^+, \alg{C}^\Delta$, and $\alg{C}^\nabla$ have subalgebras having non-positive MV-algebras as quotients. Therefore $\VV(\alg C^+)$, $\VV (\alg{C}^\Delta)$ and
$\VV(\alg{C}^\nabla)$ are varieties of MV-monoids that do not consist entirely of positive MV-algebras. However, let:
\[
    L= \{\alg {\L}_n^+ \in \mathsf{MV}^+ \mid n \in \N \setminus \{0\}\}.
\]

\begin{lemma}\label{fingen} If $T$ is a finite subset of $L$, then
    $$
    \VV(T) = \II\SU\PP(T).
    $$
    Hence, $\VV(T)$ consists entirely of positive MV-algebras.
\end{lemma}
\begin{proof}
From J\'onsson's Lemma, in its version for a finite set of finite algebras \cite[Corollary 3.4]{Jonsson1967}, we have $\VV(T) = \II\PP_{\SU}\HH\SU(T)$, where with $\PP_{\SU}(K)$ we denote the class of all subdirect products of algebras in $K$.
Furthermore, any algebra in $T$ is hereditarily simple (\cref{lem:LnSimple}), and thus
\[
\VV(T) = \II\PP_{\SU}\SU(T) \subseteq \II\SU\PP\SU(T) \subseteq \II\SU\SU\PP(T) = \II\SU\PP(T).
\]
Since the class of positive MV-algebras is a quasivariety, it is closed under subalgebras and products.
Therefore, $\VV(T)$ consists entirely of positive MV-algebras.
\end{proof}

\begin{corollary}\label{c:finitely gen var of MVs}
    Every variety generated by finitely many finite positive MV-algebras consists of positive MV-algebras.
\end{corollary}
\begin{proof}
    Let $\vv K$ be a finite set of finite positive MV-algebras and let $\vv V \coloneqq \VV(\vv K)$.
    By \cref{l:fin-is-sub}, every finite positive MV-algebra is a finite subdirect product of algebras in $L$.
    Therefore, there is a finite subset $T$ of $L$ such that $\VV(\vv K) = \VV(T)$.
    By \cref{fingen}, the class $\VV(\mathcal{K}) = \VV(T)$ consists entirely of positive MV-algebras.
\end{proof}

\begin{lemma}\label{ispu}
    If $T$ is an infinite subset of $L$, then $\VV(T)$ coincides with $\VV([0,1]^+)$ and hence it is not a variety of positive MV-algebras.
\end{lemma}
\begin{proof}
    The inclusion $\VV(T) \subseteq \VV([0,1]^+)$ is obvious. To prove the converse inclusion,  set $I \coloneqq \{n \in \mathbb{N} \setminus \{0\} \mid \alg {\L}_n^+ \in T\}$.
    Let $\alg S$ be the subalgebra of $\prod_{n \in I} \alg{\L}_n$ consisting of the converging sequences.
    Let $\varphi \colon \alg{D} \to [0,1]^+$ be the homomorphism that maps a sequence to its limit.
    It is not difficult to see that, since $I$ is infinite, $\bigcup_{n \in I}{\alg{\L}_n^+}$ is dense in $[0,1]$.
    Therefore, $\varphi$ is surjective, which proves
    $[0,1]^+ \in \HH\SU\PP(T)$. This proves the inclusion $\VV([0,1]^+) \subseteq \VV(T)$.
    By \cref{r:01}, it follows that $\VV(T)$ is not a variety of positive MV-algebras.
\end{proof}

Note that the first part of the statement of \cref{ispu} holds also for MV-algebras \cite[Proposition~8.1.2]{CignoliDOttavianoEtAl2000}.

\begin{corollary}
    There is no largest variety of positive MV-algebras.
\end{corollary}

\begin{proof}
    If one such variety $\VV$ existed, then by \cref{fingen} we would have $\alg {\L}^+_n \in \vv V$ for all $n$, and then by \cref{ispu} the variety $\vv V$ would not consist entirely of positive MV-algebras, a contradiction.
\end{proof}	

\cref{c:finitely gen var of MVs} states that the variety generated by finitely many finite positive MV-algebras is a variety of positive MV-algebras.
We will now show that these are the only varieties of positive MV-algebras.
In other words, the varieties of positive MV-algebras are exactly the varieties generated by a finite set of finite positive MV-algebras.
The proof relies on \cref{t: cancellative}, the proof of
which uses the following lemma.

\begin{lemma}\label{l: y=k}
    Let $\alg{M}$ be a nontrivial totally ordered unital commutative $\ell$-monoid such that every quotient of every subalgebra of $\alg{M}$ is cancellative. Let $\varphi$ be the unique homomorphism from $\alg{M}$ to $\R$ (which exists by \cref{Holder}). For all $y \in M$ and $k \in \Z$, if $\varphi(y) = k$, then $y = k$.
\end{lemma}

\begin{proof}
    Let $S = \{x \in M \mid \varphi(x) \in \Z\}$.
    By \cref{l: T subalgebra}, $\alg{S}$ is a subalgebra of $\alg M$.
    We shall prove that the unique homomorphism $\Z \to \alg{S}$ is surjective. By way of contradiction, suppose this is not the case.
    Then, $\alg{S}$ is not isomorphic to $\Z$ (by uniqueness of the homomorphism from $\Z$ to $\alg{S}$), and so we can apply \cref{t:for-almost-minimal} to $\alg{S}$.
    Therefore, one of the following conditions holds: (i) the image of the unique homomorphism $\varphi \colon \alg{S} \to \Z$ is not $\Z$, (ii) $\Z \overrightarrow{\times} \alg C_2^{\Delta*}$, (iii) $\Z \overrightarrow{\times} \alg C_2^{\nabla*}$.
    But (i) contradicts the definition of $\alg{S}$, and (ii) and (iii) are excluded since $\Z \overrightarrow{\times} \alg C_2^{\Delta*}$ and $\Z \overrightarrow{\times} \alg C_2^{\nabla*}$ are not cancellative.
\end{proof}

\begin{theorem}\label{t: cancellative}
    Let $\alg{M}$ be a nontrivial totally ordered unital commutative $\ell$-monoid such that every quotient of every subalgebra of $\alg{M}$ is cancellative.
    Then, the unique homomorphism from $\alg{M}$ to $\R$ is injective, and thus $\alg{M}$ is isomorphic to a subalgebra of $\R$.
\end{theorem}

\begin{proof}
    Let $\varphi$ be the unique homomorphism from $\alg{M}$ to $\R$, which exists by \cref{Holder}.
    Let $x, y \in M$ be such that $\varphi(x) = \varphi(y)$.
    We shall prove $x = y$.
    We have two cases:
    \begin{enumerate}
        \item \label{i:Q}
        $\varphi(x) \in \Q$;

        \item \label{i:notQ}
        $\varphi(x) \notin \Q$.
    \end{enumerate}

    Case \eqref{i:Q}.
    In this case, there are $n \in \mathbb N \setminus\{0\}$ and $k \in \Z$ such that $\varphi(x) = \varphi(y) = \frac{k}{n}$, i.e., $n\varphi(x) = n \varphi(y) = k$, which implies $\varphi(nx) = \varphi(ny) = k$, and thus $nx = ny = k$ by \cref{l: y=k}, and hence $x = y$ by \cref{l:torsion-free-monoid}.

    Case \eqref{i:notQ}.
    Without loss of generality, we may suppose $x \leq y$.
    We set $W \coloneqq \{n \varphi(x) + k \in \R \mid n \in \N, k \in \Z\}$.
    Note that $\alg{W}$ is a subalgebra of the unital commutative $\ell$-monoid $\R$.
    We set
    \[
    R \coloneqq \varphi^{-1}[W].
    \]
    Since $\varphi(x)$ is irrational, for every $z \in R$ there is exactly one pair $(n,k) \in \N \times \Z$ such that $\varphi(z) = n\varphi(x) + k$.
    We denote the two elements of this pair by $N_z$ and $C_z$, respectively.
    Note that $\varphi(z_1) = \varphi(z_2)$ implies $N_{z_1} = N_{z_2}$ and $C_{z_1} = C_{z_2}$.
    Moreover, for all $z_1, z_2 \in R$, we have:
    \[
    (N_{z_1}+N_{z_2})\varphi(x) + (C_{z_1} + C_{z_2}) = N_{z_1}\varphi(x) + C_{z_1} + N_{z_2}\varphi(x) + C_{z_2} = z_1 + z_2,
    \]
    which proves that $z_1 + z_2 \in R$ with $N_{z_1 + z_2} = N_{z_1} + N_{z_2}$ and $C_{z_1 + z_2} = C_{z_1} + C_{z_2}$.
    Let
    \[
    T' \coloneqq \{z \in R \mid N_z x + C_z \leq z\}.
    \]
    We claim that $\alg{T}'$ is a subalgebra of $\alg{M}$.
    Since $\alg{M}$ is totally ordered, its subset $\alg{T}'$ is a sublattice.
    Moreover, $\alg{T}'$ contains $0$, $1$ and $-1$, as witnessed by the pairs $(0,0), (0,1), (0, -1) \in \N \times \Z$, respectively.
    To prove closure under $+$, let $z_1, z_2 \in T'$.
    We already know that $z_1 + z_2 \in R$.
    Moreover,
    \begin{align*}
        N_{z_1 + z_2}x + C_{z_1 + z_2} & = (N_{z_1} + N_{z_2}) x + (C_{z_1} + C_{z_2})\\
        & = N_{z_1}x + C_{z_1} + N_{z_2}x + C_{z_2} \\
        & \leq z_1 + z_2.
    \end{align*}
    Therefore, $z_1 + z_2 \in T'$, proving the closure of $T'$ under $+$.
    This proves our claim that $\alg{T}'$ is a subalgebra of $\alg{M}$.

    We define the function
    \begin{align*}
    \psi \colon \alg{T}' & \longrightarrow \R \overrightarrow{\times} \alg{C}_2^{\Delta*}\\
    z & \longmapsto \begin{cases}
        (\varphi(z), 0) & \text{if }z = N_z x + C_z,\\
        (\varphi(z), \e) & \text{otherwise.}
    \end{cases}
    \end{align*}
    We show that this is a homomorphism.
    We have $\psi(0) = (0,0)$, $\psi(1) = (1,0)$, and $\psi(-1) = (-1,0)$. Thus $\psi$ preserves $0$, $1$ and $-1$.
    We prove that $\psi$ is order-preserving: let $z_1, z_2 \in T'$ be such that $z_1 \leq z_2$. Then $\varphi(z_1) \leq \varphi(z_2)$, and so $K_{z_1}\varphi(x) + C_{z_1} \leq K_{z_2}\varphi(x) + C_{z_2}$. If $K_{z_1}\varphi(x) + C_{z_1} < K_{z_2}\varphi(x) + C_{z_2}$, then $\psi(x) < \psi(y)$ by inspection on the first coordinates of $\psi(x)$ and $\psi(y)$.
    If $K_{z_1}\varphi(x) + C_{z_1} = K_{z_2}\varphi(x) + C_{z_2}$, either $z_2 = K_{z_1}x + C_{z_1}$ and then $z_1 = K_{z_1}x + C_{z_1}$, since $z_1 \leq z_2$ and $\varphi(z_1) = \varphi(z_2) = K_{z_1}\varphi(x) + C_{z_1}$, or $z_2 > K_{z_1}x + C_{z_1}$ and thus $\psi(z_1) = (K_{z_1}\varphi(x) + C_{z_1}, t) \leq (K_{z_1}\varphi(x) + C_{z_1}, \e) = \psi(z_2)$, for some $t \in \{0,\e\}$.
    This proves that $\psi$ is order-preserving, and thus $\psi$ is a lattice homomorphism because $\alg T'$ is totally ordered.

    We prove that $\psi$ preserves $+$. Let $z_1,z_2 \in T'$ and let $(K_{z_1}\varphi(x) + C_{z_1},t_{z_1}) = \psi(z_1)$ and $(K_{z_2}\varphi(x) + C_{z_2},t_{z_2}) = \psi(z_2)$ for some $t_{z_1},t_{z_2} \in \{0,\e\}$. We prove that
    \[
    \psi(z_1 + z_2) = (K_{z_1}\varphi(x) + C_{z_1} + K_{z_2}\varphi(x) + C_{z_2}, \max(\{t_{z_1},t_{z_2}\})) = \psi(z_1) + \psi(z_2).
    \]
    If $z_1 = K_{z_1}x + C_{z_1}$ and $z_2 = K_{z_2}x + C_{z_2}$, then
    \begin{align*}
        \psi(z_1 + z_2) & = \psi(K_{z_1}x + C_{z_1} + K_{z_2}x + C_{z_2})  \\
        & =  (K_{z_1}\varphi(x) + C_{z_1} + K_{z_2}\varphi(x) + C_{z_2}, 0)  \\
        & = \psi(z_1) + \psi(z_2).
    \end{align*}
    If $z_1 > K_{z_1}x + C_{z_1}$ or $z_2 >  K_{z_2}x + C_{z_2}$, then
    \[
    \varphi(z_1 + z_2) = \varphi(z_1) + \varphi(z_2) = K_{z_1}\varphi(x) + C_{z_1} + K_{z_2}\varphi(x) + C_{z_2}
    \]
    and $z_1 + z_2 > K_{z_1}x + C_{z_1} + K_{z_2}x + C_{z_2}$. Thus,
    \begin{align*}
    \psi(z_1 + z_2) & = \psi(K_{z_1}x + C_{z_1} + K_{z_2}x + C_{z_2})\\
    & =  (K_{z_1}\varphi(x) + C_{z_1} + K_{z_2}\varphi(x) + C_{z_1}, \e)\\
    & = \psi(z_1) + \psi(z_2).
    \end{align*}

    Note that $N_x + N_y = 1$ and $C_x = C_y = 0$.
    Therefore, $\psi(x) = (\varphi(x), 0)$, because $x = 1x + 0$.
    Furthermore, there is $\nu \in \{0, \e\}$ such that $\psi(y) = (\varphi(x),\nu)$.
    We have
    \begin{align*}
        \psi(x) + \psi(y) & = (\varphi(x), 0) + (\varphi(x),\nu) \\
        & = (2 \varphi(x), \nu) \\
        & = (\varphi(x), \nu) + (\varphi(x),\nu) \\
        & = \psi(y) + \psi(y).
    \end{align*}
    By hypothesis, the image of $\psi$ is a cancellative subalgebra of $\R \overrightarrow{\times} \alg{C}_2^{\Delta*}$.
    Thus, from $\psi(x) + \psi(y) = \psi(y) + \psi(y)$ we deduce $\psi(x) = \psi(y)$.
    Hence, $\psi(y) = \psi(x) = (\varphi(x), 0)$.
    Therefore, $y = N_y x + C_y = 1x + 0 = x$.
\end{proof}

\begin{corollary} \label{c:subalgebra}
    Let $\alg{A}$ be a nontrivial totally ordered positive MV-algebra such that for every $x,y \in A$ we have $x \oplus y = 1$ or $x \odot y = 0$.
    Moreover, suppose that every quotient of every subalgebra of $\alg{A}$ is a positive MV-algebra.
    Then $\alg{A}$ is isomorphic to a subalgebra of $[0,1]^+$.
\end{corollary}
\begin{proof}
     First we observe that $\Xi(\alg A)$ (we refer to p.~\pageref{def:quasiinverse} for the definition of $\Xi$) is nontrivial since $\alg A$ is nontrivial, is totally ordered by \cref{prop: M tot ord}, and is cancellative since $\alg A$ is a positive MV-algebra. Furthermore, by \cref{preserveandreflect}, every quotient $\alg Q$ of $\Xi(\alg A)$ is such that $\Gamma(\alg Q)$ is isomorphic to a quotient of $\alg A$ and hence to a positive MV-algebra. Thus, $\alg Q$ is cancellative for each quotient $\alg Q$ of $\Xi(\alg A)$.

    From \cref{t: cancellative} we have that $\Xi(\alg A)$ is isomorphic to a subalgebra of $\R$ and thus $\alg A \cong \Gamma(\Xi(\alg A))$ is isomorphic to $\Gamma(\R) = [0,1]^+$ still by \cref{preserveandreflect}.
\end{proof}
Using mainly \cref{c:subalgebra} we are able to prove the main result of the section, which characterizes all varieties of positive MV-algebras in terms of some of their finite subdirectly irreducible members.

\begin{theorem} \label{t:subirr-in-variety}
    Let $\mathcal{V}$ be a variety of positive MV-algebras.
    For every subdirectly irreducible member $\alg{A}$ of $\mathcal{V}$ there is $n \in \N\setminus\{0\}$ such that $\alg{A}$ is isomorphic to $\alg{\L}^+_n$.
\end{theorem}

\begin{proof}
    Let $\alg A$ be a subdirectly irreducible algebra in $\vv V$.
    By \cref{prop: M tot ord}, $\alg A$ is nontrivial, totally ordered, and such that, for all $x,y \in A$, $x \oplus y =1$ or $x \odot y =0$.
    Since $\HH\SU(\alg A) \subseteq \mathcal{V}$, every algebra in $\HH\SU(\alg A)$ is a positive MV-algebra.
    So, by \cref{c:subalgebra}, $\alg A $ is isomorphic to a subalgebra $\alg{B}$ of $[0,1]^+$.
    By \cref{l:dense-ord-dicrete-unital}, either there is $n \in \mathbb N \setminus \{0\}$ such that $\alg{B} = \alg{\L}_n^+$ or $A$ is dense in $[0,1]$. However, the latter is not possible, because otherwise by \cref{l:generates-as-much-as-01} we would have $[0,1]^+ \in \vv V$, which, by \cref{r:01}, would imply the existence of an algebra in $\vv V$ that is not a positive MV-algebra, contradicting the fact that $\vv V$ consists of positive MV-algebras.
\end{proof}

\begin{theorem}\label{t:VarofpMVs}
    The varieties of positive MV-algebras are precisely the varieties generated by a finite set of finite positive MV-algebras.
    Equivalently, they are precisely the varieties generated by a finite subset of $\{\alg {\L}_n^+ \mid n \in \N \setminus \{0\}\}$.
\end{theorem}
\begin{proof}
    Let $\vv V$ be a variety of positive MV-algebras. By \cref{t:subirr-in-variety}, any subdirectly irreducible algebra $\alg A$ in $\vv V$ is isomorphic to $\alg {\L}_n^+$ for some $n \in \mathbb{N} \setminus \{0\}$. Furthermore, the set $\{n \in \N \setminus\{0\}\mid \alg {\L}_n^+ \in \VV\}$ is finite because otherwise, by \cref{ispu}, $\vv V$ would not consist entirely of positive MV-algebras. Thus, by Birkhoff's Subdirect Representation Theorem \cite[Theorem 8.6]{BurrisSanka}, $\vv V$ is generated by a finite subset of $\{\alg {\L}_n^+ \mid n \in \N \setminus \{0\}\}$.

    Clearly, every variety generated by a finite subset of $\{\alg {\L}_n^+ \mid n \in \N \setminus \{0\}\}$ is generated by a finite set of finite positive MV-algebras. Moreover, by \cref{c:finitely gen var of MVs}, every variety generated by a finite set of finite positive MV-algebras is a variety of positive MV-algebras.
\end{proof}

\begin{corollary}\label{c:pMVvarieties}
   The varieties of positive MV-algebras form an unbounded sublattice (really, an ideal) of the lattice of subvarieties of $\mathsf{MVM}$.
\end{corollary}

\begin{proof}
   We need only to show that the join (as varieties of MV-monoids) of two varieties of positive MV-algebras is still a variety of positive MV-algebras. But this
    follows from \cref{t:VarofpMVs} since the join of two varieties is the variety generated by the union of the two generating classes.
\end{proof}

\begin{definition}\label{d:div-closed}
    We call \emph{divisor-closed finite set} a finite downset of the divisibility poset $(\N\setminus \{0\}, {\mid})$, i.e., a finite subset $I$ of $\N \setminus \{0\}$ such that, for every $n \in I$, every divisor $k \in \N \setminus \{0\}$ of $n$ belongs to $I$.
\end{definition}
For a divisor-closed finite set $I$, we set
\[
    \mathcal{K}_I\coloneqq\{\alg \L_n^+ \in \mathsf{MV}^+\mid n \in I\}.
\]

We can refine the statement and proof of \cref{t:VarofpMVs} to get the characterization of the varieties of positive MV-algebras in the following result.
We recall that $\alg{\L}^+_n$ is simple for any $n \in \N \setminus \{0\}$ (\cref{lem:LnSimple}). Moreover, it is not difficult to show that $\alg{\L}^+_n \in \II\SU(\alg{\L}^+_m)$ if and only if $n$ divides $m$.

\begin{theorem}\label{cor:vredsets}
    The set $\Lambda(\mathsf{MV}^+)$ of varieties of positive MV-algebras is in bijection with the set $\mathcal{J}$ of divisor-closed finite sets, as witnessed by the inverse functions:
    \begin{align*}
        f \colon \mathcal{J} & \longrightarrow \Lambda(\mathsf{MV}^+)
        & & &
        g \colon \Lambda(\mathsf{MV}^+)& \longrightarrow\mathcal{J}\\
        I & \longmapsto \VV(\mathcal{K}_I)
        & & &
          \mathcal{V} & \longmapsto \{ n \in \N\setminus \{0\} \mid \alg{\L}^+_n \in \mathcal{V}\}.
    \end{align*}
\end{theorem}

\begin{proof}
    Let $f \colon \mathcal{J} \to \Lambda(\mathsf{MV}^+)$ and $g \colon \Lambda(\mathsf{MV}^+) \to \mathcal{J}$ be the functions in display.
    The function $f$ is well-defined by \cref{c:finitely gen var of MVs}.
    The function $g$ is well-defined since for all $\vv V \in \Lambda(\mathsf{MV}^+)$, if $m \mid n \in g(\vv V)$ then $\alg \L^+_m \in \II\SU(\alg \L^+_n) \subseteq \vv V$ and thus $m \in g(\vv V)$.

    We prove that $f$ and $g$ are inverse functions. The composite $f \circ g \colon \Lambda(\mathsf{MV}^+) \to \Lambda(\mathsf{MV}^+)$ is the identity because every variety is generated by its subdirectly irreducible members, by Birkhoff's Subdirect Representation Theorem \cite[Theorem 8.6]{BurrisSanka}, and every subdirectly irreducible member of a variety of positive MV-algebras is isomorphic to $\alg{\L}^+_n$ for some $n \in \N \setminus \{0\}$ (\cref{t:subirr-in-variety}).

    We prove that the composite $g \circ f \colon \mathcal{J} \to \mathcal{J}$ is the identity.
    Let $I \in \mathcal{J}$.
    Since $\mathcal{K}_I \subseteq \VV(\mathcal{K}_I)$, it is immediate that $I \subseteq g(f(I))$.
    For the converse inclusion, let $n \in g(f(I))$.
    Then, $\alg{\L}^+_n \in \mathcal{V}(\mathcal{K}_I)$.
    Since $\alg \L^+_n$ is subdirectly irreducible, via an application of J\'onsson's Lemma \cite[Corollary 3.4]{Jonsson1967} we get $\alg \L^+_n \in \HH\SU(\mathcal{K}_I)$. Therefore, there is $m \in I$ such that $\alg{L}^+_n$ is a homomorphic image of a subalgebra of $\alg{L}^+_m \in \mathcal{K}_I$.
    Since $\alg \L^+_m$ is hereditary simple by \cref{lem:LnSimple} and since $\alg{\L}^+_n$ is nontrivial, $\alg{\L}_n^+$ is isomorphic to a subalgebra of $\alg \L^+_m$.
    Thus, $n$ divides $m$ and hence $n \in I$ since $I$ is a divisor-closed finite set.
    This proves the inclusion $g(f(I)) \subseteq I$ and so we have $I = g(f(I))$.
\end{proof}

\begin{remark}
    A \emph{reduced set} is a finite antichain of the divisibility poset $(\N\setminus \{0\}, \mid)$, i.e., a finite subset $I$ of $\N \setminus \{0\}$ such that, for all distinct $n, m \in I$, $n \nmid m$.
    The set $\mathcal{J}$ in \cref{cor:vredsets} of divisor-closed finite sets is clearly in bijection with the set of reduced sets since, in posets in which principal downsets are finite, finite downsets correspond to finite antichains.
\end{remark}

By \cref{cor:vredsets}, the set $\Lambda(\mathsf{MV}^+)$ of all varieties of positive MV-algebras is countably infinite, since $\Lambda(\mathsf{MV}^+)$ is in bijection with an infinite subset of the subset of all finite members of the power set $\mathcal{P}(\N)$. \cref{cor:vredsets} will also be the main tool used to find axiomatizations of the varieties of positive MV-algebras in the next section.

\section{Axiomatizations}\label{sec:axiomatization}

In this section, we present an axiomatization of all varieties of positive MV-algebras and all almost minimal varieties of MV-monoids. For the varieties of positive MV-algebras our strategy will be as follows: by \cref{cor:vredsets} each variety of positive MV-algebras is finitely generated; since any such variety is congruence distributive, we can apply Baker's Finite Basis Theorem \cite{Baker1977} to conclude that any such variety has a finite axiomatization. In particular, for any divisor-closed finite set $I$ there is a finite set of equations holding in $\VV(\mathcal{K}_I)$ and implying (relatively to MV-monoids) the quasiequation
\begin{equation}\label{e:cancMv}
  (x\oplus z \app y \oplus z) \meet (x \odot z \app y \odot z)  \ \Longrightarrow \  x \app y.
\end{equation}
Our idea to find such a finite set of equations is as follows: for each $n \in \N \setminus \{0\}$ we find a set of equations $\Phi_n$ (\cref{d:defsi}) that axiomatizes $\VV(\alg{\L}_n)$ within the variety of MV-monoids (\cref{thm:axiomsofln}).
Then $\Phi_{\lcm(I)}$ (where we set $\lcm(\emptyset) = 1$) is a finite set of equations holding in $\VV(\mathcal{K}_I)$ and implying (relatively to MV-monoids) the quasiequation \eqref{e:cancMv}.

Next, we consider further axioms to distinguish between two varieties $\vv V(\mathcal{K}_I)$ and $\vv V(\mathcal{K}_J)$ with $I \not= J$ and $\lcm(I) = \lcm(J)$, for $I$ and $J$ divisor-closed finite sets.
Letting $m$ denote the maximum of $I$ (or setting $m = 0$ if $I$ is empty) we will define $\Sigma_I$ as the set of equations given by the equation
\begin{equation}\label{eq:splittingmult}
    (m+1)x \app mx
\end{equation}
together with the equations of the form
\begin{equation}\label{eq:splittingdiv}
     m((k-1)x)^k \app (kx) ^m
\end{equation}
for all $1 \leq k \leq m$ such that $k \notin I$.
To give a rough idea of the meaning of the equations in \eqref{eq:splittingmult} and \eqref{eq:splittingdiv}, we may say that they are constructed to exclude
selected $\alg \L_n^+$'s from the generating set of a variety of positive MV-algebras. Namely, \eqref{eq:splittingmult} excludes the subdirectly irreducible generators $\alg{\L}^+_{n}$'s with $n > m$. Similarly, \eqref{eq:splittingdiv} excludes the $\alg{\L}^+_{n}$'s with $n \leq m$ and $n$ not in $I$.

In conclusion, for any divisor-closed finite set $I$, we will prove that a finite equational axiomatization of $\VV(\vv{K}_I)$ consists of the set of axioms of MV-monoids together with $\Phi_{\mathrm{lcm}(I)} \cup \Sigma_I$ (\cref{thm:axioms}).

We now start to deal with the sets $\Phi_{\lcm(I)}$ of equations holding in $\VV(\mathcal{K}_I)$ and implying the cancellation law.
In order to define $\Phi_{\lcm(I)}$ we introduce a few definitions regarding MV-monoids and unital commutative $\ell$-monoids.

\begin{definition}
    Let $\alg A$ be an MV-monoid.
    An element $x \in A$ is said to be \emph{idempotent} if $x \oplus x = x$ and $x \odot x = x$.
\end{definition}

\begin{definition}\label{def:integer}
    Let $\alg M$ be a unital commutative $\ell$-monoid.
    An element $x \in M$ is said to be \emph{integer} if $x = k$ for some $k \in \Z$.
\end{definition}
It is easy to observe that any integer element is invertible.

We are now ready to introduce the set of equations $\Phi_n$ in the language of MV-monoids. The idea behind those equations is that for any subdirectly irreducible MV-monoid $\alg A$ satisfying $\Phi_n$ and for any $x \in A$, $nx$ should be an integer element of $\Xi(\alg A)$. This is the key tool that will allow us to prove the cancellation law for the MV-monoids satisfying $\Phi_n$.

\begin{definition}
    For $n \in \N$ and $k \in \Z$ we define a unary term $\tau_{n,k}(x)$ inductively on $n$.
    The idea is that $\tau_{n,k}(x)$ equals $((nx - k) \lor 0) \land 1$ (where the computations are done in the enveloping unital commutative $\ell$-monoid).
    The base case ($n = 0$), is set as follows:
    \[
     \tau_{0, k}(x) \coloneqq \begin{cases}
         1 & \text{if }k \leq -1,\\
         0 & \text{if }k \geq 0.
     \end{cases}
    \]
    The inductive case is as follows:
    \begin{equation} \label{eq:tau-1}
        \tau_{n +1, k}(x) = \tau_{n,k-1}(x) \odot (x \oplus \tau_{n, k}(x)),
    \end{equation}
    or, equivalently,
    \begin{equation} \label{eq:tau-2}
        \tau_{n +1, k}(x) = (\tau_{n,k-1}(x) \odot x) \oplus \tau_{n, k}(x)).
    \end{equation}
\end{definition}

\begin{remark}
    To see that \eqref{eq:tau-1} and \eqref{eq:tau-2} are equivalent it is enough to observe that $(\tau_{n,k-1}(x), \tau_{n, k}(x))$ is a good pair (\cref{d:goodpair}); indeed, if $(x_0, x_1)$ is a good pair in an MV-monoid $\alg{A}$, and $y \in {A}$, then $x_0 \odot (y \oplus x_1) = (x_0 \odot y) \oplus x_1$ (\cite[Lemma~4.35]{Abbadini2021}).
    The fact that $(\tau_{n,k-1}(x), \tau_{n, k}(x))$ is a good pair follows from \cref{l:taunk} below (in the proof of which we use only the first definition) and from the fact that, for each $x$ in a unital commutative $\ell$-monoid $\alg{M}$, the pair
    \[
    \big(((x -1) \lor 0) \land 1,(x \lor 0) \land 1\big)
    \]
    is a good pair \cite[Proposition~4.64]{Abbadini2021}.
\end{remark}

\begin{example}
    We compute the terms $\tau_{n,k}(x)$ for $n \in \{ 1,2,3\}$. For $n = 1$ we get
    \[
    \tau_{1,k}(x) = \begin{cases}
        1 & \text{if }k \leq -1,\\
        x & \text{if }k = 0,\\
        0 & \text{if }k \geq 1,
    \end{cases}
    \]
    for $n =2$ we have
    \[
    \tau_{2,k}(x) = \begin{cases}
        1 & \text{if }k \leq -1,\\
        x \oplus x & \text{if }k = 0,\\
        x \odot x & \text{if }k =1,\\
        0 & \text{if }k \geq 2,
    \end{cases}
    \]
    and for $n =3$ we obtain
    \[
    \tau_{3,k}(x) = \begin{cases}
        1 & \text{if }k \leq -1,\\
        x \oplus x \oplus x & \text{if }k = 0,\\
        (x \oplus x) \odot (x \oplus (x \odot x)), & \text{if }k =1,\\
        x \odot x \odot x & \text{if }k =2,\\
        0 & \text{if }k \geq 3.
    \end{cases}
    \]
\end{example}
Note that, for every $n$, we have $\tau_{n,k}(x) = 1$ for all $k \leq -1$, and  $\tau_{n,k}(x) = 0$ for all $k \ge n$. Thus, the only interesting computations happen for $0 \leq k \leq n-1$. Moreover, $\tau_{n,1}(x) = nx$ and $\tau_{n,n-1}(x) = x^n$.

In this section we will often consider equalities where one side is computed using the operations of commutative $\ell$-monoid (and the result is in the interval $[0,1]$) and the other using the operations of MV-monoids. An example of this fact is given by the next lemma with the equality $((nx - k) \lor 0) \land 1 = \tau_{n, k}(x)$. This should not cause any confusion since given a unital commutative $\ell$-monoid $\alg M$ we have $((nx - k) \lor 0) \land 1 \in \Gamma(\alg{M})$ for all $x \in \Gamma(\alg{M})$, since by \cref{eq:Gamma} $\Gamma$ restricts $M$ to its unit interval and $0 \leq ((nx - k) \lor 0) \land 1 \leq 1$.

\begin{lemma} \label{l:taunk}
    Let $\alg{M}$ be a unital commutative $\ell$-monoid. For every $x \in \Gamma(\alg{M})$, $n \in \N$ and $k \in \Z$, we have
    \[
        ((nx - k) \lor 0) \land 1 = \tau_{n, k}(x),
    \]
    where the left-hand side is computed in $\alg{M}$ and the right-hand side in $\Gamma(\alg{M})$.
\end{lemma}
\begin{proof}
    We prove this by induction on $n \in \N$.

    The base case $n = 0$ is immediate.

    For the inductive case, let $n \geq 0$, and suppose that, for all $k \in \Z$, we have
    \[
        ((nx - k) \lor 0) \land 1 = \tau_{n, k}(x).
    \]
    We should prove
    \[
        (((n+1)x - k) \lor 0) \land 1 = \tau_{n + 1, k}(x),
    \]
    i.e.,
    \[
        (((nx - k) + x) \lor 0) \land 1 = \tau_{n + 1, k}(x).
    \]
    Set $z \coloneqq nx - k$.
    We shall prove
    \[
        ((z + x) \lor 0) \land 1 = \tau_{n + 1, k}(x).
    \]
    Set $z_{-1} \coloneqq ((z - 1) \lor 0) \land 1$ and $z_0 \coloneqq (z \lor 0) \land 1$.
    By \cite[Lemma 4.65]{Abbadini2021}, we have
    \begin{equation}\label{eq:z+x}
        ((z + x) \lor 0) \land 1 = z_{-1} \odot (x \oplus z_0).
    \end{equation}
    By the inductive hypothesis, we have
    \begin{equation} \label{eq:z0}
        z_0 \coloneqq (z \lor 0) \land 1 = ((nx - k) \lor 0) \land 1 = \tau_{n,k}(x),
    \end{equation}
    and
    \begin{equation} \label{eq:z1}
        z_{-1} \coloneqq ((z - 1) \lor 0) \land 1 = ((nx - (k -1)) \lor 0) \land 1 = \tau_{n,k-1}(x).
    \end{equation}
    By \eqref{eq:z+x}, \eqref{eq:z0} and \eqref{eq:z1}, we have
    \[
        ((z + x) \lor 0) \land 1 = z_{-1} \odot (x \oplus z_0) = \tau_{n,k-1}(x) \odot (x \oplus \tau_{n,k}(x)) = \tau_{n+1,k}. \qedhere
    \]
\end{proof}

\begin{definition}\label{d:defsi}
    For every $n \in \N$, let $\Phi_{n}$ be the following set of equations, for $k$ ranging in $\{0, \dots, n-1\}$:
    \[
    \tau_{n,k}(x) \oplus \tau_{n,k}(x) \approx \tau_{n,k}(x)
    \]
    and
    \[
    \tau_{n,k}(x) \odot \tau_{n,k}(x) \approx \tau_{n,k}(x).
    \]
    In other words, we impose that $\tau_{n, k}(x)$ is idempotent.
\end{definition}

The following is a key fact connecting integer elements of a totally ordered unital commutative $\ell$-monoid $\alg M$ and the satisfaction in $\Gamma(\alg M)$ of the equations in $\Phi_{n}$.

\begin{proposition}\label{p:totord-then-integer}
    Let $\alg{M}$ be a totally ordered unital commutative $\ell$-monoid, let $x \in \Gamma(\alg{M})$ and let $n \in \N$. The following conditions are equivalent:
    \begin{enumerate}
        \item \label{i:equals-sum}
        the element $nx$ is integer;
        \item \label{i:idempotent}
        for every $k \in \{0, \dots, n-1\}$, $\tau_{n,k}(x)$ is an idempotent element of $\Gamma(\alg{M})$, i.e., the element $x$ satisfies $\Phi_n$.
    \end{enumerate}
\end{proposition}

\begin{proof}
    \eqref{i:equals-sum} $\Rightarrow$ \eqref{i:idempotent}. This implication is immediate since if $nx$ is integer then $nx = k'$, for some $k' \in \Z$. Thus $((nx - k) \lor 0) \land 1 = ((k' - k) \lor 0) \land 1$ which is $1$ if $k'>k$ and $0$ otherwise. Hence, by \cref{l:taunk}, either $\tau_{n, k}(x) = 1$ or $\tau_{n, k}(x) = 0$, and therefore $\tau_{n, k}(x)$ is idempotent.

    \eqref{i:idempotent} $\Rightarrow$ \eqref{i:equals-sum}.
    We first prove that every idempotent element of $\Gamma(\alg{M})$ is either $0$ or $1$.
    Let $y$ be an idempotent element of $\Gamma(\alg{M})$.
    Then $y \oplus y = y$ and $y \odot y = y$.
    By \cref{prop: M tot ord}, since $\alg{M}$ is totally ordered, $y \oplus y = 1$ or $y \odot y = 0$, i.e., either $y = 1$ or $y = 0$.
    This proves that every idempotent element of $\Gamma(\alg{M})$ is either $0$ or $1$.

    Therefore, by \eqref{i:idempotent}, for every $k \in \Z$ the element $\tau_{n,k}(x)$ is either $0$ or $1$.
    By \cref{l:taunk}, for every $k \in \Z$ we have $\tau_{n,k}(x) = ((nx - k) \lor 0) \land 1$, and thus $((nx - k) \lor 0) \land 1$ is either $0$ or $1$.
    Since the sequence $\{((nx - k) \lor 0) \land 1)\}_{k \in \Z}$ is increasing and $\alg M$ is totally ordered, it follows that there is $k_0 \in \Z$ such that for every $k < k_0$ we have $((nx - k) \lor 0) \land 1 = 1$ and for every $k \geq k_0$ we have $((nx - k) \lor 0) \land 1 = 0$.
    By \cref{abbadinith}, $nx = k_0$.
\end{proof}

We now characterize for which $m$ and $n$ the algebra $\alg \L_m^+$ satisfies $\Phi_n$.

\begin{lemma}\label{l:lnSigman}
    Let $m,n \in \N$ with $n >0$. Then $\alg \L_m^+$ satisfies $\Phi_n$ if and only if $m \mid n$.
\end{lemma}

\begin{proof}
Recall that $\alg{\L}_m^+ = \Gamma(\frac{1}{n}\Z)$.

For $(\Rightarrow)$,
suppose $\alg{L}_m^+$ satisfies $\Phi_n$.
Then, by \cref{p:totord-then-integer} applied to the unital commutative $\ell$-monoid $\frac{1}{m}\Z$ and to the element $x = \frac{1}{m}$, the element $n\frac{1}{m}$ is an integer element of $\frac{1}{m}\Z$, i.e., $\frac{n}{m} \in \Z$.
This implies that $m$ divides $n$.

For $(\Leftarrow)$, suppose $m$ divides $n$.
For every $x \in \Gamma(\frac{1}{m}\Z) = \alg{\L}_m^+$, the element $nx \in \frac{1}{m}\Z$ belongs to $\Z$ because $m$ divides $n$, and hence is an integer element.
By \cref{p:totord-then-integer} $\tau_{n,k}(x)$ is idempotent for all $x \in \L_m^+$ and $k \in \{0,\dots,n-1\}$. Hence, $\alg \L^+_m$ satisfies $\Phi_n$.
\end{proof}

\begin{proposition} \label{p:is-one-of-us}
    Let $\alg{M}$ be a nontrivial totally ordered unital commutative $\ell$-monoid, let $n \in \N \setminus \{0\}$, and suppose that, for every $x \in M$, the element $nx$ is integer.
    Then there is a divisor $k$ of $n$ such that $\alg{M} \cong \frac{1}{k}\Z$.
\end{proposition}

\begin{proof}
    First of all, we prove that $\alg{M}$ is cancellative.
    Let $x,y,z \in M$ be such that $x + z = y + z$.
    Using that $n \neq 0$, we get
    \[
    x + nz = (x + z) +(n-1)z = (y + z) +(n-1)z = y + nz.
    \]
    The element $nz$ is integer and hence invertible; thus, $x = y$.
    This proves that $\alg{M}$ is cancellative.

    We prove that the unique homomorphism $\varphi$ from $\alg{M}$ to $\R$ is injective.
    Let $x,y \in M$ be such that $\varphi(x) = \varphi(y)$.
    Since $nx$ and $ny$ are integer elements, there are $a, b \in \Z$ such that $nx = a$ and $ny = b$.
    We have
    \[
    a = \varphi(a) = \varphi(nx) = n\varphi(x) = n \varphi(y) = \varphi(ny) = \varphi(b) = b.
    \]
    Therefore, $nx = ny$.
    Since $\alg{M}$ is cancellative, it is torsion-free by \cref{l:torsion-free-monoid}.
    Therefore, using that $n \neq 0$, we obtain $x = y$.
    This proves that $\varphi$ is injective.

    For every element $x$ in the image $B$ of $\alg{M}$, $nx$ is an integer element of $\R$.
    Therefore, $\alg{B}$ is a subalgebra of $\frac{1}{n}\Z$.
    Therefore, by \cref{l:dense-or-discrete}, there is a divisor $k$ of $n$ such that $B = \frac{1}{k}\Z$ and hence $\alg{M}$ is isomorphic to $\frac{1}{k}\Z$.
\end{proof}

\begin{theorem} \label{t:precisely-the-divisors}
    For every $n \in \N \setminus \{0\}$, the subdirectly irreducible MV-monoids satisfying $\Phi_{n}$ are precisely the algebras isomorphic to $\alg{\L}_k^+$ for some divisor $k$ of $n$.
\end{theorem}

\begin{proof}
    By the right-to-left implication in \cref{l:lnSigman}, for every divisor $k$ of $n$, $\alg{\L}_k^+$ satisfies $\Phi_n$.

    Conversely, let $\alg{A}$ be a subdirectly irreducible MV-monoid satisfying $\Phi_{n}$.
    For every $x \in \Xi(\alg{A})$, there are $y_1, \dots, y_l \in \Gamma(\Xi(\alg{A}))$ and $m \in \Z$ such that $x = m + \sum_{i=1}^l y_l$ (\cite[Proposition~4.68]{Abbadini2021}); by \cref{p:totord-then-integer}, for every $i \in \{1, \dots, l\}$ the element $ny_i$ is integer; it follows that $nx = nm + \sum_{i=1}^l ny_l$ is an integer element.
    Therefore, by \cref{p:is-one-of-us}, there is a divisor $k$ of $n$ such that $\Xi(\alg{A}) \cong \frac{1}{k}\Z$.
    Therefore, $\alg{A} \cong \alg{\L}_k^+$.
\end{proof}

\begin{corollary}\label{thm:axiomsofln}
    For every $n \in \N \setminus \{0\}$, the variety $\VV (\alg{\L}^+_{n})$ is axiomatized by $\Phi_n$ relatively to the variety of MV-monoids.
\end{corollary}

From the previous theorem we immediately get the next corollary.

\begin{corollary}\label{thm:sigman}
    Let $\VV$ be a variety of MV-monoids and suppose that $\VV$ satisfies $\Phi_n$ for some $n \in \N \setminus \{0\}$. Then $\VV$ is a variety of positive MV-algebras.
\end{corollary}
\begin{proof}
    By \cref{thm:axiomsofln} we have $\VV \subseteq \VV(\alg{\L}_n^+)$, and the latter consists entirely of positive MV-algebras by \cref{c:finitely gen var of MVs}.
\end{proof}
We now start to deal with the sets of equations $\Sigma_I$ whose purpose is to set $\VV(\mathcal{K}_I)$ apart from every variety $\VV(\mathcal{K}_J)$ with $\lcm(I) = \lcm(J)$ but $I \neq J$ (for $I$ and $J$ divisor-closed finite sets).

\begin{definition}\label{d:sigman}
    Let $I$ be a divisor-closed finite set, and let $m$ be the maximum of $I$ (with the convention that $m = 0$ if $I = \emptyset$).
    We define $\Sigma_I$ as the set of equations given by the single equation
    \begin{equation}\label{eq:splittingmult-rename}
        (m+1)x \app mx
    \end{equation}
    together with the equations of the form
    \begin{equation}\label{eq:splittingdiv-rename}
         m((k-1)x)^k \app (kx) ^m
    \end{equation}
    for all $1 \leq k \leq m$ such that $k \notin I$.
\end{definition}
The equations $\Sigma_I$ are inspired by the ones used by Di Nola and Lettieri in \cite{DinolaLett} to axiomatize the varieties of MV-algebras.
Our equation \eqref{eq:splittingmult-rename} plays a similar role to the role played by the equation (1) in \cite[p.~466]{DinolaLett}; cf.\ our \cref{l:killerbigger} below with \cite[Theorems~2 and~3]{DinolaLett}.
The equations in \eqref{eq:splittingdiv-rename} are simplified versions of the equations (2) in \cite[p.~466]{DinolaLett}: we take the dual of their equation $(kx^{k-1})^{m+1} \app (m+1)x ^k$ and we observe that we can replace $m+1$ with $m$. The main difference in our setting is the fact that varieties of positive MV-algebras are generated by a finite set of finite subdirectly irreducible elements (\cref{t:VarofpMVs}).

\begin{lemma}\label{l:killerbigger}
    Let $m, n \in \N$ with $n > 0$. Then $\alg \L^+_n$ satisfies
    \[
    (m+1)x \app mx,
    \]
    if and only if $n \leq m$.
\end{lemma}
\begin{proof}
    If $n \leq m$, it is not difficult to see that $\alg \L^+_n$ satisfies $(m+1)x \app mx$. Indeed, for all $a \in \L^+_n \setminus \{0\}$ we have $ma = 1$ and so the evaluation of both the left-hand side and the right-hand side at $a$ is $1$, while for $a = 0$ they are both $0$.

    If $n > m$, the equation $(m+1)x \app mx$ fails in $\alg \L^+_n$, as witnessed by $\frac{1}{n} \in \L^+_n$.
\end{proof}
The following proposition is similar to {\cite[Theorem~8]{DinolaLett}}. In our case we use $m((k-1)x)^k \app (kx) ^m$ instead of $(kx^{k-1})^{m+1} \app (m+1)x ^k$ as characterizing equation. The former is the dual of the latter up to replacing $(m+1)$ with $m$. We give a proof of the proposition since both the hypothesis and the equation are slightly different from those of \cite[Theorem~8]{DinolaLett}.

\begin{proposition}\label{p:eq2}
    Let $m \in \N\setminus\{0\}$ and $1 \leq k \leq m$.
    For every $1 \leq n \leq m$, the algebra $\alg{\L}^+_n$ satisfies
    \begin{equation*}
        m((k-1)x)^k \app (kx) ^m
    \end{equation*}
    if and only if $k$ does not divide $n$.
\end{proposition}

\begin{proof}
    Fix $1 \leq n \leq m$.
    We investigate which elements $a \in \L^+_n$ satisfy the equation in the statement.
    Let $a \in \L^+_n$. We consider three cases.

    Case $a < \frac{1}{k}$. In this case, $ka < 1$ and $({k-1})a < \frac{k-1}{k}$. From $ka < 1$ we obtain $ka \leq \frac{n-1}{n} \leq \frac{m-1}{m}$, and so $(ka)^m = 0$.
    From $(k -1)a < \frac{k -1}{k}$ we obtain $((k-1)a)^k = 0 = m((k-1)a)^k$.
    Therefore, the equation $m((k-1)x)^k \app (kx) ^m$ holds for $a < \frac{1}{k}$ (for every $n$).

    Case $a > \frac{1}{k}$. In this case, $ka = 1$ and $(k-1)a > \frac{k - 1}{k}$. From $ka = 1$ we obtain $(ka)^m = 1$.
    From $(k-1)a > \frac{k-1}{k}$ we obtain $((k-1)a)^k > 0$, which implies $((k-1)a)^k \geq \frac{1}{n} \geq \frac{1}{m}$, and hence $m((k-1)a)^k = 1$. Therefore, the equation $m((k-1)x)^k \app (kx) ^m$ holds for $a > \frac{1}{k}$ (for every  $n$).

    Case $a = \frac{1}{k}$. In this case, $ka = 1$ and $(k-1)a = \frac{k-1}{k}$. From $ka = 1$ we obtain that $(ka)^m = 1$. From $(k-1)a = \frac{k-1}{k}$ we obtain $((k-1)a)^k = 0$, and so $m((k-1)a)^k = 0$.
    Thus, the equation $m((k-1)x)^k \app (kx) ^m$ fails for $a = \frac{1}{k}$. We have $\frac{1}{k} \in \L^+_n$ if and only if $k$ divides $n$.
    Therefore, the equation $m((k-1)x)^k \app (kx) ^m$ holds if and only if $k$ does not divide $n$.
\end{proof}

Using the previous proposition and the definition of divisor-closed finite sets we immediately get the next corollary.

\begin{corollary}\label{c:killersmaller}
    Let $I$ be a divisor-closed finite set. Let $m \coloneqq \max(I)$ if $I$ is nonempty and $m \coloneqq 0$ otherwise, and let $\Psi_I$ be the set of equations of the form \eqref{eq:splittingdiv-rename} in $\Sigma_I$.
    For every $1 \leq i \leq m$, $i \in I$ if and only if $\alg{\L}_i^+$ satisfies all equations in $\Psi_I$.
\end{corollary}

\begin{proof}
    Recall that, by \cref{d:sigman}, $\Psi_I$ consists of all equations of the form $m((k-1)x)^k \app (kx) ^m \in \Psi_I$ where $k \notin I$.

    We first prove $(\Rightarrow)$.
    Suppose $i \in I$.
    Consider an equation in $\Psi_I$; this will be of the form $m((k-1)x)^k \app (kx) ^m$ for some $k \notin I$.
    Since $I$ is closed under divisors, $k$ does not divide $i$.
    Therefore, by \cref{p:eq2}, $\alg{\L}_i^+$ satisfies $m((k-1)x)^k \app (kx) ^m$.
    This proves that $\alg{\L}_i^+$ satisfies all equations in $\Psi_I$.

    Then, we prove the contrapositive of $(\Leftarrow)$.
    Suppose $i \not\in I$ with $1 \leq i \leq m$. Then, by \cref{d:sigman}, $m((i-1)x)^i \app (ix) ^m \in \Psi_I$ and, by \cref{p:eq2}, $\alg{\L}_i^+$
    does not satisfy $m((i-1)x)^i \app (ix) ^m$.
\end{proof}
We are now ready to prove the main result of the section describing the axiomatizations of all the varieties of MV-algebras relatively to the variety of MV-monoids.

\begin{theorem}\label{thm:axioms}
    Let $I$ be a divisor-closed finite set; then $\vv V(\mathcal{K}_I)$ is axiomatized by $\Phi_{\mathrm{lcm}(I)} \cup \Sigma_I$ relatively to the variety of MV-monoids.
\end{theorem}

\begin{proof}
    It is enough to prove that the subdirectly irreducible algebras satisfying $\Phi_{\lcm(I)} \cup \Sigma_I$ are precisely the algebras isomorphic to some algebras in $\mathcal{K}_I$.

    For one direction, we note that any algebra in $\mathcal{K}_I$ satisfies $\Phi_{\lcm(I)} \cup \Sigma_I$ (see \cref{t:precisely-the-divisors} for $\Phi_{\lcm(I)}$, the right-to-left implication of \cref{l:killerbigger} and the left-to-right implication in \cref{c:killersmaller} for $\Sigma_I$).

    For the converse direction, let $\alg{A}$ be a subdirectly irreducible algebra satisfying $\Phi_{\lcm(I)} \cup \Sigma_I$.
    Then, by \cref{t:precisely-the-divisors}, there is $n \in \N\setminus \{0\}$ such that $\alg{A}\cong\alg{\L}_n^+$.
    By the left-to-right implication in \cref{l:killerbigger}, we have $n \leq \max(I)$, and then by the right-to-left implication in \cref{c:killersmaller} we deduce $n \in I$. Thus, $\alg{A}$ is isomorphic to an algebra in $\mathcal{K}_I$, namely $\alg{\L}_n^+$.
\end{proof}

Using \cref{thm:axioms} we produce some easy examples of axiomatizations of varieties of positive MV-algebras.

\begin{example}
    Let $I = \{1,2,3\}$ and $J = \{1,2,3,6\}$.
    By \cref{thm:axioms}, $\vv V(\mathcal{K}_{I})$ is axiomatized by the axioms of MV-monoids, $\Phi_6$ and $4x \app 3x$, while  $\vv V(\mathcal{K}_{J})$ is axiomatized by the axioms of MV-monoids, $\Phi_6$, $7x \app 6x$, $6(3x)^4 \app (4x)^6$ and $6(4x)^5 \app (5x)^6$. Note that the failure of $4x \app 3x$ in $\alg{\L}^+_{6}$ is witnessed by $\frac{1}{6}$.
\end{example}

We conclude the section by providing axiomatizations for all almost minimal varieties of MV-monoids. By \cref{thm:axiomsofln} we already produced an axiomatization of $\vv V(\alg \L_n^+)$ for each $n \in \N\setminus \{0\}$. Therefore, by \cref{t:char-almost-minimal}, we are only left to axiomatize $\vv V(\alg{C}_2^{\Delta})$ and $\vv V(\alg{C}_2^{\nabla})$.

\begin{theorem} \label{t:axiomVarC2}
	The varieties $\vv V(\alg{C}_2^{\Delta})$ and $\vv V(\alg{C}_2^{\nabla})$ are axiomatized, within the variety of MV-monoids, by $x \oplus x \app x$ and $x \odot x \app x$, respectively.
\end{theorem}

\begin{proof}
	We only prove the statement for $\vv V(\alg{C}_2^{\Delta})$; the other one is analogous.
		
	We claim that any subdirectly irreducible MV-monoid satisfying $x \oplus x \app x$ is isomorphic to $\alg{\L}_1^{+}$ or $\alg{C}_2^{\Delta}$.
	Let $\alg{A}$ be one such. By \cref{subdMVM}, $\alg{A}$ is nontrivial, totally ordered, and for all $x,y \in A$ either $x \oplus y = 1$ or $x \odot y = 0$. Thus, for all $x,y \in A$ we have
	\[
	x \lor y \leq x \oplus y \leq (x \lor y) \oplus (x \lor y) = x \lor y,
	\]
	and so $\alg{A}$ satisfies
	\[
	x \oplus y = x \lor y.
	\]
	Moreover, by \cref{c:Holder-MV}, there is a unique homomorphism $\varphi \colon \alg{A} \to [0,1]^+$.
	The image of $\varphi$ is $\{0,1\}$, since $\varphi(x) \oplus \varphi(x) = \varphi(x \oplus x) = \varphi(x)$ for all $x \in A$, and since every $y \in [0,1]^+$ satisfying $y \oplus y = y$ belongs to $\{0,1\}$.
	Furthermore, the only element of $A$ mapped to $1$ by $\varphi$ is $1$. Indeed, if $\varphi(x) = 1$, then
	\[
		\varphi(x) \odot \varphi(x) = 1 \odot 1 = 1 \neq 0 = \varphi(0),
	\]
	which implies $x \odot x \neq 0$, which implies $x \oplus x = 1$ and thus $x = x \oplus x = 1$.	Moreover, for every $a \in A \setminus \{1\}$, we have $a \odot a = 0$ because $a \oplus a = a \neq 1$. Using the fact that $\alg{A}$ is totally ordered, it follows that for all $a,b \in A \setminus \{1\}$ we have $a \odot b = 0$.
	
	To prove the claim, we suppose $\alg{A} \not\cong \alg{\L}_1^{+}$ and show $\alg{A} \cong \alg{C}_2^{\Delta}$.
	Let $a \in A \setminus \{0,1\}$. Let $\sim_l$ be the lattice-congruence obtained by identifying all elements in $\{x \in A \mid x \leq a\}$ and let $\sim_u$ be the lattice-congruence obtained by identifying all elements in $\{x \in A \mid a \leq x < 1\}$.
	Both $\sim_l$ and $\sim_u$ are MV-monoid congruences; indeed, they clearly respect $\oplus$ (because $\oplus = \lor$), and they respect $\odot$ because for all $a,b \in A \setminus \{1\}$ we have $a \odot b = 0$.
	Clearly, the meet of $\sim_l$ and $\sim_r$ is the trivial congruence.
	Then, since $\alg{A}$ is subdirectly irreducible, either $\sim_l$ or $\sim_r$ are the trivial congruence.
	Since $a \neq 0$, the congruence $\sim_l$ is not trivial, and thus $\sim_r$ is.
	Hence, $a$ is the maximum of $A \setminus \{0,1\}$.
	Since this holds for any $a \in A \setminus \{0,1\}$, $A$ must have three elements, say $0 < \e < 1$, and they satisfy $\e \oplus \e = \e$ and $\e \odot \e = 0$; thus, $\alg{A} \cong \alg{C}_2^{\Delta}$.
	
	This proves our claim that any subdirectly irreducible MV-monoid satisfying $x \oplus x \app x$ is isomorphic to $\alg{\L}_1^{+}$ or $\alg{C}_2^{\Delta}$. Then, since $\alg{C}_2^{\Delta}$ satisfies $x \oplus x \app x$ and since $\alg{\L}_1^{+} \in \II\SU(\alg{C}_2^{\Delta})$,  $x \oplus x \app x$ axiomatizes $\vv V(\alg{C}_2^{\Delta})$.
\end{proof}

\begin{remark}
	From the claim in the proof of \cref{t:axiomVarC2} and from its dual we immediately obtain that the equations $x \oplus x \app x$ and $x \odot x \app x$ axiomatize $\vv V(\alg{\L}_1^{+})$ within the variety of MV-monoids (which we already knew from \cref{thm:axioms,t:axiomVarC2}).
\end{remark}

To conclude, in the following table we summarize our axiomatizations of the almost minimal varieties of MV-monoids and the varieties of positive MV-algebras.
\begin{center}
    \begin{tabular}{|l|l|}
        \hline
        \textbf{Variety} & \textbf{Axiomatization (within MV-monoids)}\\
        \hline
        $\VV(\alg{C}_2^\Delta)$ & $x \oplus x \approx x$ \\
        \hline
        $\VV(\alg{C}_2^\nabla)$ & $x \odot x \approx x$\\
        \hline
        $\VV(\alg{L}_1^+)$ & $x \oplus x \approx x$ and $x \odot x \approx x$\\
        \hline
        $\VV(\alg{\L}_n^+)$ & $\tau_{n,k}(x) \oplus \tau_{n,k}(x) \approx \tau_{n,k}(x)$ \ \hfill(for $0 \leq k \leq n-1$)\\
        &
        $\tau_{n,k}(x) \odot \tau_{n,k}(x) \approx \tau_{n,k}(x)$ \ \hfill(for $0 \leq k \leq n-1$)\\
        \hline
        $\VV(\{\alg \L_n^+ \mid n \in I\})$ & \ \hfill (setting $l \coloneqq \lcm{I}$ and $m \coloneqq \max{I}$)\\
        ($I$ div.-closed fin.\ set)& $\tau_{l,k}(x) \oplus \tau_{l,k}(x) \approx \tau_{l,k}(x)$ \ \hfill (for $0 \leq k \leq l-1$)\\
        & $\tau_{l,k}(x) \odot \tau_{l,k}(x) \approx \tau_{l,k}(x)$ \ \hfill (for $0 \leq k \leq l-1$)\\
        & $(m + 1)x \approx mx$\\
        & $m((k-1)x)^k \app (kx) ^m$ \ \hfill (for $1 \leq k \leq m$ s.t.\ $k \notin I$)\\
        \hline
    \end{tabular}
\end{center}

\black

\section{Conclusions}
In the Introduction, we claimed that the main motivation for our investigation was to show that MV-monoids and positive MV-algebras can be studied (with some success) by employing techniques that are usually implemented for varieties of logic.
We believe we have substantiated our claim by proving several structural results as strong as in the context of MV-algebras, within a more relaxed framework.

In detail, we characterized the almost minimal varieties of MV-monoids (\cref{t:char-almost-minimal}) and we proved two versions of H\"older's theorem for unital commutative $\ell$-monoids (\cref{Holder,c:Holder}).

In the cancellative setting, we characterized the varieties of positive MV-algebras as precisely the varieties generated by finitely many reducts of finite nontrivial MV-chains (\cref{t:VarofpMVs}). We also proved that such reducts coincide with the subdirectly irreducible finite positive MV-algebras (\cref{{thm:fin-is-sub}}).
Furthermore, we provided an axiomatization for each variety of positive MV-algebras and each almost minimal variety of MV-monoids (\cref{thm:axioms,}).

Where can we go from here? We believe that there are at least a couple of paths worth exploring. First one could go on applying universal algebraic techniques
to those classes, in order to get a better general knowledge of their behavior. For instance, what are the structurally complete or primitive sub(quasi)varieties? Can one characterize to some degree the projective algebras or the splitting algebras?

Another possibility is to weaken further the axioms of MV-algebras while trying to maintain the possibility of finding a reasonable characterization of some of their most descriptive features (subdirectly irreducible algebras, minimal varieties, finitely generated varieties and so forth).
As an example, we may relax the lattice order to be a partial order.
In fact, the lattice order plays a peripheral role in various circumstances. For instance, the set $\Xi(\alg{A})$ (in the construction of the quasi-inverse $\Xi$ of $\Gamma$) consists of all good sequences in $\alg{A}$, a notion that does not involve the lattice order.
Moreover, only the function symbols $\oplus$, $\odot$, $0$ and $1$ are involved in our axiomatizations (within MV-monoids) of the varieties of positive MV-algebras and the almost minimal varieties of MV-monoids.
Dropping the requirement that the order is a lattice may clarify the role played by $\lor$ and $\land$, adding to the study of partially ordered monoids and \emph{commutative bimonoids} \cite{GalatosPrenosil2023}.

\section*{Acknowledgments}

We appreciate the referee's careful reading of the paper and their many suggestions, which have helped improve it.
We are also grateful to those who contributed to the organization of Algebra Week 2023 at the University of Siena, which provided the opportunity for this collaboration to begin.

\subsection*{Funding} The first author's research was funded by UK Research and Innovation (UKRI) under the UK government’s Horizon Europe funding guarantee (grant number EP/Y015029/1, Project ``DCPOS''). The ``Horizon Europe guarantee'' scheme provides funding to researchers and innovators who were unable to receive their Horizon Europe funding (in this case, a Marie Skłodowska-Curie Actions (MSCA) grant) while the UK was in the process of associating. The third author's research was partially funded by the Austrian Science Fund FWF P33878 and by the PRIMUS/24/SCI/008 of the Charles University Prague.

\bibliographystyle{plain}

\bibliography{Biblio}

\end{document}